\documentclass[11pt]{amsart}

\usepackage{amscd, amsmath, amssymb,verbatim,graphicx,chngcntr}
\usepackage{epsfig, color}
\usepackage{amsmath,amssymb}
\usepackage{amsthm}
\usepackage{hyperref}
\usepackage[margin=1.1in]{geometry}
\usepackage{esint}
\usepackage{stackrel}

\numberwithin{equation}{section}

\newtheorem{prop}{Proposition}
\newtheorem{lemma}[prop]{Lemma}

\newtheorem{thm}[prop]{Theorem}
\newtheorem{cor}[prop]{Corollary}

\numberwithin{prop}{section}

\theoremstyle{definition}
\newtheorem{defn}[prop]{Definition}

\newtheorem{rmk}[prop]{Remark}

\newcommand{\del}{\partial}

\newcommand{\brs}[1]{\left| #1 \right|}

\newcommand{\gD}{\Delta}
\newcommand{\gd}{\delta}

\newcommand{\gl}{\lambda}

\newcommand{\ga}{\alpha}
\newcommand{\gb}{\beta}
\renewcommand{\ge}{\epsilon}
\newcommand{\N}{\nabla}

\newcommand{\rv}{{\rm v}}

\newcommand{\Ric}{{\rm Ric}}
\newcommand{\Vol}{{\rm Vol}}
\newcommand{\diam}{{\rm diam}}

\newcommand{\Lip}{{\rm Lip}}

\newcommand{\cC}{\mathcal{C}}

\newcommand{\cH}{\mathcal{H}}

\newcommand{\cL}{\mathcal{L}}

\newcommand{\cN}{\mathcal{N}}

\newcommand{\cR}{\mathcal{R}}
\newcommand{\cS}{\mathcal{S}}

\newcommand{\cV}{\mathcal{V}}

\newcommand{\dR}{\mathbb{R}}

\renewcommand{\SS}{\mathcal S}
\newcommand{\til}[1]{\widetilde{#1}}

\renewcommand{\bar}[1]{\overline{#1}}

\newcommand{\hook}{\mathbin{\hbox{\vrule height2.4pt width4.5pt depth-2pt
\vrule height5pt width0.4pt depth-2pt}}}

\newcommand{\gL}{\Lambda}

\newcommand{\IP}[1]{\left<#1\right>}

\newcommand{\RC}{\mathcal {RC}}

\DeclareMathOperator{\Sym}{Sym}
\DeclareMathOperator{\Rc}{Ric}

\DeclareMathOperator{\Rm}{Rm}

\DeclareMathOperator{\supp}{supp}

\begin{document}

\title[Codimension four regularity of generalized Einstein structures]{Codimension four regularity of generalized Einstein structures}

\begin{abstract} We establish codimension $4$ regularity of noncollapsed sequences of metrics with bounds on natural generalizations of the Ricci tensor.   We obtain a priori $L^2$ curvature estimates on such spaces, with diffeomorphism finiteness results and rigidity theorems as corollaries.
\end{abstract}

\date{\today}

\author{Xin Fu}
\address{Xin Fu\\Rowland Hall\\
         University of California\\
         Irvine, CA 92617}
\email{\href{mailto:fux6@uci.edu}{fux6@uci.edu}}

\author{Aaron Naber}
\address{Aaron Naber\\Lunt Hall\\ Northwestern University\\ Evanston IL 60208}
\email{\href{mailto:anaber@math.northwestern.edu
}{anaber@math.northwestern.edu
}}
\author{Jeffrey Streets}
\address{Jeffrey Streets\\Rowland Hall\\
         University of California\\
         Irvine, CA 92617}
\email{\href{mailto:jstreets@uci.edu}{jstreets@uci.edu}}

\maketitle

\section{Introduction}

Given a smooth manifold $M$, a pair $(g, H)$ of a Riemannian metric and closed three-form is generalized Einstein if
\begin{align*}
\Rc^g - \tfrac{1}{4} H^2 =&\ 0,\\
d^*_g H =&\ 0,
\end{align*}
where $H^2(X,Y) = \IP{X \hook H, Y \hook H}_g$.  These equations arise in physical theories of supergravity and renormalization group flow (cf. for instance \cite{CSW, Rgflow}).  Mathematically, they arise as critical points of a natural extension of the Einstein-Hilbert action to generalized geometry \cite{GRFbook}, and describe canonical metrics in complex geometry \cite{IP}.  The generalized Einstein condition can be rephrased in terms of the Ricci curvature of a natural connection with torsion.  Specifically, we set
\begin{align*}
\N = D + \tfrac{1}{2} g^{-1} H,
\end{align*}
where $D$ denotes the Levi-Civita connection.  We refer to this as a \emph{Bismut connection} due to its natural appearance in complex geometry \cite{Bismut}.  It turns out that
\begin{align*}
\Rc^{\N} = \Rc^g - \tfrac{1}{4} H^2 + \tfrac{1}{2} d^*_g H,
\end{align*}
so that the generalized Einstein condition in this case can be rephrased as
\begin{align*}
\Rc^{\N} \equiv 0.
\end{align*}

Going further, various physical theories, such as the Hull-Strominger system \cite{Hull,Str} (cf. \cite{GF}) suggest systems of equations coupling the Einstein equation to differential forms of varying degrees.  With these in mind we define a broader class of equations.  Let $H = \bigoplus_{i=0}^n H_i$ now denote a linear combination of differential forms of different degrees.  We say that a pair $(g, H)$ is generalized Einstein if
\begin{align*}
\Rc^g - \tfrac{1}{4} H^2 =&\ 0\, ,\\
d H + Q_1(H) =&\ 0\, ,\\
d^*_g H + Q_2(H) =&\ 0\, ,
\end{align*}
where now $H^2(X,Y) := \IP{X \hook H, Y \hook H}_g$ with $X \hook H$ a sum of differential forms.  The $Q_i(H)$ are differential forms which are quadratic expressions in $H$ defined using the metric $g$, interior product, and wedge product.  To simplify notation, we introduce a tensor
\begin{align*}
\RC = \left( \Rc^g - \tfrac{1}{4} H^2, dH + Q_1(H), d^*_g H + Q_2(H) \right) \in \Sym^2(T^* M) \oplus \Lambda^* \oplus \Lambda^*.
\end{align*}
Thus the generalized Einstein condition is equivalent to
\begin{align*}
\RC \equiv 0.
\end{align*}
Furthermore, we say that a pair $(g, H)$ has bounded generalized Ricci curvature if there exists a constant $\gL$ such that
\begin{align}
\brs{\RC} = \brs{\Rc^g - \tfrac{1}{4} H^2} + \brs{ dH + Q_1(H)} + \brs{d^*_g H + Q_2(H)} \leq \gL.
\end{align}
As discussed above, in the case $H \in \Lambda^3, Q_1(H) = Q_2(H) = 0$, this is a bound for the Ricci curvature of the Bismut connection, but in general the generalized Einstein equation might not represent the vanishing of the Ricci curvature of a connection. 

The generalized Einstein equation is elliptic for the pair $(g, H)$.  Aiming towards understanding limiting behavior, construting new examples, and understanding rigidity phenomena, it is natural to develop a structure theory for singularity formation.  In this work, building upon the well-developed regularity theory for Ricci curvature \cite{A89,ChC1,ChC2,CCTi_eps_reg,CheegerNaber_Ricci,Codim4,JiNa,CJN}, we establish various sharp structural results for the generalized Ricci curvature.  The first main result is an analogue of the classical Einstein situation in \cite{Codim4}, that the singular set of a Gromov-Hausdorff limit of manifolds with bounded generalized Ricci curvature has codimension $4$.

\begin{thm} \label{t:codim4} Let $(M^n, d_j, p_j) \to (X, d, p)$ be a Gromov-Hausdorff limit of manifolds with $\brs{\RC_{M_j^n}} \leq n-1$, and $\Vol(B_1(p_j)) > \rv > 0$.  Then the singular set $\SS$ satisfies
\begin{align*}
\dim \SS \leq n-4.
\end{align*}
\end{thm}

The classic examples of Eguchi-Hanson \cite{EH} show that this result is sharp.  Furthermore, recent work of the third author and Ustinovksiy \cite{SU} gives an extension of the Gibbons-Hawking ansatz to generalized K\"ahler geometry, yielding families of metrics solving a soliton version of the generalized Einstein equation.  While not generalized Einstein, these metrics have $H$ nontrivial, and have bounded Bismut Ricci curvature even as the poles `collide,' forming a codimension $4$ orbifold singularity as in the classic Gibbons-Hawking/Eguchi-Hanson picture.  

Going further, we exhibit effective volume estimates for the singular set, which, combined with a detailed analysis of a suitable chosen neck region as in \cite{JiNa} leads to a sharp $L^2$ estimate for the Riemannian curvature, and sharp $L^p$ estimates for $H$ and $\N H$ as well.

\begin{thm}\label{t:L2curvature}
There exists $C=C(n,\rv)$ such that if $(M^n,g, H)$ satisfies $\brs{\RC_{M^n}} \le n-1$ and $\Vol(B_1(p))>\rv>0$, then we have
\begin{enumerate}\label{e:main_Rm_est2}
\item $\fint_{B_1(p)}\brs{\Rm}^2\leq C\, $
\item $\fint_{B_1(p)}\brs{H}^4\leq C\, $
\item $\fint_{B_1(p)}\brs{\nabla H}^2\leq C\, $
\end{enumerate}
\end{thm}

These analytic tools have a few immediate geometric and topological consequences.  First, following \cite{Codim4} we establish diffeomorphism finiteness of noncollapsed $4$-manifolds with bounded generalized Ricci curvature and diameter.  This result is sharp, in the sense that none of these conditions can be removed while maintaining finitely many diffeotypes.

\begin{cor}\label{t:main_finite_diffeo}
There exists $C=C(\rv,D)$ such that if  $M^4$ satisfies $|\RC_{M^4}|\leq 3$, $\Vol(B_1(p))>\rv>0$ 
and $\diam (M^n)\leq D$, then  $M^4$ can have one of at
 most $C$ diffeomorphism types.
\end{cor}

Using the sharp $L^2$ estimate we can also obtain a geometric rigidity result.  Note that compact generalized Bismut-Ricci flat structures exhibit an interesting rigidity phenomenon in four dimensions.  In particular, if $H$ is nonzero, then there is only one example: $M= S^3 \times S^1$, $g = g_{S^3} \oplus g_{S^1}$, with $H = dV_{g_{S^3}}$ (cf. \cite{GRFbook} Theorem 8.26).  Using our structure theory, we can extend this rigidity phenomena to the case of complete generalized Bismut-Ricci flat structures of maximal volume growth.

\begin{cor} \label{c:rigidity} Let $(M^4, g, H)$ be a complete manifold with $H \in \Lambda^3$ such that $\Rc^{\N}_{M^4} \equiv 0$ and $g$ has maximal volume growth.  Then $H \equiv 0$ and $g$ is Ricci-flat.
\end{cor}

\vspace{.3cm}

\noindent \textbf{Acknowledgements:} The third author was supported by DMS-1454854.  The authors thank Mario Garcia-Fernandez, Wenshuai Jiang and Yury Ustinovskiy for helpful conversations.

\vspace{.5cm}

\section{An $\epsilon$-regularity result and codimension $2$ regularity}

A classic result of Anderson says that a noncollapsed Riemannian unit ball with two sided Ricci bound will have uniform harmonic radius lower bound if it is also sufficient close to the Euclidean ball in Gromov-Hausdorff topology. 
The main purpose of this section is to prove a regularity result, Lemma \ref{regularity} below, for noncollapsed Rimannian manifolds with a two sided generalized Ricci curvature bound, which generalizes the result of Anderson \cite{Anderson_Einstein}. The main point is that bounded generalized Ricci assumption will give us an elliptic system for $g$ and $H$ in harmonic coordinates.  Furthermore, the form of the generalized Ricci curvature implies that if it is bounded, then the classical Ricci curvature is automatically bounded below.

To begin we recall some classical stratification theory for limiting spaces of noncollapsed Riemannian manifold with Ricci curvature bounded below.  Let $(M^n_i,g_i,p_i)\stackrel {d_{GH}}{\longrightarrow} (X,d,p)$ under the conditions
\begin{align*}
&\Ric_{M^n_i}\geq\ -(n-1)\, ,\\
&\Vol(B_1(p_i))>\ \rv>0\, .
\end{align*}

\noindent The classical stratification theory is based on the idea of separating out points of limiting metric space $X$ based the symmetries of tangent cones.  In the context of Ricci limit spaces this was first carried out by Cheeger and Colding \cite{ChC2}.  To introduce it let us recall the appropriate notion of symmetry involved.  Let $C(Y)$ denote the metric cone on the metric space $Y$, then we recall:

\begin{defn}
\label{d:ksymm}
The metric space $X$ is called {\it $k$-symmetric} if 
$X$ is isometric to $\dR^k\times C(Z)$ for some $Z$.  Furthermore, we say $X$ is \emph{$k$-symmetric at $x\in X$} if  there is an isometry of $X$ with $ \dR^k\times C(Z)$ which carries $x$ to a vertex of the cone $ \dR^k\times C(Z)$.
\end{defn}

The stratification of the singular set $\SS\subseteq X$ is built by considering the filtration
\begin{align}
\label{e:stratification}
\emptyset\subset \SS^0\subseteq\cdots\subseteq \SS^{n-1}:= \SS\subseteq X^n\, , 
\end{align}
where
\begin{align}
\label{e:sing_set}
\SS^k:= \{x\in X\, :\, \text{no tangent cone at } x \,\, {\rm is\,\,} (k+1)\text{-symmetric}\}\, .
\end{align}
The set $\SS^k\setminus \SS^{k-1}$ 
is called the \emph{$k$-th stratum} of the singular set.  A key result of \cite{ChC2} is the Hausdorff dimension bound
\begin{align}\label{e:sing_set_haus_est}
\dim\, \SS^k\leq k\, ,\;\;\;\;\; \text{for all}\,\, k\, .
\end{align}

In  \cite{ChC2},
\cite{Codim4}, by showing that $\SS^{n-1}\setminus \SS^{n-2}=\emptyset$, respectively $\SS^{n-1}\setminus
\SS^{n-4}=\emptyset$, the following sharper estimates were proved:
\begin{align}
\label{e:de}
&\dim \, \SS \leq n-2 ,\, \text{ if \,\,}\Ric_{M^n_i}\geq -(n-1)\, . \\
&\dim \, \SS \leq n-4, \,  \text{ if \,\,}|\Ric_{M^n_i}|\leq (n-1)\, .
\end{align}

We next turn to understanding the regular set of $X$.  A natural starting place is the harmonic radius of a point:

\begin{defn}\label{d:harmonic_radius}
 For $x\in X$, we define the \emph{harmonic radius} $r_h(x)$ so that $r_h(x)=0$ if no neighborhood of $x$ is a Riemannian manifold. 
 Otherwise, we define $r_h(x)$ to be the largest $r>0$ such that there exists a mapping $\Phi:B_r(0^n)\to X$ such that:
\begin{enumerate}
\item $\Phi(0)=x$ with $\Phi$ is a diffeomorphism onto its image.
\vskip1mm

\item $\Delta_g x^\ell = 0$, where $x^\ell$ are the coordinate functions and $\Delta_g$
 is the Laplace Beltrami operator.
\vskip1mm

\item If $g_{ij}=\Phi^*g$ is the pullback metric, then
\begin{align*}
||g_{ij}-\delta_{ij}||_{C^0(B_r(0^n))}+r||\partial_k g_{ij}||_{C^0(B_r(0^n))} \leq 10^{-3}\, .
\end{align*}
\end{enumerate}
\end{defn}

\begin{defn}\label{1} Given $x \in X$ we say the \emph{regularity scale} $\tilde r_x$ is
\begin{align*}
\tilde r_x \equiv \max_{0 < r \leq 1} \left\{ \sup_{B_r(x)} \brs{\Rm} \leq r^{-2} \right\}.
\end{align*}
If $x \in \SS$ then $\tilde r_x \equiv 0$.
\end{defn}

With the above preliminaries we can now state the main result of this section:
 
\begin{lemma}\label{regularity} Given $n \in \mathbb N$, $\rv > 0$ there exists $\ge(n, \rv) > 0, r(n,\rv) > 0$ so that if $(M^n, g, H)$ satisfies
\begin{enumerate}
\item $\brs{\RC} \leq \ge$,
\item $\Vol(B_1(p)) > \rv$,
\item $d_{GH}(B_2(p), B_2(0)) \leq \ge$, \end{enumerate}
where $B_2(0)$ is Euclidean ball, then $r_h(p) \geq r$.  If $(M^n, g, H)$ is generalized Einstein, then $\tilde r_p \geq r$.
\end{lemma}

\begin{rmk} The arguments below rely on the blowup/contradiction/compactness argument.  In the proofs, for notational convenience, we state the natural scaling law for the case $H \in \Lambda^3$, where $\til{g} = \gl g, \til{H} = \gl H$.  If more generally $H = \sum_{i=1}^n H_i, H_i \in \Lambda^i$, the natural scaling law which preserves the generalized Ricci tensor is
\begin{align*}
\til{g} = \gl g, \qquad \til{H} = \sum_{i=1}^n \gl^{\frac{i-1}{2}} H_i.
\end{align*}

\end{rmk}

In the following lemma, we prove the regularity estimates for structures with bounded generalized Ricci curvature in harmonic coordinates.

\begin{lemma}\label{H-control-regular} Given $d, \gL > 0$, if $(M^n,g,H,x)$ satisfies
\begin{enumerate}
    \item $\brs{\RC} + \brs{H}^2 <\gL$,
    \item $r_h(x)>d$,
\end{enumerate}
then in harmonic coordinates around $y\in B_{2^{-1}d}(x)$ one has for any $0 < \ga < 1$ a constant $C = C(n,d,\gL,\ga)$ so that 
\begin{align*}
    \brs{\brs{H}}_{C^{\alpha}} + \brs{\brs{g_{ij} - \gd_{ij}}}_{C^{1,\alpha}} + \brs{\brs{\gD g_{ij}}}_{L^{\infty}} \leq C.
\end{align*}
\begin{proof} Since $\brs{H}^2+\brs{\RC}<\gL$, then $|dH|+|\Rc^g|+|d_g^*H|<C(\gL,n,d)$.  In harmonic coordinates, it follows by \cite{Anderson_Einstein} that
$$\Rc^g_{ij}=\Delta g_{ij}+Q(g,\partial g)_{ij}$$
Hence this gives the bound $||\gD g_{ij}||_{L^{\infty}}$. At last using $W^{2,p}$ theory of elliptic systems and the Sobolev inequality, we have $||g_{ij}-\delta_{ij}||_{C^{1,\alpha}}\leq C(\alpha,\gL,d,n)$ for any $\alpha$.  We can express locally $H = (d + d^*) B$.  By the generalized Ricci curvature bound and the assumed estimate for $H$, we obtain a uniform estimate for $||\gD_d B||_{L^{\infty}}$.  By $W^{2,p}$ estimates and Sobolev inequality we obtain a uniform $C^{1,\alpha}$ estimate for $B$, and thus the required $C^{\alpha}$ estimate for $H$.
\end{proof}
\end{lemma}

Using the above regularity we can show that the $C^0$ norm of $H$ has a priori control of a certain kind inside the harmonic radius.  The form of the quantity below allows $H$ to blow up at the scale invariant rate at the boundary of the given harmonic coordinates.

\begin{lemma}\label{H-control} Given $d, \gL > 0$ there exists $C > 0$ so that if $(M^n,g,H,x)$ satisfies
\begin{enumerate}
    \item $\brs{\RC}<\gL$,
    \item $r_h(x)>d$,
\end{enumerate}
then
\begin{align*}
    \sup_{p \in B_d(x)} \brs{H} (d - d(p,x)) < \infty.
\end{align*}
\begin{proof} If the claim were false, we can choose a sequence of points $x_i \in B_d(x)$ such that
\begin{align*}
    \brs{H}(x_i) (d - d(p,x_i)) \to \infty.
\end{align*} 
By a standard point-picking argument, we can ensure that this sequence of points also satisfies
\begin{align*}
    \sup_{B_{ (\brs{H}(x_i))^{-1}}(x_i)} \brs{H} \leq 2 \brs{H}(x_i).
\end{align*}
Now let $\gl_i := \brs{H}^2(x_i)$, and consider the blowup sequence $(\til{g}_i, \til{H}_i) = (\gl_i g, \gl_i H)$.  By construction of the sequence $x_i$ these are defined on a $\til{g}$-ball of radius $1$, and furthermore
\begin{align*}
    1 = \brs{\til{H_i}}_{\til{g}_i}(x_i) \leq \sup_{B_1(x_i,\til{g}_i)} \brs{\til{H}}_{\til{g_i}} \leq 2.
\end{align*}
Furthermore, since the original harmonic coordinates for $g$ are defined on $B_d(x)$, there are harmonic coordinates for $\til{g}_i$ on a ball of radius $1$, with the $C^1$ norm going to zero as $i \to \infty$.  By Lemma \ref{H-control-regular} we obtain a uniform $C^{1,\alpha}$ estimate for $\til{g}_i$ and a uniform $C^{\alpha}$ estimate for $\til{H}_i$.  By Arzela-Ascoli we obtain a limiting structure which is weakly generalized Ricci flat, and $\til{H}_{\infty}(x_{\infty}) = 1$.  However, as the $C^1$ norm of $\til{g}_i$ goes to zero, it follows that $g_{\infty}$ is flat.  Since the limit is weakly generalized Ricci flat, we can trace this to yield $R^{g_{\infty}} = \tfrac{1}{4} \brs{H_{\infty}}^2$, and it follows that $H_{\infty} \equiv 0$, giving a contradiction.
\end{proof}
\end{lemma}

With the above preparations, we can prove the main Lemma \ref{regularity} of this section.
\begin{proof}[Proof of Lemma \ref{regularity}]Under the same assumption, it suffices to prove a stronger estimate $$d(x,\partial B_2(0))^{-1}r_h(x)\geq r$$ for all $x\in B_2(0).$  We argue by contradiction, fix $\rv > 0$ and suppose there exists a sequence $\ge_i \to 0$ and $(M^n_i, g_i, H_i)$, with $p_i \in M_i$ such that
\begin{enumerate}
\item $\brs{\RC_{M_i}} \leq \ge_i$,
\item $\Vol(B_1(p_i)) > \rv$,
\item $d_{GH}(B_2(p_i), B_2(0)) \leq \ge_i$,
\item $d(x_i,\partial B_2(p_i))^{-1}r_h({x_i})\leq \frac{1}{i}$,
\end{enumerate}
where the minimum of the function $d(x,\partial B_2(p_i))^{-1}r_h(x)$ is achieved at $x_i$. 
Let
\begin{align*}
d_i:=d(x_i,\partial B_2(p_i)), \qquad r_i:=r_h(x_i).
\end{align*}
By the choice of $x_i$, for all points $y\in B_{3d_i/4}(x_i)$ we have $r_h(y)\geq \frac{r_i}{4}$. Consider the rescaled sequence $\{(B_{3d_i/4}(x_i), x_i, r_i^{-2} g_i, r_i^{-2} H_i)\}$. By Lemmas \ref{H-control-regular} and \ref{H-control}, we obtain that 
 $\{(B_{d_i/2}(x_i), x_i, r_i^{-2} g_i, r_i^{-2} H_i)\}$ converges to  a generalized Einstein manifold $Y$ in
 $C^{1,\alpha}$ topology for any $\alpha>0$. By the fact that $\Ric_Y\geq 0$, Colding's volume convergence theorem \cite{Co1}, and volume comparison, $Y$ is indeed the Euclidean space.  Since the harmonic radius is continuous in $C^{1,\alpha}$ topology, this  contradicts the fact that $r_h(x_i)=1$ along the blow up sequence.
\end{proof}


A corollary of Lemma \ref{regularity} is that the singular set $\SS$ of the limiting space of  a sequence of non collapsed Riemannian manifolds with two sided generalized Ricci bound is closed.
\begin{cor} Given a sequence of $(M^n_{i}, g_{i}, H_{i})$ of Riemannian manifolds satisfying
\begin{enumerate}
\item $\brs{\RC_{M_i}} \leq (n-1)$,
\item $\Vol(B_1(p_i)) > \rv$,
\end{enumerate}
then $(M^n_{i}, g_{i}, p_{i})$ converges subsequentially in the $C^{1,\alpha}$  topology outside the closed set $\cS$ to a $C^{1,\alpha}$ manifold $(X,d,p)$ for any $0 < \alpha<1$. 
\begin{proof}
Since bounded generalized Ricci tensor implies a lower bound of the Ricci tensor, we can directly apply Cheeger-Colding theory to extract a limit $(X,d,p)$ of $(M^n_{i},g_{i},p_{i})$. We only need to prove that the singular set $\cS$ is closed. Now we fix a smooth point $q\in X\setminus \cS$, for any fixed constant $\epsilon'$ small,  we have a sequence of metric ball $B_{\epsilon'}(q_{i},g_{i})$ contained in $M_i$ converging to $B_{\epsilon'}(q,d)$ 
in the Gromov-Hausdorff sense.  Now since a tangent cone at $q$ is $\mathbb R^n$,
 there is a large constant $\delta(\epsilon)$, which depends on the constant $\epsilon$ in Lemma \ref{regularity}, such that 
$$d_{GH}(B_{2}(q,\delta(\ge)^{-1}d), B_2(0)) \leq \frac{\ge}{2},
$$
where $B_2(0)$ is the Euclidean ball.  Hence we can find a sequence of balls $B_{2}(q_i,\delta(\ge)^{-1}g_i)$ contained in $M_i$ such that
$$d_{GH}(B_{2}(q_i,\delta(\ge)^{-1}g_i), B_2(0)) \leq \ge
$$
when $i$ is sufficiently large. By Lemma \ref{regularity}, it follows that there is a neighborhood $U$ of $q$ contained in $X\setminus \mathcal S$ and outside $\cS$, we have the $C^{1,\alpha}$ convergence. This completes the proof.
\end{proof}
\end{cor}

\vspace{.5cm}

\section{Codimension $4$ regularity}

In this section we prove Theorem \ref{t:codim4}, establishing that the singular sets have codimension $4$.  In the two subsections below we rule out codimension $2$ and $3$ singularities, respectively.

\vspace{.3cm}
\subsection{Nonexistence of codimension $2$ singularities}

In this subsection we prove that the singular set has codimension at least $3$.  The proof relies on analyzing the possible tangent cones, and there are two rigidity results for generalized Ricci-flat spaces key to this analysis.  First we show in Proposition \ref{cone} below that generalized Ricci-flat cones have vanishing $H$ and are thus Ricci-flat.  Furthermore, we show in Proposition \ref{p:2dimrigidity} that the only two-dimensional examples are either flat space with vanishing $H$ or the round two-sphere with $H$ a multiple of the volume form.  Given this, the proof follows ideas from the classical Ricci case, relying principally on the Slicing Theorem of Cheeger-Naber (\cite{Codim4} Theorem 1.23) to rule out certain codimension $2$ tangent cones.

\begin{prop}\label{cone} Let $(M^{n}, g, H)$ be a smooth manifold which is isometric to an open set in a generalized Einstein cone.  Then $\Rc^g \equiv 0, H \equiv 0$.
\begin{proof} By assumption we have $M \subset \mathbb R_+ \times \Sigma$, and on $M$ we can express
\begin{align*}
g = dt^2 + t^2 g_{\Sigma}.
\end{align*}
We know that the sectional curvatures of the mixed planes are zero, so the Ricci curvature of the $\del_t$ direction is zero, thus
\begin{align*}
0 = \Rc^g(\del_t, \del_t) = \tfrac{1}{4} H^2(\del_t, \del_t) = \tfrac{1}{4} \brs{\del_t \hook H}^2.
\end{align*}
Thus $\del_t \hook H = 0$.  Note that this also implies that for the quadratic polynomial $Q_1(H)$ we also have $i_{\del_t} Q_1(H) = 0$.  We thus obtain from the Cartan formula that
\begin{align*}
    L_{\del_t} H =&\ d i_{\del_t} H + i_{\del_t} d H\\
    =&\ i_{\del_t} Q_1(H)\\
    =&\ 0.
\end{align*}
As the Ricci tensor of a cone is expressed $\Rc^g = \Rc^{g_{\Sigma}} - (n-2) g_{\Sigma}$, it follows easily that $L_{\del_t} \Rc^g = 0$.  Furthermore, from the expression for $g$ we easily obtain $L_{\del_t} g =\ 2 t g_{\Sigma}$.  Putting these facts together, it follows from the generalized Ricci-flat equation then that
\begin{align*}
    0 =&\ \left( L_{\del_t} \RC \right)(X,X) = \left( L_{\del_t}H^2 \right)(X,X) \sim t \IP{X \hook H}^2_{g_{\Sigma}}.
\end{align*}
It follows that $X \hook H \equiv 0$.  Recalling that already  $\del_t\hook H=0$, it follows that $H \equiv 0$, and then $\Rc^g \equiv 0$.
\end{proof}
\end{prop}

\begin{prop} \label{p:2dimrigidity} Let $(\Sigma^2, g, H)$ be a smooth generalized Einstein manifold.  Then there exists $\gl \in \mathbb R$ such that $H = \gl dV_{g}$.  Furthermore, either $\gl = 0$ and $g$ is flat, or $\gl \neq 0$ and $g$ is homothetic to $g_{S^2}$.
\begin{proof} Since $\Sigma$ has dimension $2$, it follows that $H$ can only have components of degree $1$ or $2$.  Let us express $H = \phi + \psi$ where $\phi \in \Lambda^1, \psi \in \Lambda^2$.  One has $\phi^2 = \phi \otimes \phi$, and $\psi^2 = \tfrac{1}{2} \brs{\psi}^2 g$, and therefore the first equation of the generalized Einstein system yields
\begin{align} \label{p:2dim10}
0 =&\ \Rc^g - \tfrac{1}{4} H^2 = \tfrac{1}{2} R^g g - \tfrac{1}{4} \phi \otimes \phi - \tfrac{1}{8} \brs{\psi}^2 g
\end{align}
It follows that the tensor $\phi \otimes \phi$ must be a scalar multiple of $g$, which forces $\phi \equiv 0$.  It follows then that the third equation of the generalized Einstein system implies $d^*_g \psi = 0$, which implies $\psi = \gl dV_g$ for some constant $\gl$.  If $\gl = 0$, equation (\ref{p:2dim10}) implies $R^g \equiv 0$, and $g$ is flat.  If $\gl \neq 0$, then equation (\ref{p:2dim10}) implies $R^g$ is a positive constant, and then $g$ is homothetic to $g_{S^2}$.
\end{proof}

\end{prop}

Now we closely follow the argument of \cite{Codim4}.  Our first goal is to show an $\epsilon$-splitting result, Lemma \ref{h_splitting}, which converts Gromov-Hausdorff line-splitting to the existence of splitting functions with strong analytic estimates.  We first recall the definition of $\ge$-splitting maps.
\begin{defn}\label{d:intro}
An \emph{$\epsilon$-splitting map}  $u=(u^1,\ldots,u^{k}):B_r(p)\to\mathbb R^{k}$ is a harmonic map 
such that:

\begin{enumerate}

\item $|\nabla u|\leq 1+\epsilon$,
\vskip2mm

\item $\fint_{B_r(p)} |\langle \nabla u^\alpha,\nabla u^\beta\rangle-\delta^{\alpha\beta}|^2<\epsilon^2$,
\vskip2mm

\item $r^2\fint_{B_r(p)} |\nabla^2u^\alpha|^2<\epsilon^2$.
\end{enumerate}
\end{defn}
\noindent

The equivalence of harmonic splittings and geometric splittings follows from \cite{ChC1}.  The gradient estimate (1) follows from \cite{Codim4}:

\begin{lemma}[\cite{ChC1}]  \label{h_splitting}

For every $\epsilon,R>0$ there exists $\delta =\delta(n,\epsilon,R)>0$ such that if 
$\Ric^g_{M^n}\geq -(n-1)\delta$ then:
\begin{enumerate}
\item
 If $u:B_{2R}(p)\to\mathbb R^k$ is a $\delta$-splitting map, then there exists a 
map $f:B_{R}(p)\to u^{-1}(0)$ such that 
$$
(u,f):B_{R}(p)\to \mathbb R^k\times u^{-1}(0)\, ,
$$
is an $\epsilon$-Gromov-Hausdorff map, where $u^{-1}(0)$ is given the induced metric.

\vskip1mm

\item  If
\begin{align*}
d_{GH}(B_{\delta^{-1}}(p),B_{\delta^{-1}}(0))<\delta,
\end{align*}
where $0\in \mathbb R^k\times Y$, then there exists an $\epsilon$-splitting map $u:B_{R}(p)\to \mathbb R^k$.  
\end{enumerate}
\end{lemma}

For convenience we also record the Transformation Theorem (\cite{Codim4} Theorem 1.23)

\begin{thm}\label{transformation}
For each $\epsilon > 0$, there  exists $\delta(n,\ge)$ such that if $M^n$ satisfies $\Ric_{M^n}\geq-(n-1)\delta$ and $u: B_2(p)\rightarrow \mathbb R^{n-2}$ is a harmonic  $\delta$-splitting map, then there exists a subset $G_{\ge}\subset B_1(0^{n-2})$ that satisfies the following:
\begin{enumerate}
    \item $\Vol(G_{\ge}) > \Vol(B_1(0^{n-2})) - \ge$.
    \item If $s\in G_{\ge}$, then $u^{-1}(s)$ is nonempty
    \item For each $x\in u^{-1}(G_{\ge})$ and $r\le1$, there exists a lower triangular matrix $A(x,r)\in GL(n-2)$ with positive diagonal entries such that $A\circ u: B_r(x)\rightarrow \mathbb R^{n-2}$ is an $\ge$-splitting map
    \end{enumerate}
\end{thm}
We remark that the Transformation Theorem above only requires Ricci lower bound, which is satisfied if we assume there is a lower bound for the generalized Ricci tensor.  Putting together the above tools of Lemma \ref{regularity}, Proposition \ref{cone}, and Theorem \ref{transformation}, we prove the main technical result of this section following \cite{Codim4}.  In the statement below $S^1_{\beta}$ denotes the circle of circumference $\gb \leq 2 \pi$.

\begin{thm}\label{codim2}
Let $(M^n_i,g_i,H_i,p_i)$ be a sequence of Riemannian manifolds satisfying $|\RC_{M^n_i}|\to 0, \Vol(B_1(p_i))>\rv>0$, and such that 
$$(M_i,g_i,p_i)\stackrel{GH}{\longrightarrow}(\mathbb R^{n-2}\times C(S^1_{\beta}),d,p).$$
Then $\beta=2\pi$ and $\mathbb R^{n-2}\times C( S^1_{\beta})=\mathbb R^n$.
\end{thm}
\begin{proof}
We will prove the result by contradiction. Assuming the statement is false,  there exists a sequence $(M^n_i,g_i,H_i,p_i)$ of 
Riemannian manifolds satisfying $|\RC_{M^n_i}|\rightarrow 0 $, $\Vol(B_1(p_j))>\rv>0$, and such that
\begin{align*}
(M_i^n,d_i,p_i)\stackrel{GH} {\longrightarrow} (\mathbb R^{n-2}\times C(S^1_\beta),d,p) ,
\end{align*}
with $\beta<2\pi$ and $p$ a vertex.  

We first observe that by the noncollapsing assumption we have $\beta\geq \beta_0(n,\rv)$.  Now, by Lemma \ref{h_splitting}, there exist $\delta_i$-splitting maps 
$u_i:B_2(p_i)\rightarrow\mathbb R^{n-2}$ with $\delta_i\rightarrow 0$.  Fix some sequence $\epsilon_i\to 0$ which is 
tending to zero so slowly compared to $\delta_i$,  that Theorem \ref{transformation} holds
for $u_i:B_2(0)\to \mathbb R^{n-2}$ with $\epsilon_i$.  Let $G_{\epsilon_i}\subset B_1(0^{n-2})$ be
 the corresponding good values of $u_i$, and let $s_i\in G_{\epsilon_i}\cap B_{10^{-1}}(0^{n-2})$
 be fixed regular values.

Observe that $\mathbb R^{n-2}\times C(S^1_\beta)$ is smooth outside of the singular set 
$\mathcal S= \mathbb R^{n-2}\times\{0\}\subset \mathbb R^{n-2}\times C(S^1_\beta)$.  In particular on
$\mathbb R^{n-2}\times C(S^1_\beta)$ we have $r_h(x)\approx c d(x,\mathcal S)$ where $c$ depends on the link of the metric cone.  By Lemma \ref{regularity}, it follows
that the convergence of $M^n_i$ is in $C^{1,\alpha}\cap W^{2,q}$ away from $\mathcal S$, for every 
$\alpha<1$ and $q<\infty$.  Let $f_i:B_{\epsilon^{-1}_i}(p)\to B_{\epsilon^{-1}_i}(p_i)$ be the 
$\epsilon_i$-Gromov Hausdorff maps, and let us denote $\mathcal S_j\equiv f_j(\mathcal S)\subset M^n_j$.  
Then by the previous statements, for every $\tau>0$, all $j$ sufficiently large, and 
$x\in B_1(p_j)\setminus T_\tau(\mathcal S_j)$, we have $r_h(x)\geq \frac{\tau}{2}$.

Consider again the submanifold $u^{-1}_{j}(s_j)\cap B_1(p_j)$.  Define the scale
\begin{align*}
r_j= \min\{r_h(x):x\in u^{-1}_{j}(s_j)\cap B_1(p_j)\}\, .
\end{align*}
By the discussion of the previous paragraph, this minimum is actually obtained at some 
$x_j\in u^{-1}_{j}(s_j)\cap B_1(p_j)$, with $x_j\to \mathcal S_j\cap B_{10^{-1}}(p_j)$. Moreover, since $S^1_\beta$,  the cross-section of
the cone factor,  satisfies $0<\beta<2\pi$, it follows  that $r_j\to 0$.  According to Theorem \ref{transformation},
there exists a lower triangular matrix $A_j\in GL(n-2)$ with positive diagonal entries,
such that 
$v_j\equiv A_j\circ \big(u_j-s_j\big):B_{r_j}(x_j)\to \mathbb R^{n-2}$ is an $\epsilon_j$-splitting map.
Note that we have normalized so that each of our regular values is the zero level set.

Now let us consider the sequence $(M^n_j,r_j^{-1}d_j,x_j)$.  After passing to a subsequence if necessary, 
which we will continue to denote by $  (M^n_j,r_j^{-1}d_j,x_j)$, we have 
\begin{align*}
(M^n_j,r_j^{-1}d_j,x_j)\stackrel{d_{GH}}{\longrightarrow} (X,d_X,x),
\end{align*}
in the pointed Gromov-Hausdorff sense, where $X$ splits off $\mathbb R^{n-2}$ isometrically.  Note that, by our noncollapsing 
assumption, we have $\Vol(B_1(x_j))>c(n,\rv)>0$, and hence, in the rescaled spaces, we have 
$\Vol(B_r(x_j))>c(n,\rv) r^n$ for all $r\leq R_j\to \infty$.  In particular, $X$ has Euclidean volume
growth at $\infty$, i.e. 
$\Vol(B_r(x'))>c(n,\rv)\,r^n$ for all $r>0$.  

After possibly passing to another subsequence, we can limit the functions $v_j$ to a function $v:X\to \mathbb R^{n-2}$. 
 Note that by our normalization, we have
$v_j:B_{2}(x_j)\to \mathbb R^{n-2}$ are $\epsilon_j$-splittings, and that by the proof of the Theorem \ref{transformation} (cf. Cheeger-Naber \cite{Codim4} pages 1118-1121), we have for each $R>2$ that the maps $v_j:B_{R}(x_j)\to \mathbb R^{n-2}$ 
are $C(n,R)\epsilon_j$-splittings.  In particular, we can conclude that 
\begin{align*}
X= \mathbb R^{n-2}\times S,
\end{align*}
where $v:X\to\mathbb R^{n-2}$ is the projection map and  $S= v^{-1}(0)$.

Now by construction, in the rescaled spaces we have for any $y\in v^{-1}_j(0)$ that $r_h(y)\geq 1$.  
Therefore, the limit $X$ is $C^{1,\alpha}\cap W^{2,q}$ in a neighborhood of $v^{-1}(0)$, and hence 
$S= v^{-1}(0)$ is a nonsingular surface.  Thus, since $X=\mathbb R^{n-2}\times S$ it follows  that $X$ is
 at least a $C^{1,\alpha}\cap W^{2,q}$ manifold with $r_h\geq 1$.  Since the generalized  Ricci curvature is uniformly 
bounded, in fact tending to zero, we have by Lemma \ref{regularity} that the 
convergence $(M^n_j,r_j^{-1}d_j,x_j)\to (X, d_X,x)$ is in 
$C^{1,\alpha}\cap W^{2,q}$.  Because the convergence is in $C^{1,\alpha}\cap W^{2,q}$ we have that 
$r_h$ behaves continuously in the limit (cf. \cite{Anderson_Einstein}).  In particular, we have  $r_h(x'_j)\to r_h(x')$ and so, $r_h(x')=1$. 

On the other hand, since $|\RC_{M^n_j}|\to 0$ and $X$ is $C^{1,\alpha}\cap W^{2,q}$ it follows that $(X, d)$ is a 
smooth generalized Ricci-flat metric manifold, hence $S$ is also generalized Ricci flat.  By  Proposition \ref{p:2dimrigidity}, $S$ is actually Ricci flat and hence flat.  Since we have already shown that
$X$ has Euclidean volume growth, this implies that $X=\mathbb R^n$.
However, we have also already concluded that $r_h(x')=1$, which gives us our desired contradiction.
\end{proof}

Using this result, we obtain that a noncollapsed limit space is smooth away from 
a set of codimension $3$. We will use this in the next subsection to show $(n-3)$-symmetric splittings cannot arise as limits.

\begin{cor}\label{c:codim3}
Let $(M^n_i,g_i,H_i,p_i)$ denote a sequence of Riemannian manifolds satisfying $|\RC_{M^n_i}|\leq n-1$, $\Vol(B_1(p_j))>\rv>0$ and such that
\begin{align*}
(M_j^n,d_j,p_j) \to (X,d,p)\, .
\end{align*}
Then there exists a subset, $\mathcal S\subset X$, with $\dim\mathcal S\leq n-3$, such that for each $x\in X\setminus \mathcal S$,
 we have $r_h(x)>0$.  In particular, $X\setminus \mathcal S$ is a $C^{1,\alpha}$ Riemannian manifold.
\end{cor}

\begin{proof}
Recall the standard stratification of $X$.  In particular, if we consider the subset $\mathcal S^{n-3}\subset X$ 
we have that $\dim \mathcal S^{n-3}\leq n-3$, and for every point $x\not\in \mathcal S^{n-3}$ there exists {\it some} 
tangent cone at $x$ which is isometric to $\mathbb R^{n-2}\times C(S^1_\beta)$.  That is, there exists $r_a\to 0$ such that
\begin{align*}
(X,r_a^{-1}d,x)\to \mathbb R^{n-2}\times C(S^1_\beta)\, .
\end{align*}
However by Theorem \ref{codim2} we then have  $\beta= 2\pi$, which is to say that 
\begin{align*}
(X,r_a^{-1}d,x)\to \mathbb R^n.
\end{align*}
Thus, for $a\in \mathbb N$ sufficiently large, we can apply Lemma \ref{regularity} to balls in the rescaled limiting sequence,
to see that a neighborhood of $x$ is a $C^{1,\alpha}$ Riemannian manifold, which proves the corollary.
\end{proof}

\subsection{Nonexistence of Codimension $3$ singularities}\label{ss:codim3}

In this section we we prove that 
$(n-3)$-symmetric metric spaces cannot arise as limits of manifolds with bounded generalized Ricci curvature, finishing the proof of codimension $4$ regularity.  Specifically, we prove the following:

\begin{thm}\label{t:codim3}
Let $(M^n_i,g_i,H_i,p_i)$ denote a sequence of Riemannian manifolds satisfying $|\RC_{M^n_i}|\to 0$, $\Vol(B_1(p_i))>\rv>0$ and such that
\begin{align}
(M_i^n,d_i,p_i) \to \mathbb R^{n-3}\times C(Y)\, ,
\end{align}
in the pointed Gromov-Hausdorff sense, where $Y$ is some compact metric space.  
Then $Y$ is isometric to the unit $2$-sphere and hence $\mathbb R^{n-3}\times C(Y)=\mathbb R^n$.
\end{thm}
\begin{proof}
Fix a limit space $\mathbb R^{n-3}\times C(Y)$ as in the statement.  The first observation is that by Corollary \ref{c:codim3}, it follows that $Y$ is a smooth surface. 
 Indeed, if there were a point $y\in Y$ such that $r_h(y)=0$ then since $X= \mathbb R^{n-3}\times C(Y)$ 
it would follow that there is a set of codimension at least $2$ such that $r_h\equiv 0$, 
which cannot happen by Corollary \ref{c:codim3}.

Since  $|\RC_{M^n_j}|\to 0$,
$X$ is a metric cone with $\RC=0$ on the smooth part of $X$. By Proposition \ref{cone} 
$X$ is actually Ricci flat.  It follows that $Y$ is a smooth Einstein manifold satisfying $\Ric^g_Y= g$.  
Because $Y$ 
is a surface, this means in particular that $Y$ has constant sectional curvature $\equiv 1$.  
Thus,  either $Y=\mathbb R\mathbb P^2$ or $Y= S^2$, the unit $2$-sphere, and  in the latter case we are done. 

The remainder of the argument follows Cheeger-Naber \cite{CheegerNaber_Ricci}.  So let us study the case $Y= \mathbb R\mathbb P^2$.  Note that away from the singular set, $\mathcal S\equiv \mathbb R^{n-3}\times \{0\}$, we have that the $M_i^n$ 
converge to $\mathbb R^{n-3}\times C(Y)$ in $C^{1,\alpha}$.  For $\epsilon>0$ small, choose 
$u_i:B_2(p_i)\to \mathbb R^{n-3}$ to be an $\epsilon$-splitting as in Lemma \ref{h_splitting}.  
If $f_i:B_2(p)\to B_2(p_j)$ denotes the associated Gromov-Hausdorff 
maps, we put $\mathcal S_i=f_i(\mathcal S)$.  Then for $\tau>0$ small but fixed, we  have for $i$ 
sufficiently large, that on $B_1(p)\setminus T_\tau(\mathcal S_i)$, the estimates $|\nabla u_i|>\frac{1}{2}$ and 
$|H_i|\leq 1$ hold.  

Consider Poisson approximation $h_i$ 
to the square of distance
 function $d^2(x,p_i)$ on $B_2(p_i)$.  That is, $\Delta h_i = 2n$ and $h_i=4$ on $\partial B_2(p_i)$.  
We have (see for instance \cite{ChC1}) that $|h_i-d^2(\cdot,p_i)|\to 0$ uniformly in $B_2(p_i)$.  Furthermore, using the $C^{1,\alpha}$ convergence, we have for $i$ sufficiently large that 
$|\nabla h_i|>\delta$ and $|\nabla^2 h_i|\leq 4n$ on $B_1(p_i)\setminus B_\tau(\mathcal S_i)$.  Once again,
appealing to the $C^{1,\alpha}$ convergence, for all $i$ sufficiently large and any regular value 
$s_i\in B_1(0^{n-3})$ (which exists by Sard's Theorem) we have that $u_i^{-1}(s_i)\cap h_i^{-1}(1)$ is diffeomorphic to $\mathbb R\mathbb P^2$. Then for $i$ sufficiently large, 
$u^{-1}_i(s_i)\cap\{h\leq 1\}$
is a smooth $3$-manifold, whose boundary is diffeomorphic to $\mathbb R\mathbb P^2$.  However, the second Stiefel-Whitney
number of $\mathbb R\mathbb P^2$ is nonzero, and in particular, $\mathbb R\mathbb P^2$ 
does not bound a smooth $3$-manifold.  This contradicts $Y= \mathbb R\mathbb P^2$, finishing the proof.
\end{proof}

\begin{thm}\label{s:eps_reg}
Given $n \in \mathbb N$ and $\rv > 0$, there exists $\epsilon(n,\rv)>0, r(n,\rv) > 0, C(n,\rv) > 0$ such that if $(M^n, g, H, p)$ satisfies 
$|\RC_{M^n}|<\epsilon, \Vol(B_1(p))>\rv$ and 
$$d_{GH}(B_2(p),B_2(0))<\epsilon,$$
where $0$ is a vertex of the cone $\mathbb R^3\times C(Y)$
for some metric space $Y$, then we have $$r_h(p)>r.$$  Furthermore, if $(g,H)$ is generalized Einstein, we have 
$$\sup_{B_{1}(p)}\brs{\Rm}\leq C.$$
\end{thm}

\begin{proof}
Given $n$ and $\rv > 0$, choose $r$ as in Lemma \ref{regularity} and assume no such $\epsilon$ exists.  Then there exists a sequence of spaces $(M^n_i, g_{i},H_i,p_i)$ such that $|\RC_{M^n_i}|<\epsilon_i\to 0, \Vol(B_1(p_i))>\rv$ and
$$d_{GH}(B_2(p_i),B_2(0_i))<\epsilon_i\to 0$$
where $0_i\in \mathbb R^3\times C(Y_i)$ is a vertex, but $r_h(p_i)<r$.  After passing to a subsequence, we have $$B_2(p_i)\to B_2(0)$$ 
where $0\in\mathbb R^3\times C(Y)\equiv X$ a vertex.  By Theorem \ref{t:codim3}, $Y$ is isometric to the unit $2$-sphere, and so
$$B_2(p_j)\to B_2(0^n)\subset \mathbb R^n.$$
Applying the $\epsilon$-regularity of Lemma \ref{regularity}, we conclude $r_h(p_j)\geq r$ which is a contradiction.
\end{proof}
\begin{rmk}By a scaling argument, it is easy to see that the above theorem also holds if we only assume $|\RC|<n-1$.
\end{rmk}
\vspace{.3cm}

\subsection{Quantitative Stratification and Effective Estimates}\label{s:qs_estimates}

Having shown in Theorem \ref{t:codim4} that noncollapsed limits of sequences of manifolds with bounded generalized
 Ricci curvature are smooth away from a closed codimension $4$ subset, we will now give an application. 
 In particular, we will use the ideas of quantiative stratification first introduced in \cite{CheegerNaber_Ricci} 
in order to improve the codimension estimates on singular sets of limit spaces
to curvature estimates on manifolds with bounded generalized Ricci curvature.  We begin here by reviewing the quantitative stratification and the main results from \cite{CheegerNaber_Ricci}. 
These will play a crucial role in our estimates. 

\begin{defn}\label{d:quant_symmetry}
Given a metric space $Y$ with $y\in Y$, $r>0$ and $\epsilon>0$, we say
\begin{enumerate}
\item $B_r(y)$ is \emph{$(k,\epsilon)$-symmetric} if there exists a pointed $\epsilon r$-GH map $\iota:B_r(0^k,z_c)\subseteq \dR^k\times C(Z)\to B_r(y)\subseteq Y$ with $\iota(0^k,z_c)=y$.
\item $B_r(y)$ is \emph{$(k,\epsilon)$-symmetric} with respect to $\cL^k_\epsilon\subseteq B_r(y)$ if $\cL^k_\epsilon \equiv \iota\big(\dR^k\times \{z_c\}\big)\cap B_r(y)$.
\end{enumerate}
\end{defn}

To state the definition in words, we say 
that $Y$ is $(k,\epsilon,r)$-symmetric if the ball $B_r(x)$ looks $\ge$-close to having $k$-symmetries.  
The quantitative stratification is then defined as follows:

\begin{defn}
For each $\epsilon>0$, $0<r<1$ and $k\in\mathbb N$, define the closed quantitative $k$-stratum, $\mathcal S^k_{\epsilon,r}(X)$, by
\begin{align}
\mathcal S^k_{\epsilon,r}(X)\equiv \{x\in X:\text{ for no $r\leq s\leq 1$ is $x$ a $(k+1,\epsilon,s)$-symmetric point}\}\, .
\end{align}
\end{defn}

\begin{thm}[Quantitative Stratification, \cite{CheegerNaber_Ricci}]\label{t:quant_strat}
Let $M^n$ satisfy $\Rc^g \geq -(n-1)$ with $\Vol(B_1(p))>{\rv}>0$.  Then for every $\epsilon,\eta>0$ 
there exists $C=C(n,\rv,\epsilon,\eta)$ such that
\begin{align*}
\Vol\left(T_r\left(\mathcal S^k_{\epsilon,r}(M)\cap B_1(p)\right)\right)\leq C r^{n-k-\eta}\, .
\end{align*}
\end{thm}

We combine Theorem \ref{s:eps_reg} and Theorem \ref{t:quant_strat}
 in order to prove an a priori $L^q$ estimate for curvature, $q < 2$.
 
\begin{thm}\label{t:main_estimate}
Given $n \in \mathbb N$, $q < 2$, $\rv > 0$, there exists $C=C(n,\rv,q)$ such that if $(M^n, g, H)$ satisfies $\brs{\RC_{M^n}}\leq n-1$ and $\Vol(B_1(p))>\rv>0$
then one has
\begin{align}\label{e:main_Rm_est}
\fint_{B_1(p)} \brs{\Rm}^q\leq C\, .
\end{align}
\end{thm}

\begin{proof}
Let $(M^n,g,H,p)$ satisfy $\brs{\RC_{M^n}}\leq n-1$
 and $\Vol(B_1(p))>\rv>0$.  We will first show that for every $q<2$ there exists $C=C(n,v,q)>0$ such that
\begin{align*}
\fint_{B_1(p)} r_h^{-2q} \leq C\, .
\end{align*}
Simultaneously, we will show that if $M^n$ is generalized Einstein, then this can be improved to
\begin{align*}
\fint_{B_1(p)} \til r_x^{-2q} \leq C\, ,
\end{align*}
where $\til r_x$ denotes the regularity scale at $x$.

Let $q<2$ be fixed and set $\eta=4-2q$, and fix $\ge > 0$ as in Theorem \ref{s:eps_reg}.  With this choice of constants, by Theorem \ref{t:quant_strat} there exists $C(n,v,q)$ such that
\begin{align}\label{e:curv_ein1}
\Vol(T_r(\{x\in \mathcal S^{n-4}_{\epsilon,2r}\cap B_1(p)\})) < Cr^{4-\eta}\, .
\end{align}
Note that by rescaling, we may regard the $\epsilon$-regularity theorem (Theorem \ref{s:eps_reg})
as stating that if $x$ is $(n-3,\epsilon,2r)$-symmetric, then $r_h> r$, and if $M^n$ is generalized
Einstein then $r_x>r$.  In fact, we have that if $x$ is $(n-3,\epsilon,s)$-symmetric for any 
$s\geq 2r$, then $r_h>r$.  Thus, if $x\not\in \mathcal S^{n-4}_{\epsilon,2r}$, then $r_h>\frac{s}{2}> r$.  
The contrapositive gives the inclusion
\begin{align}
\{x\in B_1(p):r_h\leq r\} \subseteq \mathcal S^{n-4}_{\epsilon,2r}\cap B_1(p)\, ,
\end{align}
which by (\ref{e:curv_ein1}) implies the desired estimate
\begin{align} \label{f:curvest10}
\Vol(T_r(\{x\in B_1(p):r_h\leq r\})) < Cr^{4-\eta}\leq Cr^{2q}\, .
\end{align}
If $M^n$ is generalized Einstein, then Theorem \ref{s:eps_reg} allows us to replace $r_h$ with $\tilde r_x$, as claimed.

Now, for $q<2$, let us prove the $L^q$ bound on the curvature.  
For this note that if $r_h(x)>r$ then by definition there exists a harmonic coordinate system,
$\Phi:B_r(0^n)\to M$ 
with $\Phi(0)=x$ and such that the components of $g$ in these coordinates satisfy
\begin{align*}
||g_{ij}-\delta_{ij}||_{C^0(B_r(0))}+r||\partial_k g_{ij}||_{C^0(B_r(0))}<10^{-3}.
\end{align*}
Since the generalized Ricci curvature satisfies the 
bound $|\RC_{M^n}|\leq n-1$, by Lemma \ref{H-control}, combined with harmonic radius lower bound, we have $\brs{\Rc_{M^n}}<C(\rv,n)r^{-2}$, where the $r^{-2}$ factor is due to scaling. This implies that
\begin{align*}
|\Delta_x g_{ij}| < C(n)r^{-2}\, ,
\end{align*}
where $\Delta_x$ denotes the Laplacian written in coordinates.
Then for every $\alpha<1$ and 
$s<\infty$, we have the scale-invariant estimates
\begin{align*}
&r^{1+\alpha}||\partial_k g_{ij}||_{C^\alpha(B_{\frac{3r}{4}}(0))}\leq C(n,\alpha),\\
&r^{2}||g_{ij}||_{W^{2,s}(B_{\frac{3r}{4}}(0))}\leq C(n,s).
\end{align*}
In particular, applying this to $s=q$ we get 
\begin{align}\label{e:main_est:1}
r^{2q}\fint_{B_{r/2}(x)}\brs{\Rm}^q\leq C(n)r^{2q}\fint_{B_{3r/4}(0)}|\Phi^*\Rm|^q < C(n,q).
\end{align}

We use this estimate together with a covering argument to finish the proof.
Put $\eta=2-q$.  Then $q+\frac{\eta}{2}<2$ and it follows from (\ref{f:curvest10}) that 
\begin{align}\label{e:main_est:2}
\Vol(T_r(\{x\in B_1(p):r_h\leq r\})) < Cr^{2q+\eta}\, ,
\end{align}
for $C(n,v,q)>0$.  Consider the covering $\{B_{\tfrac{1}{2} r_h(x)}(x)\}$ of $B_1(p)$, and using the Vitali covering lemma we refine this to a subcovering denoted $\{B_{r_j}(x_j)\}$ where
\begin{enumerate}
\item $B_1(p)\subseteq \bigcup B_{r_j}(x_j)$, where $r_j=\frac{1}{2}r_h(x)$.
\vskip1mm
\item $\{B_{r_j/4}(x_j)\}$ are disjoint.
\end{enumerate}
By using \eqref{e:main_est:2} and volume comparison, we see that for each $\alpha\in\mathbb N$, we have
\begin{align}
\sum_{2^{-\alpha-1}<r_j\leq 2^{-\alpha}} \Vol(B_{r_j}(x_j)) \leq C r_j^{2q+\eta} = C\,r_j^{2q}\, 2^{-\eta\alpha}\, .
\end{align}
By summing over $\alpha$, this gives 
\begin{align}
\sum r_j^{-2q}\Vol(B_{r_j}(x_j))\leq C\sum 2^{-\eta\alpha}\leq C(n,v,q)\, .
\end{align}
Finally, combining this with \eqref{e:main_est:1} we get 
\begin{align*}
\fint_{B_1(p)}\brs{\Rm}^q &\leq \rv^{-1}\sum \int_{B_{r_j}(x_j)}\brs{\Rm}^q \notag\\
&\leq C(n,v,q)\sum r_j^{-2q}\Vol(B_{r_j}(x_j))\\
&\leq C(n,v,q)\, , 
\end{align*}
which finishes the proof of Theorem \ref{t:main_estimate}. 
\end{proof}

\vspace{.5cm}

\section{$L^2$ curvature estimate}

In this section we extend Theorem \ref{t:main_estimate} to yield an a priori $L^2$ curvature estimate in the case of bounded generalized Ricci tensor. Before we proceed, we outline the rough ideas of the proof of Theorem \ref{t:L2curvature}: 
\begin{enumerate}
\item The first ingredient is the neck decomposition theorem proved by \cite{JiNa}, where the unit ball is decomposed into the neck part and $\epsilon$-regular part and a key measure estimate is also proved (see Theorem \ref{t:neck_decomposition}).  Then we bound the $L^2$ curvature on the neck region and on the $\epsilon$-regular part separately. We remark that the curvature estimate on the $\epsilon$-regular part is already implicit in the work above, and the main new step is to obtain the $L^2$ curvature estimate on the neck region.
\item 
The second ingredient is a tool from \cite{JiNa} (see Lemma \ref{l:curvature_level_sets}), which roughly enables us to estimate $\brs{\Rm}$ pointwisely by $\brs{\Ric}$ on the neck region modulo the so-called $\cH$ volume, which can be controlled with only Ricci lower bound assumption on the neck region. The estimate is summarized in Theorem \ref{t:L2_neck}.
\item The third ingredient is the main new material when compared with \cite{JiNa}. By item $(2)$ above, to estimate the $L^2$ integral of $\brs{\Rm}$ on the neck region, we only need to estimate the $L^2$ integral of $\brs{\Ric}$ on the neck region. By using the bounded generalized Ricci condition, this is equivalent to estimate the $L^4$ integral of $H$ tensor on the neck region. This is achieved in Theorem \ref{t:splitting_neck_summable_H} by deriving a superconvexity estimate for $H$. 
\item In the key estimate of $\int \brs{H}b^{-3}$ in Theorem \ref{t:splitting_neck_summable_H}, a term of the form $\int\brs{\Rm}\brs{H}b^{-1}$ appears. If we have $\brs{\Rm}\le\epsilon^2b^{-2}$ on the neck region, which is indeed true if we are dealing with generalized Einstein case, then we have $\int \brs{\Rm}\brs{H}b^{-1}\le \int\epsilon^2\brs{H}b^{-3}$, which can be absorbed. In the bounded generalized Einstein case,  we can use a H\"older inequality to deal with $\int\brs{\Rm}|H|b^{-1}$, which leads to a term $\int \brs{\Rm}^2\approx \int |H|^4$ by item (2). It seems that we land in a loop, i.e. using $L^4$ norm of $|H|$ to bound $\int b^{-3}|H|$.  However through a detailed analysis we show that there is always a small multiple in front of $\int |H|^4$, and hence it can be absorbed. In the end, we use an induction argument on scales to show that $L^4$ norm of $H$ is small on the neck region. 
\end{enumerate}

\subsection{Preliminaries on neck regions}

In this subsection, we define  neck regions and collect some recent structural results about them from \cite{JiNa} and \cite{CJN}. A crucial point for us is the Ahlfors regularity estimate, which holds on the neck region  with only a Ricci lower bound.\\

\begin{defn}\label{d:neck}
We call $\cN\subseteq B_2(p)$ a \emph{$\delta$-neck region} if there exists a closed subset $\cC = \cC_0\cup \cC_+=\cC_0\cup\{x_i\}\subset B_2(p)$  with $p\in \cC$ and a radius function $r:\cC\to \dR^+$ on $\cC_+$ with $r_x=0$ on $\cC_0$ such that $\cN \equiv B_2\setminus \overline B_{r_x}(\cC)$ satisfies
\begin{enumerate}
	\item[(n1)] $\{B_{\tau(n)^2 r_x}(x)\}$ are pairwise disjoint where $\tau(n)=10^{-10n}\omega_n$ and $|\cV_{\delta^{-1}}(x)-\cV_{\delta r_x}(x)|\le\delta^2.$
	\item[(n2)] For each $r_x\leq r\leq 1$ there exists a $\delta r$-GH map $\iota_{x,r}:B_{\delta^{-1}r}(0^{n-4},y_c)\subseteq\dR^{n-4}\times C(S^3/\Gamma)\to B_{\delta^{-1}r}(x)$, where $\Gamma\subseteq O(4)$ is nontrivial  and $\iota_{x,r}(0^{n-4},y_c)=x$. 
	\item[(n3)] For each $r_x\leq r$ with $B_{2r}(x)\subseteq B_2(p)$ we have that $\cL_{x,r}\equiv \iota_{x,r}\big(B_r(0^{n-4})\times\{y_c\}\big)\subseteq B_{\tau r}(\cC)$ and $\cC\cap B_r(x)\subset B_{\delta r}(\cL_{x,r})$.
	\item[(n4)] $|\Lip\, r_x|\leq \delta$.
\end{enumerate}	
For each $2\tau(n)\leq s$ we define the region $\cN_s\equiv B_2\setminus \overline B_{s\cdot r_x}(\cC)$. 
\end{defn}

\begin{rmk} In $(n1)$,
 we define the {\it volume ratio} by
 \begin{align}
 \label{e:volratio}
 \cV^\delta_r(x) := \frac{\Vol(B_r(x))}{\Vol^{-\delta}(B_r)}\, ,
 \end{align}
 where $\Vol^{-\delta}(B_r)$ denote the volume of an $r$-ball in a simply connected space of constant curvature $M^n_{-\delta}$.
\end{rmk}
\begin{rmk}
 By the nontriviality of $\Gamma$ in (n2), we know that when $\delta$ is small and $$d_{GH}(B_{\delta^{-1}r}(0^{n-4},y_c), B_{\delta^{-1}r}(x))<\delta r,$$ $B_{\delta^{-1}r}(x)$ is not $(n-3,\eta(n,\rv))$ symmetric (see Definition \ref{d:quant_symmetry}). Therefore fixing a small $\eta(n,\rv)$, our definition of neck region is compatible with the definition of $(n-4,\delta,\eta)$ neck region in \cite[Definition 2.4]{CJN}. 
\end{rmk}
It will be important to measure the codimension 4 measure of the singular set $\cC$ in the neck region. So we define the packing measure as:
\begin{defn}
	Let $\cN \equiv B_2\setminus \overline B_{r_x}(\cC)$ be a neck region, then we define the associated \emph{packing measure}
	\begin{align}
    \mu=\mu_\cN \equiv \sum_{x\in \cC_+} r_x^{n-4}\delta_{x} + \lambda^{n-4}|_{\cC_0}\, ,
\end{align}
where $\lambda^{n-4}|_{\cC_0}$ is the $n-4$-dimensional Hausdorff measure restricted to $\cC_0$.
\end{defn}
We next record a key Ahlfors regularity estimate  for neck regions:
\begin{thm}\label{Ahlfors regularity}\cite[Thm 2.9]{CJN}
	Let $(M^n,g,p)$ satisfy $\Vol(B_1(p_j))>\rv>0,\ \Ric^g \geq -\delta$, and suppose $\cN = B_2(p)\setminus \bigcup_{x\in \cC} \overline B_{r_x}(x)$ is a $\delta$-neck region with $\delta<\delta'(n,\rv)$. There exists a constant $B(n)$ so that for each $x\in \cC$ and $r_x<r$ with $B_{2r}(x)\subseteq B_2$ we have
	\begin{equation}\label{Ahlfors}
	B(n)^{-1}r^{n-4} < \mu\big(B_r(x)\big)< B(n)r^{n-4}\, .
	\end{equation}
\end{thm}
\vspace{.2cm}

The main structure theorem we will prove for the neck region is:   

\begin{thm}\label{t:neck_region}
Let $(M^n_j,g_j,H_j,p_j)\to (X,d,p)$ be a Gromov-Hausdorff limit with $\Vol(B_1(p_j))>\rv>0$, and fix $\ge > 0$.  Then for $\delta\leq \delta'(n,\rv,\epsilon)$, if $|\RC_{M^n_j}|\leq\delta$ and $\cN = B_2(p)\setminus \overline B_{r_x}(\cC)$ is a $\delta$-neck region, then the following hold:
\begin{enumerate}
	\item For each $x\in \cC$ and $r>r_x$ such that $B_{2r}(x)\subseteq B_2(p)$ the induced packing measure $\mu$ satisfies the Ahlfors regularity condition $ B(n)^{-1}r^{n-4}<\mu(B_r(x)) <B(n)r^{n-4} $.
	\item $\cC_0$ is $n-4$ rectifiable.
	\item $X$ is a $W^{2,p}$ manifold on $\cN$. 
	\item $\int_{\cN\cap B_1(p)} \brs{\Rm}^2(x)\,dx < \epsilon$.
\end{enumerate}
\end{thm}

Note that items $(1),(2)$ and  $(3)$ are the main theorems of \cite{CJN}, while item $(4)$ is new and will be proved at the end of this section. 

\vspace{.2cm}

\subsection{Cutoff Functions on Neck Regions}

In this subsection we build a natural cutoff function associated with a neck region.  We will record some of the basic properties and estimates associated to it, which will be useful throughout the paper.

\begin{lemma}\label{l:neck:cutoff}\cite{JiNa} 
	Let $(M^n,g,p)$ satisfy $\Vol(B_1(p))>\rv>0$ with ${\Ric^g}>-\delta$ and $\cN\subseteq B_2(p)$ a $\delta$-neck region.  Then for  $\delta<\delta'(n,\rv)$ there exists a cutoff function $\phi_1\cdot \phi_2 \equiv \phi_\cN: B_2\to [0,1]$  such that
	\begin{enumerate}
	\item $\phi_1(y) \equiv 1$ if $y\in B_{18/10}(p)$ with $\supp(\phi_1) \subseteq B_{19/10}(p)$.
	\item $\phi_2(y)\equiv 1$ if $y\in \cN_{10^{-3}}$ with $\supp(\phi_2)\cap B_2 \subseteq \cN_{10^{-4}}$.
	\item $|\nabla \phi_1|$, $|\Delta \phi_1|\leq C(n)$.
	\item $\supp(|\nabla \phi_2|)\cap B_{19/10}\subseteq  A_{10^{-4}r_x,10^{-3}r_x}\big(\cC\big)$ with $r_x|\nabla \phi_2|$, $r^2_x|\nabla^2\phi_2|<C(n,\rv)$ in each $B_{r_x}(x)$.
	\end{enumerate}
\end{lemma}
\begin{proof}
Recall from \cite{ChC1} that for each annulus $A_{s,r}(x)$ we can build a cutoff function $\phi$ such that
\begin{gather}\label{e:cutoff}
\begin{split}
    &\phi\equiv 1 \text{ on }B_s(x)\, ,\;\;\phi\equiv 0 \text{ outside }B_r(x)\, \\
	&|\nabla \phi|\leq \frac{C(n)}{|r-s|}\, ,\;\;|\Delta \phi|\leq \frac{C(n)}{|r-s|^2}\,\, .
\end{split}
\end{gather}
It will be important to recall that $\phi$ is built as a composition $\phi=c\circ f$, where $c$ is a smooth cutoff on $\dR$ with $c=1$ on $[0,s^2]$ and $c=0$ outside $[0,r^2]$, and $f$ satisfies $\Delta f=2n$ on $B_{5r}(x)$ and is uniformly equivalent to the square distance $d^2_x$.

Thus we can let $\phi_1$ be the cutoff associated to the annulus $A_{\frac{18}{10},\frac{19}{10}}(p)$.  To build $\phi_2$ let us begin by defining for each $x\in \cC$ the cutoff $\phi_x=c_x\circ f_x$ as in \eqref{e:cutoff} associated to the annulus $A_{10^{-4}r_x,10^{-3}r_x}(x)$.  Using elliptic estimates and Bochner formula we can get the pointwise estimate $|\nabla^2 f_x|<C(n) r_x^{-2}$ on $A_{10^{-4}r_x,10^{-3}r_x}(x)\cap \cN_{10^{-4}}$ .  In particular we can get the estimate
\begin{align}\label{e:cutoff:1}
	r_x|\nabla \phi_x|,\, r_x^2|\nabla^2 \phi_x|<C(n) \text{ on }A_{10^{-4}r_x,10^{-3}r_x}(x)\cap \cN_{10^{-4}}\, .
\end{align}

Now with $\tau<10^{-5}$ in the definition of neck regions let us define the cutoff function
\begin{align*}
\phi_2 \equiv \prod \Big(1-\phi_x\Big)\, .	
\end{align*}
We have that $\text{supp}\,\phi_2 \cap B_2 \subseteq \cN_{10^{-4}}$ and for each point $y\in B_2(p)$ we have by $(n1)$ that there are at most $C(n,\rv,\tau(n))$ of the cutoff functions $\phi_x$ which are nonvanishing at $y$.  In particular, the estimates \eqref{e:cutoff:1} then easily imply the required estimates on $\phi_2$.
\end{proof}

\vspace{.2cm}
\subsection{Scale Invariant  Curvature Estimate on the Neck Region}\label{sss:weak_sobolev_convergence}

The main goal of this subsection is to collect several curvature estimates which will play an important role in subsequent sections.  Let us begin with the following, which tells us that the harmonic radius at a point in the neck region is roughly the distance of that point to the effective singular set.  The proof is immediate using $(n2)$ and $(n3)$ of a neck region and regularity Theorem \ref{s:eps_reg} proved in section 3, but it is worth mentioning the result explicitly:

\begin{lemma}\label{l:neck:distance_harm_rad}
	Let $(M^n,g,H,p)$ satisfy $\Vol(B_1(p))>\rv>0$ and $\brs{\RC}\le n-1$.  Then there exists $C(n,\rv)$ and $\gd'(n, \rv) > 0$ small such that if $\cN = B_2(p)\setminus \overline B_{r_x}(\cC)$ is a $\delta$-neck region with $\gd < \gd'(n,\rv)$, then for each $y\in \cN_{10^{-6}}$ one has
	\begin{align*}
		C(n,\rv)d(y,\cC)\leq r_h(y)\leq C(n,\rv)^{-1}d(y,\cC)\, .
	\end{align*}
Furthermore, for any $\epsilon>0$, if $\delta\le\delta'(n,\rv,\epsilon)$, it holds that	$$d(y,\cC)^{4-n}\int_{B_{d(y,\cC)}(y)\setminus B_{1/2d(y,\cC)}(y)}\brs{\Rm}^2\,dz\,\le\epsilon.$$
\end{lemma}
\begin{rmk}
 It is worth emphasizing that the $L^2$ curvature estimate above is not good enough for a global $L^2$ bound, as there is no $n-4$ control on the number of such balls needed to cover a neck region.  The proof of the $L^2$ bound on the whole neck region is much more subtle and requires a superconvexity estimate.
\end{rmk}

\begin{proof} 
We prove the lower bound first.  By $\Vol(B_1(p))>\rv>0$, we have an upper bound of the order of group $\Gamma$ in the definition of neck region.  Now suppose we don't have such  $C(n,\rv),\delta(v,\rv)$. Then we will have a sequence of neck regions $(\mathcal N_i,\delta_i)$ with $C_i,\delta_i\to 0$ and $y_i\in\mathcal N_{10^{-6}} $ such that $r_h(y_i)<C_id(y_i,\cC)$. By the scaling invariant property of condition $(n2),(n3)$, we may assume that for some $x_i\in\mathcal C_i, r_{x_i}=1/2, d(y_i,x_i)= (1/2) 10^{-6}$, this implies $r_h(y_i)\to 0$. On the other hand $B(x_i,3/4)\to\dR^{n-4}\times C(S^3/\Gamma)$ in the Gromov-Hausdorff sense and $y_i\to y$ with $d(y,\dR^{n-4}\times 0):=\epsilon>0$. This implies that on the flat cone $\dR^{n-4}\times C(S^3/\Gamma)$, $B(y,\epsilon/2)$ is a Euclidean ball. 
 Then by Theorem \ref{s:eps_reg}, $r_h(y_i)$ will stay away from $0$, which is a contradiction. Now by the scale invariance of second statement of the theorem, we argue by contradiction and may assume that $r_y=1$. Then the annulus $B_1\setminus B_{1/2}$ will be $10^6\delta$-close to an annulus on the flat metric space $\mathbb R^{n-4}\times C(Y/\Gamma)$. Using Theorem \ref{s:eps_reg} again yields the result.
\end{proof}
\vspace{.2cm}
The above lemma gives a scale invariant estimate of $L^2$  curvature on a regular ball. By using Ahlfors regularity, in the next two lemmas, we prove summable $L^2$  curvature estimates in two different cases, which will play a role in subsequent sections.  First, we prove in the following lemma that the $L^2$ curvature on annuli with centers $x$ in $\cC$ and size comparable to $r_x$ is small.  \\

\begin{lemma}\label{t:neck:annulusRm}
	Let $(M^n,g,H,p)$ satisfy $\Vol(B_1(p))>\rv>0$ and $\brs{\RC}\le n-1$. For any $\epsilon>0$, 
	if $\delta\le \delta'(n,\epsilon,\rv)$ and $\cN = B_2(p)\setminus \overline B_{r_x}(\cC)$ is a $\delta$-neck region, then it holds that
\begin{align*}
\int_{(\cN_{\tau(n)}\setminus \cN)\cap B_1(p)} \brs{\Rm}^2(z)\,dz\le C(n,B)\epsilon.
\end{align*} \end{lemma}
\begin{proof}
Recall that by our definition, $\cN_{s}\equiv B_2\setminus \overline B_{s\cdot r_x}(\cC)$, then we have 

\begin{align*}
\int_{(\cN_{\tau(n)}\setminus \cN)\cap B_1(p)} \brs{\Rm}^2(z)\,dz\le&\ \sum_{x\in B_1(p)\cap \cC}\int_{B_{r_x}(x)\setminus \cup_{y\in \cC}B_{\tau(n)r_y}(y)}\brs{\Rm}^2(z)\,dz\\
\le&\ \sum_{x\in B_1(p)\cap \cC} C(n,B) \epsilon r_x^{n-4}\le C(n,B)\epsilon,
\end{align*}

where we have used Lemma \ref{l:neck:distance_harm_rad} for second inequality and Ahlfors regularity Theorem \ref{Ahlfors regularity} for the last inequality.
\end{proof}

The following lemma shows that the $L^2$ curvature is not only scale invariantly small on the neck region but also integrably small at a fixed scale on the neck region.

\begin{lemma}\label{l:neck:scaleinvariantsmall}
Let $(M^n,g,H,p)$ satisfy $\Vol(B_1(p))>\rv>0$ with $|{\RC}|\le n-1$.
For any $\epsilon>0$, if  $\delta<\delta'(n,\epsilon,\rv)$ are such that $\cN = B_2(p)\setminus \bigcup_{x\in \cC} \overline B_{r_x}(x)$ is a $\delta$-neck region, then
\begin{align*}
    \int_{\mathcal N\cap \{2^{-i}\le d(x,\mathcal C)\le 2^{-i+1}\}}\brs{\Rm}^2\le\epsilon
\end{align*}
\end{lemma}
\begin{proof}
It follows from Lemma \ref{l:neck:distance_harm_rad} that if $2^{-i}\le d(y,\mathcal C)\le 2^{-i+1}$, then 
\begin{equation}\int_{B(y,2^{-i-1})}\brs{\Rm}^2\,dz\,\le\epsilon 2^{-i(n-4)}.
\end{equation}

Choose a Vitali covering $\{B_{2^{-i+3}}(x_\alpha),x_\alpha\in \cC,\alpha=1,\cdots, N\}$ of $\text{supp }\phi_\cN \cap \{2^{-i}\le d(x,\cC)\le 2^{-i+1}\}$ such that $\{B_{2^{-i}}(x_\alpha),x_\alpha\in \cC,\alpha=1,\cdots,N\}$ are pairwise disjoint. We get 
\begin{align}\label{e:dcC-321}
 C(n,B)\geq \mu(B_{19/10}(p))\geq \mu(\text{supp }\phi_\cN)\geq\sum_{\alpha=1}^N\mu(B_{2^{-i}}(x_\alpha))\geq N B^{-1}2^{-(n-4)i}6^{-(n-4)}.
\end{align}
We should remark that the last inequality in \eqref{e:dcC-321} always holds even $2^{-i}\le r_{x_\alpha}$ or $2^{-i}\geq d\big(x_\alpha,\partial B_2(p)\big)$. Actually, for the first case ($2^{-i}\le r_{x_\alpha}$) the inequality is obvious from the definition of $\mu$, for the second case ($2^{-i}\geq d\big(x_\alpha,\partial B_2(p)\big))$ we can find $B_{2^{-i}/3}(\tilde{x}_\alpha)\subset B_2(p)\cap B_{2^{-i}}(x_\alpha)$ with $\tilde{x}_\alpha\in \cC$.
Therefore, we have $N\le C(n,B)2^{(n-4)i}$. Thus 
\begin{gather}\label{e:dcC-32}
\begin{split}
\Vol(\text{supp }\phi_\cN \cap \{2^{-i}\le d(x,\cC)\le 2^{-i+1}\})\le&\ \sum_{\alpha=1}^N\Vol(B_{2^{-i+3}}(x_\alpha))\\
\le&\ C(n,B)2^{(n-4)i} C(n)2^{n(-i+3)}\\
\le&\ C(n,B)2^{-4i}.
\end{split}
\end{gather}

Now another  standard Vitali covering argument  using volume comparison  will show that the set $(\text{supp }\phi_\cN \cap \{2^{-i}\le d(x,\cC)\le 2^{-i+1}\})$ can be covered by  at most $C(n,\rv)2^{-i(4-n)}$ number of balls with center on $(\text{supp }\phi_\cN \cap \{2^{-i}\le d(x,\cC)\le 2^{-i+1}\})$ and radius $2^{-i-1}$. Therefore,
 \begin{equation}\int_{\mathcal N\cap \{2^{-i}\le d(x,\mathcal C)\le 2^{-i+1}\}}\brs{\Rm}^2\le\epsilon C(n,\rv)2^{-i(n-4)}2^{i(n-4)}\le \epsilon C(n,\rv),
 \end{equation}
 which completes the proof.
\end{proof}

\vspace{.2cm}

\subsection{Estimates of Green's Functions on Neck Regions}\label{ss:splitting:Greens_standard}

In this section we discuss Green's functions on neck regions.  We will use these estimates in the next section to discuss Green's functions with respect to the packing measure $\mu$.  
The proofs are nearly verbatim as in \cite{JiNa} except in the proof of Lemma \ref{l:green_function_cone_point}, where we use two-sided generalized Ricci bound instead of two-sided regular Ricci to conclude the $C^1$ convergence of metric.  To begin we recall fundamental heat kernel estimates on a general smooth Riemannian manifolds, which follow from the results in \cite{LiYau_heatkernel86}, \cite{SY_Redbook}, \cite{SoZha},\cite{Hamilton_gradient}, \cite{Kot_hamilton_gradient}: 
\begin{thm}[Heat Kernel Estimates]\label{t:heat_kernel}
	Let $(M^n,g,x)$ be a pointed Riemannian manifold with $\Vol(B_1(x))\geq \rv>0$ and ${\Ric}^g \geq -(n-1)$. Then for any $0<t\le 10$ and $\epsilon>0$, the heat kernel $\rho_t(x,y)=(4\pi t)^{-n/2}e^{-f_t}$ satisfies for all $y\in M$ that
	\begin{enumerate}
	\item $C(n,\epsilon) t^{-n/2} e^{-\frac{d(x,y)^2}{(4-\epsilon)t}}\leq \rho_t(x,y)$. 
	\item $\rho_t(x,y)\leq C(n,\rv,\epsilon) t^{-n/2} e^{-\frac{d(x,y)^2}{(4+\epsilon)t}}$.
	\item $t|\nabla f_t|^2\leq C(n,\rv)\Big(1+\frac{d^2(x,y)}{t}\Big)$ . 
	\item $-C(n,\rv,\epsilon)+\frac{d^2(x,y)}{(4+\epsilon)t}\leq f_t \leq C(n,\rv,\epsilon)+\frac{d^2(x,y)}{(4-\epsilon)t }$
	\end{enumerate}
\end{thm}
\vspace{.2cm}

Now we define Green's functions with respect to the packing measure $\mu$  by using standard Green's function.
\begin{defn}\label{d:Green_neck}
	Given $(M^n,g,p)$ a compact Riemannian manifold with $\cN\subseteq B_2(p)$ a $\delta$-neck region with packing measure  $\mu_\cN$, let us consider $\bar \mu=\mu_{\cN}|_{B_{39/20}(p)}$ and we define the following:
	\begin{enumerate}
		\item We denote by $G_x(y)$ a Green's function at $x$.  That is, $-\Delta_y G_x(y) = \delta_x, G_x(y)|_{\partial B_4(x)}=f(y)$, where $\delta_x$ is the Dirac delta at $x$ and $f$ is a smooth function to be specified later.   
		\item We denote by $G_{\bar \mu}(y)$ the function $G_{\bar \mu}(y) = \int_{B_{39/20}(p)} G_x(y)\, d\mu(x)$  a Green's function which solves $-\Delta G_{\bar\mu} = \bar \mu$.
		\item We denote by $b(y)\equiv G^{-1/2}_{\bar \mu}$ the $\bar \mu$-Green's distance function to the center points $\cC$. 
	\end{enumerate}
\end{defn}

The intuition is that $b$ should behave in a manner which is comparable to the Green's function from the singular set in $\dR^{n-4}\times C(S^3/\Gamma)$, which itself is a multiple of $d(\cS^{n-4},\cdot)^{-2}\approx b^{-2}$.  Our main goal in this subsection is to prove estimates to this effect:

\begin{prop} \label{p:green_function_distant_function}
	Let $(M^n,g,H,p)$ satisfy $\Vol(B_1(p))>\rv>0$ with $\brs{\RC} \leq \gd$.  Then for given $B>0$ there exists  $\delta'(n,B,\rv)$,  $C(n,B,\rv)>0$ and a Green function $G_{\bar \mu}$ such that if $\cN = B_2(p)\setminus \overline B_{r_x}(\cC)$ is a $\delta$-neck region, $\delta<\delta'$, then
	\begin{enumerate}
	\item For $x\in \cN_{10^{-6}}\cap B_{39/20}(p)$, we have $C^{-1}\,d(x,\cC)\leq b(x)\leq C\,d(x,\cC)$.
	\item For $x\in \cN_{10^{-6}}\cap B_{39/20}(p)$, we have $C^{-1}\leq |\nabla b|\leq C$.
	\end{enumerate}
\end{prop}

Before giving the proof of Proposition \ref{p:green_function_distant_function}, we record a series of lemmas.  First, let us point out that the advantage of $\mu$ after restricting to $B_{39/20}(p)$ is the following uniform Ahlfors regularity:

\begin{lemma} \label{l:Ahlfors}
Let $(M^n,g,H,p)$ satisfy $\Vol(B_1(p))>\rv>0$ with $\brs{\RC} \leq \gd$.  If $\cN = B_2(p)\setminus \overline B_{r_x}(\cC)$ is a $\delta$-neck region, $\delta<\delta'(n,\rv)$, then for all $x\in \cC\cap B_{39/20}(p)$, $r_x\le r\le 1$ we have 
\begin{align}
C(n)^{-1}B^{-1}r^{n-4}\le \bar\mu(B_r(x))\le C(n)Br^{n-4}.
\end{align}
\end{lemma}
\begin{proof}
For the upper bound, if $r\le 1/40$  it follows directly by inequality (\ref{Ahlfors}) since $B_{2r}(x)\subset B_2(p)$. If $r\geq 1/40$, it follows from $\bar \mu(B_r(x))\le \mu(B_{39/20}(p))\le C(n)B$ where the last inequality is based on inequality (\ref{Ahlfors}) and a Vitali covering $\{B_{1/40}(x_\alpha),x_\alpha\in \cC\cap B_{39/20}(p)\}$ of $B_{39/20}(p)$. 
On the other hand,  the lower bound follows directly by inequality (\ref{Ahlfors}) and the fact that there exists $B_{r/3}(\tilde{x})\subset B_r(x)\cap B_2(p)$ with $\tilde{x}\in \cC$.
\end{proof}

Now let us collect together
a list of useful computations from Lemma 5.9 of \cite{JiNa} before proving Proposition \ref{p:green_function_distant_function}.
\begin{lemma}\label{l:green_formula}
	Given $(M^n,g,p)$ with $\cN\subseteq B_2(p)$ a $\delta$-neck region and $b(y)$ an associated Green's distance function, then for a smooth compactly supported function $f:B_4(p)\rightarrow \dR$ the following hold:
\begin{enumerate}
 \item ${\bar \mu}[f] = -\int_{b\leq r} G_{{\bar \mu}}\Delta f(z)\,dz +2 r^{-3}\int_{b=r}f|\nabla b|(x)\, dx + r^{-2}\int_{b\leq r}\Delta f(z)\,dz$.
 \item ${\bar \mu}[f] = -\int_{b\leq r} G_{{\bar \mu}}\Delta f(z)\,dz +2 r^{-3}\int_{b=r}f|\nabla b|(x)\, dx + r\frac{d}{dr}\left(r^{-3}\int_{b=r}f|\nabla b|(x)\,dx\right)$,
 \item $\frac{d}{dr}\Big(r^{-3}\int_{b=r} f\, |\nabla b|(z)\,dz\Big) = \, r^{-3}\int_{b\leq r}\Delta f(z)\,dz\,$
 \item $r\frac{d^2}{dr^2} \Big(r^{-3}\int_{b=r} f\, |\nabla b|(z)\,dz\Big)+3\frac{d}{dr}\Big(r^{-3}\int_{b=r} f\, |\nabla b|(z)\,dz\Big) = \, r^{-2}\int_{b=r} \Delta f \,|\nabla b|^{-1}(z)\,dz$.
\end{enumerate}
\end{lemma}

\begin{proof}
Let us observe that $\big(b^{-2}-r^{-2}\big)$ vanishes on $b=r$, is smooth in a neighborhood of this set, and $\text{supp}\Big\{\Delta\big(b^{-2}-r^{-2}\big)\Big\}\subseteq \{b<r\}$.  Thus we can use standard properties of the distributional Laplacian to compute
\begin{align*}
\int_{b\leq r} \Delta f\left(b^{-2}-r^{-2}\right)(z)\,dz &= -{\bar \mu}\big[f\big]-\int_{b=r}f\langle \frac{\nabla b}{|\nabla b|},\nabla b^{-2}\rangle(z)\,dz \notag\\
	&= -{\bar \mu}\big[f\big] + 2r^{-3}\int_{b=r} f|\nabla b|(z)\,dz\, ,
\end{align*}
which proves the first formula.  To compute the second formula let us first compute
\begin{align}
\Delta b = 3b^{-1}|\nabla b|^2\, ,	
\end{align}
and recall that the mean curvature of the $b=r$ level set is given by $\text{div}\big(\frac{\nabla b}{|\nabla b|}\big)$.  Therefore we can compute

\begin{align*}\label{e:greens_formula:1}
\frac{d}{dr}\left(r^{-3}\int_{b=r}f|\nabla b|(z)\,dz\right) = &\ -3r^{-4}\int_{b=r}f|\nabla b| (z)\,dz+ r^{-3}\int_{b=r}\langle \nabla f, \frac{\nabla b}{|\nabla b|}\rangle(z)\,dz\\
&\ \qquad + r^{-3}\int_{b=r} f \langle \nabla|\nabla b|, \frac{\nabla b}{|\nabla b|^2}\rangle(z)\,dz + r^{-3}\int_{b=r} f \text{div}\big(\frac{\nabla b}{|\nabla b|}\big)(z)\,dz\notag\\
=&\ r^{-3}\int_{b\le r}\Delta f(z)\,dz -3r^{-4}\int_{b=r}f|\nabla b|(z)\,dz+ r^{-3}\int_{b=r} f \langle \nabla|\nabla b|, \frac{\nabla b}{|\nabla b|^2}\rangle(z)\,dz \notag\\
&\ \qquad + r^{-3}\int_{b=r} f \Big(\frac{\Delta b}{|\nabla b|}-\langle \nabla b, \frac{\nabla|\nabla b|}{|\nabla b|^2}\rangle\Big)(z)\,dz\notag\\
=&\ r^{-3}\int_{b\le r}\Delta f(z)\,dz\, .
\end{align*}
This implies formula (3), and then (2) is a formal consequence.  For item (4) we further differentiate the above formula to obtain
\begin{align*}
r\frac{d^2}{dr^2}\Big(r^{-3}\int_{b=r} f\, |\nabla b|(z)\,dz\Big)+3\frac{d}{dr}\Big(r^{-3}\int_{b=r} f\, |\nabla b|(z)\,dz\Big) &= \,r^{-2}\frac{d}{dr} \left(\int_{b\le r} \Delta f (z)\,dz\right)\\ \nonumber
&=\,r^{-2}\frac{d}{dr} \left(\int_0^rds\int_{b=s}\Delta f\,|\nabla b|^{-1} (z)\,dz\right)\\ \nonumber
&=r^{-2}\int_{b=r}\Delta f\,|\nabla b|^{-1}(z)\,dz,
\end{align*}
as claimed.
\end{proof}

In order to prove Proposition \ref{p:green_function_distant_function}, we first consider the Green function $G_x$ for any center point $x\in\cC$ in a neck region.  Besides the basic expected estimates, we need to see that at every point $y\in\cN$ in the neck region itself there is a fixed unit vector $v_y\in T_yM$ for which $\nabla_yG_x$ has a definite lower bound on the gradient in the direction of $v_y$.  This will be used heavily when we integrate to construct $G_{\cN}$ in order to see that the gradient of $G_{{\bar \mu}}$ has a definite lower bound in the neck region.  Precisely, we have the following:

\begin{lemma} \label{l:green_function_cone_point}
Let $(M^n,g,H,p)$ satisfy $\Vol(B_1(p))>\rv>0$ and $\brs{{\RC}}<\delta$ with $\cN=B_2(p)\setminus \bigcup_{x\in \cC} \overline B_{r_x}(x)$ a $\delta$-neck region.  For each $R>0$, if $\delta\leq \delta'(n,R,\rv)$  then there exists $C(n,\rv)>0$ such that $\forall$ $x\in \cC$ there exists a Green's function $G_x$ such that:
\begin{enumerate}
\item For all $y\in \cN_{10^{-5}}$  we have $C^{-1}d_x^{2-n}(y)\le G_x(y)\le Cd_x^{2-n}(y)$.
\item For all $y\in \cN_{10^{-5}}$  we have $|\nabla_yG_x|(y)\le Cd_x^{1-n}(y)$.
\item For all $y\in \cN_{10^{-5}}$ there exists a unit vector $v_y\in T_yM$ such that, setting $r=d(y,\cC)$, for all $x\in B_{Rr}(y)\cap \cC$ one has 
\begin{equation*}
\langle\nabla_yG_x(y),v_y\rangle\,>\, c(n,R) r^{1-n}.
\end{equation*}
\item In particular, if $y\in \cN_{10^{-5}}$ with $r=d(y,\cC)$ and $x\in B_{10r}(y)\cap \cC$ then
\begin{align*}
\langle\nabla_yG_x(y),v_y\rangle\,>\, c(n) r^{1-n}.
\end{align*}
\end{enumerate}
\end{lemma}
\begin{proof}
For each $x\in \cC$, since the estimates in (1) and (2) are scale invariant, in order to estimate the Green's function on the ball $B_{\delta^{-1/2}}(p)$ with $|{\RC_M}|<\delta$, it suffices to construct and estimate the Green function on the ball $B_1(p)$ with $|{\RC_M}|<1$.  In this case, we first define

$${G}_{x,0}(y)=\int_0^1\rho_t(x,y)dt\, .$$ Then by heat kernel estimate of Theorem \ref{t:heat_kernel}, we can compute
$C^{-1}d_x^{2-n}(y)\le G_{x,0}(y)\le Cd_x^{2-n}(y)$. Additionally, we can compute that $\Delta {G}_{x,0}(y)=\int_0^1\partial_t \rho_t(x,y)dt=\rho_1(x,y)-\delta_x(y).$ Therefore let us solve $G_{x,1}$ by
\begin{align}\label{e:green_rho1}
\Delta G_{x,1}(y)=\rho_1(x,y),
\end{align}
on $B_{4}(x)$ with $G_{x,1}(y)= G_{x,0}(y)$ on $\partial B_{4}(x)$. We define our Green's function
$$G_x=G_{x,0}-G_{x,1}\, .$$
We claim that $G_x$ satisfies the properties of the lemma.  Indeed, by noting the uniform bound on the heat kernel $\rho_1(x,y)$ we may use a standard maximal principle and Cheng-Yau gradient estimates on (\ref{e:green_rho1}) as in \cite{Cheeger01},   in order to show that $|G_{x,1}|+|\nabla G_{x,1}|<C(n)$ in $B_{4/3}(p)$ is uniformly bounded. Coupling with the estimates of $G_{x,0}$, we obtain the estimates on $G_x$. Thus we have proved (1) and (2).

Now we only give a proof of item (3), since the argument is the same for item (4). We argue by contradiction. Therefore assume for some $R>0$ that there exist $\delta_i$-neck regions $\cN_i$ with $\delta_i\to 0$, $y_i\in \cN_{i,10^{-5}}$ and $r_i=d(y_i,\cC_i)$ such that
\begin{align*}
\sup_{v\in T_{y_i}M,~|v|=1}\inf _{x_{i}\in B_{Rr_i}(y_i)\cap \cC_i}\langle\nabla_yG_{x_i}(y_i),v\rangle\,<\, i^{-1}r_i^{1-n}.
\end{align*}
 Scaling $B_{Rr_i}(y_i)$ to $\tilde{B}_{R}(y_i)$ and denoting the corresponding Green function to be $\tilde{G}_{x_i}$, then
\begin{align*}
\sup_{v\in T_{y_i}M, ~|v|=1}\inf _{x_{i}\in \tilde{B}_{R}(y_i)\cap \cC_i}\langle\nabla_y\tilde{G}_{x_i}(y_i),v\rangle\,<\, i^{-1}.
\end{align*}
To deduce a contradiction, we will show that the Green function $\tilde{G}_{x_i}$  converges to a function $D\,d^{2-n}_{x}$ on  $\mathbb{R}^{n-4}\times C(S^3/\Gamma)$ with constant $C^{-1}(n,\rv)<D<C(n,\rv)$.  Since the convergence is $C^1$ on the neck region due to the harmonic radius control (here we need bounded generalized Ricci), we can take $v_y$ to be any vector which approximates the radial direction on the $C(S^3/\Gamma)$ factor in order to conclude the result.

Thus, we only need to show the Green function $\tilde{G}_{x_i}\to D\,d^{2-n}_{x}$. On the one hand, we notice that $\tilde{B}_{\delta_i^{-1/2}}(x_i)$ converges to the same limit $\mathbb{R}^{n-4}\times C(S^3/\Gamma)$ for any sequence $x_{i}\in \tilde{B}_{R}(y_i)\cap \cC_i$.  On the other hand, by the established properties $(1)$ and $(2)$ we have on $\tilde{B}_{\delta_i^{-1/2}}(x_i)$, we have
$C^{-1}d_{x_i}^{2-n}\le \tilde{G}_{x_i}\le Cd_{x_i}^{2-n}$ and $|\nabla \tilde{G}_{x_i}|\le Cd_{x_i}^{1-n}$.  By Arzel\'a-Ascoli we have that $\tilde{G}_{x_i}$ converges to a function $G_x$ on the limit space $\mathbb{R}^{n-4}\times C(S^3/\Gamma)$ which satisfies the estimates
\begin{align}\label{e:green_cone_upper_lower}
C^{-1}d_{x}^{2-n}\le {G}_{x}\le Cd_{x}^{2-n},~~~\, \mbox{ and }\,~~~\, |\nabla G_x|\le Cd_x^{1-n}\, .
\end{align}
Since $\tilde G_{x_i}$ converges smoothly on the regular part of $\mathbb{R}^{n-4}\times C(S^3/\Gamma)$ we at least have that $G_x$ is harmonic away from the singular set.  If we lift $G_x$ to a function $G_{\tilde x}$ on $\dR^n$, we get away from $\tilde x$ that $G_{\tilde x}$ is locally lipschitz and harmonic away from a set of zero capacity.  Hence, $G_x$ is harmonic away from $\tilde x$ with the bounds \eqref{e:green_cone_upper_lower}.  Now the only harmonic functions on $\dR^n\setminus \tilde x$ with estimates \eqref{e:green_cone_upper_lower} are multiples of the Green's function.  Hence we have $G_x=D\, d_x^{2-n}$ for some constant $C^{-1}(n,\rv)<D<C(n,\rv)$ as claimed, which finishes the proof.
\end{proof}

We are now ready to prove Proposition \ref{p:green_function_distant_function}.
\begin{proof}[Proof of Proposition \ref{p:green_function_distant_function}]

Noting that $G_{\bar \mu}(y) = \int_{B_{39/20}(p)} G_x(y)\, d{\bar \mu}(x)$ we will use the pointwise Green function estimates in Lemma \ref{l:green_function_cone_point} and the Ahlfors regularity of Lemma \ref{l:Ahlfors}.  Indeed, for any $y\in \cN_{10^{-5}}\cap B_{39/20}(p)$ let $r\equiv d(y,\cC)$, and let us estimate $G_{\bar \mu}(y)$ above as follows:
\begin{align*}
G_{\bar \mu}(y) \le&\ C\int_{B_{39/20}(p)} d_x^{2-n}(y)d{\bar \mu}(x)\\
=&\ C\int_{B_{10r}(y)\cap B_{39/20}(p)}d_x^{2-n}(y)d{\bar \mu}(x)+C\sum_{i=1}^\infty\int_{A_{10^ir,10^{i+1}r}(y)\cap B_{39/20}(p)}d_x^{2-n}(y)d{\bar \mu}(x) \\
\le&\  Cr^{2-n}\int_{B_{10r}(y)\cap B_{39/20}(p)}d{\bar \mu}(x)+Cr^{2-n}\sum_{i=1}^\infty 10^{-i(n-2)}\int_{A_{10^ir,10^{i+1}r}(y)\cap B_{39/20}(p)}d{\bar \mu}(x)\, ,\notag\\
\leq&\ C\cdot B r^{-2}\Big(1+\sum 10^{-2i}\Big)\equiv C(n,\rv,B) r^{-2}\, ,
\end{align*}
as claimed.  To prove the lower bound of $G_{\bar \mu}$ we can similarly compute
\begin{align*}
G_{\bar \mu}(y)&\geq C^{-1}\int_{B_{39/20}(p)} d_x^{2-n}(y)d{\bar \mu}(x)\geq C^{-1}\int_{B_{10r}(y)\cap B_{39/20}(p)}d_x^{2-n}(y)d{\bar \mu}(x)\notag\\
&\geq C^{-1} 10^{2-n}r^{2-n}{\bar \mu}\Big(B_{8r}(x_0)\cap B_{39/20}(p)\Big)\geq C^{-1}(n,\rv,B)r^{-2},
\end{align*}
where $d(x_0,y)\le 2r=2d(y,\cC)$ and $x_0\in\cC$. 
Since $b^{-2}(y)=G_{\bar \mu}(y)$, we have $C^{-1}r\le b(y)\le Cr$, or that $C^{-1}\,d(y,\cC)\leq b(y)\leq C\,d(y,\cC)$. Hence we have proven (1) of Proposition \ref{p:green_function_distant_function}.

For the gradient estimate, using the same computational strategy as above we have
\begin{align*}
|\nabla G_{\bar \mu}(y)|=\int_{B_{39/20}(p)} |\nabla G_x(y)|d{\bar \mu}(x)\le C\int_{B_{39/20}(p)} d_x^{1-n}(y)d{\bar \mu}(x)\le Cr^{-3}.
\end{align*}
By noting that  $b^{-2}=G_{\bar \mu}$, then $2 b^{-3}|\nabla b|=|\nabla b^{-2}|=|\nabla G_{\bar \mu}|\le Cr^{-3}$. By the upper bound estimate of $b$, we have $|\nabla b|\le C$.  For the gradient lower bound,  for any fixed unit vector $v\in T_yM$, we have
\begin{align*}
\langle\nabla G_{\bar \mu}(y),v\rangle&=\int_{B_{39/20}(p)} \langle\nabla G_x(y),v\rangle d{\bar \mu}(x)=\int_{B_{rR}(y)}\langle\nabla G_x(y),v\rangle d{\bar \mu}(x)+\int_{M\setminus B_{rR}(y)}\langle \nabla G_x(y),v\rangle d{\bar \mu}(x).
\end{align*}
By the gradient upper bound estimate $|\nabla G_x(y)|\le Cd_x^{1-n}(y)$, we have
\begin{align*}
\left|\int_{M\setminus B_{rR}(y)}\langle\nabla G_x(y),v\rangle d{\bar \mu}(x) \right|\le \int_{M\setminus B_{rR}(y)}|\nabla G_x(y)|d{\bar \mu}(x)\le CR^{-3}r^{-3} .
\end{align*}
On the other hand, for fixed $R$ we have by the Green function estimates of Lemma \ref{l:green_function_cone_point} that if $\delta\leq \delta(R,n,\rv)$, then there is a unit vector $v_y\in T_yM$ such that
\begin{align*}
\int_{B_{rR}(y)}\langle \nabla G_x(y),v_y\rangle d{\bar \mu}(x)\geq \int_{B_{10r}(y)}r^{1-n}c(n)d{\bar \mu}(x)\geq C(n,B)r^{-3}.
\end{align*}
Therefore, combining the estimates above, we have
\begin{align*}
\langle\nabla G_{\bar \mu}(y),v_y\rangle\,\geq&\ \int_{B_{rR}(y)}\langle\nabla G_x(y),v_y\rangle d{\bar \mu}(x)-\brs{\int_{M\setminus B_{rR}(y)}\langle \nabla G_x(y),v_y\rangle d{\bar \mu}(x)}\\
\geq&\ C_0(n,B)r^{-3}-C_1(n,B)R^{-3}r^{-3}.
\end{align*}
Choosing $R=R(\rv,n,B)$ large enough we conclude $
\langle\nabla G_{\bar \mu}(y),v_y\rangle\,\geq\, \frac{C(B,n)}{2}r^{-3}$.  In particular, this gives us the estimate $|\nabla G_{\bar \mu}|(y)\geq C r^{-3}$ and hence the desired estimate $|\nabla b|(y)\geq C$ for $y\in \cN_{10^{-5}}\cap B_{39/20}(p)$.
\end{proof}

\vspace{.2cm}
\subsection{Reduction to $L^4$ estimate for $H$}

In this subsection we reduce the $L^2$ estimate for the curvature on the neck region to an $L^4$ estimate for $H$.  The main theorem of this subsection will be:
\begin{thm}\label{t:L2_neck}
	Let $(M^n,g,H,p)$ satisfy $\Vol(B_1(p))>\rv>0$ and $\brs{\RC}\le \delta$.  For any $\epsilon>0$, if  $\delta<\delta'(n,\epsilon,\rv)$ and $\cN = B_2(p)\setminus \bigcup_{x\in \cC} \overline B_{r_x}(x)$ is a $\delta$-neck region, then for each $B_{2r}(y)\subseteq B_2$ we have the curvature estimate \begin{equation}r^{4-n}\int_{\cN_{10^{-5}}\cap B_r(y)} \brs{\Rm}^2(z)\,dz \leq \epsilon+ r^{4-n}\int_{\cN_{10^{-5}}\cap B_{{r}}(y)} \brs{H}^4(z)\,dz.\end{equation} In particular, we have that 	\begin{equation}\int_{\cN_{10^{-5}}\cap B_1(p)} \brs{\Rm}^2(z)\,dz \leq \epsilon+\int_{\cN_{10^{-5}}\cap B_1(p)} \brs{H}^4(z)\,dz.\end{equation}
\end{thm}

The key result toward the proof of Theorem \ref{t:L2_neck} is the following local $L^2$ curvature estimate:
\begin{prop}\label{p:local_L2_neck}
Let $(M^n,g,H,p)$ satisfy $\Vol(B_1(p))>\rv>0$ and $|\RC|\leq\delta$. Also let  $\cN = B_2(p)\setminus \bigcup_{x\in \cC} \overline B_{r_x}(x)$ be a $\delta$-neck region.  Then there exist constants $K(n,B(n))>1, \delta'(n,\rv)$ such that if $y\in \cN_{K(n,B(n))}$ is a point with $d(y,\cC)=d(y,z)=r$ and $B_{r/10}(z)\subset B_2(p)$ for $z\in \cC$ and if  $\delta<\delta'(n,\rv)$, we have  
	\begin{align*}
		\int_{B_{r/4}(y) } \brs{\Rm}^2(z)\,dz \leq C(n,\rv,B) \left(\int_{B_{r/2}(y)}\brs{\Ric^g}^2\, dz+ \int_{B_{r/20}(z)}\big|{\cH}_{10r^2}(x)-{\cH}_{r^2/10}(x)\big|d\mu(x) \right).
	\end{align*}
\end{prop}
In turn, this proposition will rely on a key lemma which gives a pointwise scale-invariant estimate for the curvature tensor in terms of the Ricci curvature and a Gromov-Hausdorff approximation to a cone. Now let us record  Lemma 4.4  from \cite{JiNa}.

\begin{lemma}[Curvature estimate of Level sets]\label{l:curvature_level_sets} 
Let $(M^n,g,p)$ be a Riemannian manifold with  a map $\Phi=(h,f_1,\cdots, f_{n-4}): B_{10r}(p)\to \mathbb{R}_{+}\times \mathbb{R}^{n-4}$. Assume $f_0^2=h-\sum_{i=1}^{n-4}f_i^2\geq c_0r^2>0$ and $|f_i|\le c_0^{-1}r$ on $B_{2r}(y)\subset B_{10r}(p)$. Let $A=(a_{ij})$ be a $(n-3)\times (n-3)$ symmetric matrix with $a_{ij}(x)=\langle \nabla f_i(x),\nabla f_j(x)\rangle$.  Assume further $|\det A|(x)\geq c_0>0$ and $|\nabla f_i|\le c_1$ on $B_{2r}(y)$.  Then for any $x\in B_r(y)$, we have the following scale invariant estimate
\begin{equation}r^4\brs{\Rm}^2(x)\le C(n,c_0,c_1)\left(r^4\brs{\Ric^g}^2+r^{2}|\nabla^3 h|^2+\sum_{i=1}^{n-4}r^4|\nabla^3 f_i|^2+\mathcal{F}+\mathcal{F}^2\right),\end{equation}
where $\mathcal{F}=|\nabla^2h-2g|^2+\sum_{i=1}^{n-4}r^2|\nabla^2f_i|^2$.
\end{lemma}

By noting the curvature estimate above, to prove Proposition \ref{p:local_L2_neck}, we only
need to find $n-4$ functions that satisfy the condition in Lemma \ref{l:curvature_level_sets}. The main
purpose of the following subsections is to find and control such functions.

\subsection{Estimates of $\cH$-volume}

Let $(M,g,p)$ be a manifold with $\Ric^g\geq -(n-1)\delta$ and $\Vol(B_1(p))\geq \rv$. We recall a useful monotone quantity, the so called $\cH_t^\delta$-volume, which can be used to characterize closeness to a cone metric, as follows: 

\begin{align}\label{e:Hvolume_def}
\cH^\delta_t(p)=\int_M(4\pi t)^{-n/2}e^{-\frac{d^2(x,p)}{4t}}dx-\int_0^t\frac{1}{4s}\int_M( L_\delta(x)-2n)(4\pi s)^{-n/2}e^{-\frac{d^2(x,p)}{4s}}dx ds,
\end{align}
where $L_\delta(x)=2n+2(n-1)d(p,x)\sqrt{\delta}\coth \big(\sqrt{\delta}d(p,x)\big)$.  Note that $L_0\equiv 2n$ for spaces with nonnegative Ricci curvature, so that this second term is purely a correction term.  On a first reading of this section we recommend the reader sets $\delta=0$, most of the formulas simplify quite a bit in this case.  By direct computation, we have
\begin{align}
\partial_t\mathcal{H}^\delta_t(p)
&=\int_M\left( \frac{d^2(p,x)}{4t^2}-\frac{L_\delta(x)}{4t}\right)(4\pi t)^{-n/2}e^{-\frac{d^2(x,p)}{4t}}dx\, .
\end{align}
Noting that
$
\Delta e^{-\frac{d^2(x,p)}{4t}}=\left(-\frac{1}{4t}\Delta d^2(p,x)+\frac{d^2(p,x)}{4t^2}\right)e^{-\frac{d^2(x,p)}{4t}}\, ,
$
we have
\begin{align}\label{e:derivative_cHt}
\partial_t\mathcal{H}^\delta_t(p)
&=\frac{1}{4t}\int_M\big( \Delta d^2(p,x)-{L_\delta(x)}\big)(4\pi t)^{-n/2}e^{-\frac{d^2(x,p)}{4t}}dx\le 0\, ,
\end{align}
where we use the Laplacian comparison in the last inequality.  
\begin{rmk}\label{r:scaling_cHvolume}
By rescaling, one can check that $\cH^{\delta}_t(p)=\tilde{\cH}^{\delta t}_1(p)$ where $\tilde{\cH}$ is the $\cH$-volume of $(M,\tilde{g},p)=(M,t^{-1}g,p)$. 
\end{rmk}

For simplicity of notation, we will drop the $\delta$ of $\cH_t^\delta$ when there is no confusion, but one should keep in mind the dependence of $\cH_t$ on the lower Ricci curvature bound $-(n-1)\delta$.  Let us consider the following heat flow
\begin{gather}\label{e:heat_flow}
\begin{split}
&\partial_t f_t=\Delta f_t-2n,\\
&f_0(x)=2nU(d(p,x))\, ,
\end{split}
\end{gather}
where $U(r)=\int_0^r\sinh^{-(n-1)}(\sqrt{\delta} t)\int_0^t\sinh^{n-1}(\sqrt{\delta} s)ds dt.$

\begin{rmk}\label{r:deltaf0_delta_d2}
By direct computations, we have 
\begin{align}
0\le 2n-\Delta f_0(x)=\frac{nU'(d(p,x))}{d(p,x)}\Big(L_\delta(x)-\Delta d^2(p,x) \Big)\le n\Big(L_\delta(x)-\Delta d^2(p,x)\Big).
\end{align}
\end{rmk}

We begin by recording some basic points about $f_t$ from Lemma 4.21 of \cite{JiNa}  which will be useful:

\begin{lemma}\label{l:heat_flow_estimate}
\cite{JiNa} 
	 If $(M,g)$ satisfies $\Ric^g \geq -(n-1)\delta$ and $f_t$ is as in \eqref{e:heat_flow} then we have the following:
\begin{enumerate}
\item $\Delta f_t\le 2n$ for all $t\geq 0$.
\item $-2nt\leq f_t\leq f_0 = 2nU(d(p,x))$.
\item If $Q_t\equiv e^{-2(n-1)\delta t}(|\nabla f|^2-4f-8nt)-8(n-1)\delta (tf+2nt^2)$, then $Q_t\leq 0$ for all $t\geq 0$.
\item If $P_t(x)=e^{-2(n-1)\delta t}(|\nabla f|^2-4f+8t(\Delta f-2n))-8(n-1)\delta (tf+2nt^2)+8t^2\delta(n-1)(\Delta f-2n)$ then we have $P_t\leq 0$ for all $t\geq 0$.
\end{enumerate}
\end{lemma}

\begin{proof}

Since $(\partial_t-\Delta)(f+2nt)=0$ we have 
\begin{align}
f_t(x)+2nt=\int_M f_0(y)\rho_t(x,dy).
\end{align}
This implies that 
\begin{align}
\Delta f_t=\partial_t f_t+2n=\int_Mf_0(y)\partial_t\rho_t(x,dy)=\int_M \Delta f_0(y)\rho_t(x,dy).
\end{align}
Therefore we arrive at 
\begin{align}\label{e:deltaft-2n}
(2n-\Delta f_t)(x)=\int_M (2n-\Delta f_0)(y)\rho_t(x,dy)\, .
\end{align}
By Laplacian comparison, we have $\Delta U(d(p,x))\le 1$. Thus $(2n-\Delta f_0)\geq 0$. Hence \eqref{e:deltaft-2n} implies (1).  From this we immediately get $\partial_tf\le 0$ which yields the upper bound in $(2)$.  For the lower bound, we use $(\partial_t-\Delta)(f+2nt)=0$ to get $f_t(x)\geq -2nt$.

For the gradient estimate, we have by direct computation that
\begin{align*}
(\partial_t -\Delta)(|\nabla f|^2-4f)&=-2|\nabla^2f|^2-2\Ric^g(\nabla f,\nabla f)+8n\le 2(n-1)\delta |\nabla f|^2+8n.
\end{align*}
By noting our lower bound $f_t+2nt\geq 0$ we can conclude from this
\begin{align*}
	(\partial_t-\Delta)Q_t(x)\le 0\, .
\end{align*}
Combined with $Q_0\le 0$ we obtain $(3)$.  Finally, to prove $(4)$ we compute that
\begin{align}\label{e:evolution_Pt}
(\partial_t-\Delta)P_t\le -2e^{-2(n-1)\delta t}|\nabla^2 f-2g|^2\, .
\end{align}
Combined with $P_0\leq 0$ this proves the desired result.
\end{proof}

Let us now prove a couple more refined estimates on $f_t$ that depend on the pinching of our $\cH$-volume.  Precisely, we have the following:

\begin{lemma}\label{l:cHfunctional}
Let $(M,g,p)$ satisfy $\Vol(B_1(p))\geq \rv>0$ and $\Ric^g \geq -(n-1)\delta$. Denote by $\eta(r)=|\cH_{r^2/2}(p)-\cH_{2r^2}(p)|$ the $\cH$-volume pinching at scale $r$.  Then we have the estimates:
\begin{enumerate}
\item $\fint_{0}^{r^2} \fint_{B_{4r}(p)}|\nabla^2 f_t-2g|^2dydt\le C(n,\rv)(\delta r^2+\eta).$
\item For any $t\le r^2$, we have $\sup_{y\in B_{4r}(p)}|f_t-f_0|(y)\le C(n,\rv)\epsilon(\eta)\, r^2,$ where $\epsilon(\eta)\to 0$ if $\eta\to 0$.
\item For any $t\le r^2$, we have $\fint_{B_{4r}(p)}|\Delta f_t-2n|(z)\,dz+r^{-2}\fint_{B_{4r}(p)}|\nabla f_0-\nabla f_t|^2(z)\,dz\le C(n,\rv)\eta.$
\end{enumerate}
Moreover, if we assume further $|\RC|\leq \delta$ and $r_h(x)\geq r$ for some $x\in B_{3r}(p)$, then for any $r^2/2\le t\le r^2$, we have
\begin{align}\label{e:Hessian_parabolic_pointwise}
\sup_{B_{r/2}(x)}|\nabla^2 f_t-2g|^2\le C(n,\rv)(\delta r^2+\eta+r^2\sup_{B_{r}(x)}\brs{\Ric^g}).
\end{align}
Finally, if $\phi$ denotes a Cheeger-Colding cutoff function as in \eqref{e:cutoff}, then there exists $t\in [r^2/2, r^2]$ such that
\begin{gather}\label{e:hessian_parabolic_integral}
\begin{split}
\fint_{B_{r}(x)}|\nabla^2 f_t-2g|^2+&r^2\fint_{B_{r}(x)}\phi|\nabla^3 f_t|^2(z)\,dz\\
\le&\ C(n,\rv)\Big(\delta r^2+\eta+r^4\fint_{B_{r}(x)}(\brs{\Ric^g}^2+\phi\frac{1}{10}\brs{\Rm}^2)\Big).
\end{split}
\end{gather}
\end{lemma}
\begin{proof}
By \eqref{e:deltaft-2n} we have
\begin{align*}
(2n-\Delta f_t)(x)=\int_M (2n-\Delta f_0)(y)\rho_t(x,dy)\, .
\end{align*}
Thus
\begin{align*}
\int_M (2n-\Delta f_t)(x)\rho_{r^2/2}(p,dx)=&\ \int_M\int_M(2n-\Delta f_0)(y)\rho_t(x,dy)\rho_{r^2/2}(p,dx)\\
=&\ \int_M (2n-\Delta f_0)(y)\rho_{t+r^2/2}(p,dy)\, .
\end{align*}
Therefore,  by the heat kernel estimates of Theorem \ref{t:heat_kernel} and Remark \ref{r:deltaf0_delta_d2} and noting that $\Delta f_t\le 2n$, for any $t\le r^2$, one can get
\begin{gather*}
\begin{split}
\fint_{B_{10r}(p)}|2n-\Delta f_t|(y)dy&\le C(n,\rv)\int_M (2n-\Delta f_t)(y)\rho_{r^2/2}(p,dy)\\
&\le C(n,\rv)\int_M (2n-\Delta f_0)(y)\rho_{t+r^2/2}(p,dy)\\ 
&\le C(n,\rv)\int_M \big( {L_\delta(y)}-\Delta d^2(p,y)\big)\rho_{t+r^2/2}(p,dy)\\
&\le C(n,\rv)\int_M\big( {L_\delta(y)}-\Delta d^2(p,y)\big) \Big(t+r^2/2\Big)^{-n/2}e^{-\frac{d(p,y)^2}{(4+2^{-1})(t+r^2/2)}}dy\\
&\le  C(n,\rv)\int_M\big( {L_\delta(y)}-\Delta d^2(p,y)\big) r^{-n}e^{-\frac{4d(p,y)^2}{27r^2}}dy\\ \label{e:meanvalue_cH7}
&\le C(n,\rv) |\cH_{7r^2/4}(p)-\cH_{2r^2}(p)|\\
&\le C(n,\rv)\eta\, ,
\end{split}
\end{gather*}
where the constants $C(n,\rv)$ change line by line, and we have used \eqref{e:derivative_cHt} and mean value equality to deduce the third line from bottom. 

Hence, we have
\begin{align}\label{e:L1f0_ft}
\fint_{B_{10r}(p)}|f_0-f_t|(x)dx\le \int_0^t\fint_{B_{10r}(p)}(2n-\Delta f_s)(x)dxds\le C(n,\rv)\eta\, t\, .
\end{align}
We also note that a consequence of part (3) of Lemma \ref{l:heat_flow_estimate} we have
\begin{gather*}
\begin{split}
|\nabla f|^2&\le \big(4+8(n-1)\delta\,t\, e^{2(n-1)\delta t}\big)\big(f+2nt\big)\\
&\le \big(4+8(n-1)\delta\,t\, e^{2(n-1)\delta t}\big)\big(f_0+2nt\big)\, .
\end{split}
\end{gather*}
Using this, we therefore have $\sup_{B_{10r}(p)}|f_0-f_t|\le C(n,\rv)\epsilon(\eta)\, r^2$.

Now let $\phi$ be a cutoff function as in \cite{ChC1} with support in $B_{10r}(p)$ and $\phi\equiv 1$ on $B_{8r}(p)$ and $r|\nabla \phi|+r^2|\Delta \phi|\le C(n)$. By integrating by parts, we have
\begin{align*}
\int_M\phi^2 |\nabla f_0-\nabla f_t|^2(x)\,dx \le&\  2\int_M|\nabla\phi| \cdot |f_0-f_t|\cdot \phi|\nabla f_0-\nabla f_t|(x)\,dx+\int_M \phi^2 |\Delta f_0-\Delta f_t|\cdot |f_0-f_t|(x)\,dx\\ \nonumber
\le&\ \frac{1}{2}\int_M\phi^2 |\nabla f_0-\nabla f_t|^2(x)\,dx+C(n)r^{-2}\int_{B_{10r}(p)}|f_0-f_t|^2\\
&\ + \sup_{B_{10r}(p)}|f_0-f_t|\int_{B_{10r}(p)}|\Delta f_0-\Delta f_t|(x)\,dx.
\end{align*}
Combining with \eqref{e:L1f0_ft} and $\sup_{B_{10r}(p)}|f_0-f_t|\le C(n,\rv)\epsilon(\eta)\, r^2$ and the $L^1$ estimates of $|\Delta f_t-2n|$, we have $\int_M\phi^2 |\nabla f_0-\nabla f_t|^2(z)\,dz\le C(n,\rv)\eta r^{2+n}$. Hence, we prove (2) and (3).

To prove (1), recall that $P_t(x)=e^{-2(n-1)\delta t}(|\nabla f|^2-4f+8t(\Delta f-2n))-8(n-1)\delta (tf+2nt^2)+8t^2\delta(n-1)(\Delta f-2n)$ satisfies $P_t\leq 0$ from Lemma \ref{l:heat_flow_estimate}.  Let $\varphi$ be a cutoff function as in \eqref{e:cutoff} with support in $B_{6r}(p)$ and $\varphi\equiv 1$ on $B_{5r}(p))$ and $r|\nabla \varphi|+r^2|\Delta \varphi|\le C(n)$. Multiplying equation (\ref{e:evolution_Pt}) by $\varphi$ and integrating by parts we have
\begin{equation}
\int_0^{r^2}\fint_{B_{5r}(p)} |\nabla^2 f-2g|^2(z)\,dzdt\le C(n,\rv)\fint_{B_{6r}(p)} |P_{r^2}|(z)\,dz+C(n,\rv)\fint_0^{r^2}\fint_{B_{6r}(p)}|P_t|(z)\,dzdt.
\end{equation}
To estimate $\fint_{B_{6r}(p)}|P_t|(z)\,dz $, one only needs to estimate $\fint_{B_{6r}(p)} \Big||\nabla f|^2-4f\Big|(z)\,dz$. This can be controlled by considering the evolution of $f_t^2-f_0^2$. In fact, we have
\begin{align*}
(\partial_t-\Delta)(f_0^2-f^2)
&=(4n+8)(f-f_0)-2f_0(\Delta f_0-2n)-2(|\nabla f_0|^2-4f_0)+2(|\nabla f|^2-4f).
\end{align*}
Let $\psi$ be a cutoff function as in \eqref{e:cutoff} with support in $B_{10r}(p)$ and $\psi\equiv 1$ on $B_{6r}(p))$ and $r|\nabla \psi|+r^2|\Delta \psi|\le C(n)$. By noting $0\le (4f_0-|\nabla f_0|^2)(y)\le C(n)\delta d(p,y)^4$  for $d(p,y)\le 10$, we can show
\begin{align}
\int_0^{3r^2/2}\int_{M}e^{-2(n-1)\delta t}(|\nabla f|^2-4f)\psi (z)\,dzdt \geq -C(n,\rv)(\delta r^2+\eta )r^{4+n} .
\end{align}
Using the above, the $L^1$ estimate on the Laplacian of $f_t$ and $|f_t|\le C(n)r^2$ on $B_{10r}(p)$ for $t\le 3r^2/2$, we have
\begin{align}\nonumber
\int_0^{3r^2/2}\int_{M}P_t\psi(z)\,dzdt\geq &-C(n,\rv)(\delta r^2+\eta )r^{4+n}\, .
\end{align}
Noting that $P_t\le 0$, we have
\begin{align}
\fint_0^{3r^2/2}\fint_{B_{6r}(p)}|P_t|(z)\,dzdt\le C(n,\rv)(\delta r^2+\eta )r^2\, ,
\end{align}
which finishes the proof of (1).

Now we wish to prove \eqref{e:Hessian_parabolic_pointwise}.  Indeed, under the assumption $r_h(x)\geq r$ we have that
\begin{align}\label{e:moser_iteration:1}
	r^{2q}\fint_{B_{3r/4}(x)}\brs{\Rm}^q(z)\,dz<C(n,q)\, ,
\end{align}
for all $q<\infty$.  Using this, the estimates of \eqref{e:Hessian_parabolic_pointwise} are fairly standard, and follow much the same path as the proof of Proposition \ref{p:prelim:harmonic_estimates} below, so we will only sketch the argument. Denote $H_f=\nabla^2f-2g$,  then we can compute
\begin{align}\label{e:evolution_Hf}
(\partial_t-\Delta)|H_f|^2=-2|\nabla^3f|^2+ \Rm\ast H_f\ast H_f+\Ric^g\ast H_f+\nabla(\Ric^g(\nabla f,\cdot))\ast H_f\, ,
\end{align}
where $\ast$ means tensorial linear combinations and the exact expression can be computed as in \cite{JiNa}.  
Then we can apply a standard parabolic Moser iteration using \eqref{e:moser_iteration:1} and (1)  to conclude the pointwise estimate in (\ref{e:Hessian_parabolic_pointwise}). Now we address the estimate (\ref{e:hessian_parabolic_integral}), firstly the $L^2$ bound of $|H_f|$ is consequence of $(1)$. In order to conclude the $L^2$ estimate on $\nabla^3 f$, we simply multiply \eqref{e:evolution_Hf} by a standard cutoff function $\phi$ and integrate by parts, which shows that 
\begin{align*}
\int_{t_1}^{t_2}\int_{B_{r}(x)}\phi|\nabla^3 f_t|^2(z)\,dzdt \le&\ C(n,\rv)r^{-2}\int_{t_1}^{t_2}\int_{B_r(x)}|H_f|^2(z)\, dzdt+\frac{1}{2}\int_{t_1}^{t_2}\int_{B_{r}(x)}\phi|\nabla^3 f_t|^2(z)\,dzdt\\
&\ \ \qquad +r^2\int_{t_1}^{t_2}\int_{B_{r}(x)}|\Ric^g|^2(z)\,dzdt+\frac{1}{10}r^2\int_{t_1}^{t_2}\int_{B_{r}(x)}\phi\brs{\Rm}^2(z)\,dzdt\\
&\ \ \qquad +\int_{B_{r}(x)}\phi|H_{f_{t_1}}|^2(z)\,dz+\int_{B_{r}(x)}\phi|H_{f_{t_2}}|^2(z)\,dz.
\end{align*}
where we have used the fact that $|H_f|$ is bounded by $C(n,\rv)$. Using the $L^2$ bound of  $|H_f|$ just proved and the above formula will finish the proof.
\end{proof}
\subsection{$\cH$ volume and $L^2$ curvature estimate}
The main purpose of this
subsection is to prove the local $L^2$ curvature estimate in Proposition \ref{p:local_L2_neck}. The
key ingredient is the parabolic estimate of $\cH$-volume in Lemma \ref{l:cHfunctional} and the
pointwise curvature estimate in Lemma \ref{l:curvature_level_sets}.  First, let us recall the concept
of independent points.

\begin{defn}[$(k,\rho,r)$-independent points]\label{d:indep_point}
In a metric space $(X,d)$ a set of points $U=\{x_0,\cdots,x_k\}\subset B_{r}(x)$ is \emph{$(k,\rho,r)$-independent} if for any subset $U'=\{x_0',\cdots,x_k'\}\subset \mathbb{R}^{k-1}$ we have 
\begin{align}
d_{GH}(U,U')\geq \rho \cdot r.
\end{align}
\end{defn}

\begin{rmk}\label{r:remark_independentpoints}
Let $X\subset \mathbb{R}^n$, if there exists no $(k,\rho,r)$-independent set in $B_r(x)\cap X$, then $B_r(x)\cap X\subset B_{4\rho r}(\mathbb{R}^{k-1})$ for some $(k-1)$-plane $\mathbb{R}^{k-1}\subset \mathbb{R}^n$. To see this, if $B_r(x)\cap X$ is not a subset of $B_{3\rho r}(\mathbb{R}^{k-1})$ for any $(k-1)$-plane, then one can find $(k,\rho,r)$-independent set in $B_r(x)\cap X$ by induction on $k$.
\end{rmk}
\begin{proof}[Proof of Proposition \ref{p:local_L2_neck}]
The main idea for the proof is to use the pointwise curvature estimate in Lemma \ref{l:curvature_level_sets}. The key ingredient is to find $n-3$ functions $(h,u_1,\ldots,u_{n-4})$ which satisfy the conditions of this lemma.  Intuitively, such $h$ should be a square distance function and $(u_1,\ldots,u_{n-4})$ splitting functions which form a cone map to $\mathbb{R}^{n-4}\times C(S^3/\Gamma)$.  We will construct such functions in detail in the following paragraphs. First we will show the following claim:

\vskip 0.1in

\noindent \textbf{Claim 1:} For $\delta\le \delta'(n,B(n),\tau(n))$ there exist $\rho(n,B), L(n,B), K(n,B)>0$ such that for any $y\in \cN_{K}$ with $d(y,\cC)=d(y,z)=r$ satisfying $B_{r/10}(z)\subset B_2(p)$,  we have $(n-4,\rho,r/20)$-independent points $\{x_0,\cdots, x_{n-4}\}\subset \tilde{\cC}\subset B_{r/20}(z)$, where
\begin{align*}
\tilde{\cC}=\{x\in \cC\cap B_{r/20}(z): \big|{\cH}_{10r^2}(x)-{\cH}_{r^2/10}(x)\big|\le L(n,B) \eta\}
\end{align*}
 and $\eta = \fint_{B_{r/20}(z)}\big|{\cH}_{10r^2}(x)-{\cH}_{r^2/10}(x)\big|d\mu(x)$.

\vskip 0.1in

\textbf{Proof of Claim 1:} We will choose $K(n,B)=10^{10}\rho(n,B)^{-1}>10^{10}$ where $\rho(n,B)$ will be fixed later. 
By a maximal function argument, for any $\epsilon>0$, there exists $L= L(n,\epsilon)>0$ such that the set 
\begin{align*}
\tilde{\cC}=\{x\in \cC\cap B_{r/20}(z): \big|{\cH}_{10r^2}(x)-{\cH}_{r^2/10}(x)\big|\le L(n,\epsilon) \eta\}
\end{align*}
satisfies 
\begin{align}\label{e:lowerbouncC}
\mu(\tilde{\cC})\geq (1-\epsilon)\mu(B_{r/20}(z))\geq (1-\epsilon) B^{-1}r^{n-4}20^{-(n-4)}\geq C(n)B^{-1}r^{n-4}.
\end{align}
From the definition of $(n-4,\rho,r)$-independent points and Remark \ref{r:remark_independentpoints} , it suffices to prove the following: 
\vskip 2mm
If $\rho=\rho(n,B):=\hat{C}(n)^{-1}B^{-2}$ then for all $k$ planes $P_k\subset \mathbb{R}^{n-4}\times \{y_c\}\subset \mathbb{R}^{n-4}\times C(S^3/\Gamma)$  with $k\le n-5$, one has
\begin{align}\label{e:PkcoveringtildecC}
B_{5\rho r}(\iota(P_k))\setminus \tilde{\cC}\ne\emptyset,
\end{align}
where $\iota$ is a $\delta r$-GH map from $\mathbb{R}^{n-4}\times C(S^3/\Gamma)$.  
\vskip 2mm
To see \eqref{e:PkcoveringtildecC}, note that $\iota(P_k)\cap B_{r/20}(z)$ can be covered by $C(n)\rho^{-k}$ many $B_{\rho r}$-balls, which implies that 
\begin{gather}\label{e:PkcoveringtildecC2}
\begin{split}
\mu(B_{5\rho r}(\iota(P_k))\cap B_{r/20}(z))\le&\ C(n)\rho^{-k} \cdot C(n)B (\rho r)^{n-4}\\
\le&\ C(n)B\rho^{n-4-k}r^{n-4}\\
\le&\ C(n)B\rho r^{n-4},
\end{split}
\end{gather}
where we have used the fact that for each $B_{\rho r}(w)$ with $w\in B_{r/20}(z)$ and $\rho<1/40$ that 
\begin{align}\label{e:upperpackingrhow}
\mu(B_{6\rho r}(w))\le C(n)B (\rho r)^{n-4}.
\end{align}  
Actually, to see \eqref{e:upperpackingrhow} we first note that $r\geq  K(n,B)r_z=10^{10}\rho(n,B)^{-1}r_z$ by the choice of $z$.  Since $|\text{Lip}\ r_x|\le \delta$, for any $x\in \cC\cap B_{r}(z)$ we have $r_x\le 2r$. Thus for any $w\in B_{r/20}(z)$ if $B_{6\rho r}(w)\cap \cC=\emptyset$ then \eqref{e:upperpackingrhow} holds trivially. If $B_{6\rho r}(w)\cap \cC\ne \emptyset$ then there exists $x\in B_{6\rho r}(w)\cap \cC\subset B_r(z)\cap \cC$, hence $\mu(B_{6\rho r}(w))\le \mu(B_{12\rho r}(x))\le  B (12 \rho r)^{n-4}$ which proves \eqref{e:upperpackingrhow}.  The estimates \eqref{e:lowerbouncC} and \eqref{e:PkcoveringtildecC2} imply \eqref{e:PkcoveringtildecC} if $\rho(n,B)= \hat{C}(n)^{-1}B^{-2}$ for a large $\hat{C}(n)$. Therefore, we have proved Claim 1. $\square$ \\

Now recall that $\eta=\fint_{B_{r/20}(z)}\big|{\cH}_{10r^2}(x)-{\cH}_{r^2/10}(x)\big|d\mu(x)$.  Applying Lemma \ref{l:cHfunctional} to each $x_i$, we have $n-3$ functions $f_{i,t}$ such that  $\fint_{0}^{r^2} \fint_{B_{4r}(x_i)}|\nabla^2 f_{i,t}-2g|^2dxdt\le C(n,\rv,B)(\delta r^2+\eta).$  For $1\le i\le n-4$, let \begin{equation}w_{i,t}= \left(f_{i,t}-f_{0,t}-d(x_0,x_i)^2\right)/2d(x_0,x_i).  
\end{equation}

\vskip 0.1in

\noindent \textbf{Claim 2:} There exists an $(n-4,n-4)$-matrix $D$ with $|D|\le C(n,B,\tau)$ such that if $(v_{i,r^2})\equiv (w_{i,r^2})D$ then $\bar f_{0,r^2}\equiv f_{0,1}-\sum_{i=1}^{n-4} v_{i,r^2}^2$ satisfies $\bar f_{0,r^2}\geq C(n,B,\tau)r^2>0$ on $B_{5r/6}(y)$. Further, if we denote $v_{0,r^2}=\sqrt{\bar f_{0,r^2}}$ on $B_{5r/6}(y)$, then we have $\min_{x\in B_{r/2}(y)}|\det A|(x)\geq 1/2$, where $A(x)=\langle \nabla v_{i,r^2},\nabla v_{j,r^2}\rangle(x)$ for $i,j=0,\cdots, n-4$.\\

 \textbf{Proof of Claim 2:}  We prove this claim by contradiction, therefore let us assume this is not true. Then for $\delta_a\to 0, a\in \mathbb N$, there exists $\delta_a$-neck region $\cN_a\subset B_2(p_a)\subset M_a$ with $y_a\in \cN_{a,K(n,B)}$ (cf. Definition \ref{d:neck}), $d(y_a,\cC_a)=r_a$ and $(n-4,\rho, r_a/20)$-independent points $\{x_{0,a},\cdots,x_{n-4,a}\}\subset {\cC}_a$, but there is no matrix $D_a$ with $|D_a|\le C(n,B,\tau)$ satisfying the claim for $\cN_a$, where $C(n,B,\tau)$ will be determined later.  Let us rescale each metric $g_a$ to $\tilde{g}_a$ such that $d(y_a,\cC_a)=1$. Taking the limit we have $M_a\to \mathbb{R}^{n-4}\times C(S^3/\Gamma)$ and $\cC_a\to \cC_\infty\subset \mathbb{R}^{n-4}\times \{y_c\}$ with $y_a\to y_\infty$ and $x_{i,a}\to x_{i,\infty}$, where $\{x_{i,\infty}\}$ is a $(n-4,\rho,1/20)$-independent set.   By the $C^0$ estimate of Lemma \ref{l:cHfunctional} we have that $|\tilde{f}_{i,1,a}-d_{x_{i,a}}^2|\to 0$, and hence $\tilde{f}_{i,1,a}\to \tilde{f}_{i,1,\infty}\equiv d^2_{x_{i,\infty}}$.
 On the one hand, it is then clear that $\tilde{w}_{i,1,\infty}\equiv\left(\tilde{f}_{i,1,\infty}-\tilde{f}_{0,1,\infty}-d(x_{0,\infty},x_{i,\infty})^2\right)/2d(x_{0,\infty},x_{i,\infty})$ is a linear function. As $\{x_{i,\infty}\}$ is a $(n-4,\rho,1/20)$-independent set, one can choose a matrix $D_\infty$ with $|D_\infty|\le C(n,\rho)\equiv C(n,B,\tau)$ such that $(\tilde{v}_{i,1,\infty})\equiv(\tilde{w}_{i,1,\infty})D_\infty$ represents the standard coordinate functions of $\mathbb{R}^{n-4}\times \{y_c\}\subset \mathbb{R}^{n-4}\times C(S^3/\Gamma)$ with $(0^{n-4},y_c)=x_{0,\infty}$. Then $\bar f_{0,1,\infty}\equiv\tilde{f}_{0,1,\infty}-\sum_{i=1}^{n-4}\tilde{v}_{i,1,\infty}^2$ is the distance square function $d_{\mathbb{R}^{n-4}\times \{y_c\}}^2$. Thus the $(n-3,n-3)$ matrix $A_\infty$ defined in $B_1(y_\infty)\subset \mathbb{R}^{n-4}\times C(S^3/\Gamma)$ satisfies $A_\infty=(\delta_{ij})$. However, for the rescaled metric $\tilde{g}_a$ we claim that
 $\tilde{f}_{i,1,a}\to \tilde{f}_{i,1,\infty}=d^2_{x_{i,\infty}}$ in $C^2$ sense on $B_{2/3}(y_\infty)$.  First, by the harmonic radius lower bound of $y_a$ in Lemma \ref{l:neck:distance_harm_rad}, we have  $\tilde g_a\to g_\infty$.  This implies that $\eta_a$ converges to $\eta_\infty$ and $\eta_\infty=0$, because the limit metric is a cone metric. Secondly by Lemma \ref{l:neck:Hsmall}, we can bound $\brs{\Ric^{\til{g}_a}}$ by $\epsilon'(\delta_a)$ on $B_{{3/4}}(y_a)$,
 hence $\tilde f_{i,1,a}\to \tilde f_{i,1,\infty}=d^2_{x_i,\infty} $ in $C^2$ sense when $\delta_a \to 0$  by Lemma \ref{l:cHfunctional}.  By choosing $D_a=D_\infty$ for $a$ sufficiently large, this derives our contradiction and proves Claim 2. $\square$

\vskip 0.1in

Now we plan to use such functions to prove our expected curvature estimates by using Lemma \ref{l:curvature_level_sets}. 
Using the harmonic radius lower bound $r_h(y)\geq c(n,\rv)r$ and the estimate (\ref{e:hessian_parabolic_integral}) of Lemma \ref{l:cHfunctional}, we have
\begin{gather}\label{hh}
\begin{split}
\fint_{B_{r/2}(y)}|\nabla^2h-2g|^2& +r^2\fint_{B_{r/2}(y)}\phi|\nabla^3h|^2(z)\,dz\\
\le&\ C(n,\rv,B) \left(\eta+\delta r^2+r^4\fint_{B_{r/2}(y)}\left(\brs{\Ric^g}^2+\frac{1}{10}\phi\brs{\Rm}^2\right)\,dz \right)\,.
\end{split}
\end{gather}

By Claim 2, let us consider the map $$\Phi=(h,u_1,\cdots, u_{n-4})\equiv (f_{0,r^2},v_{1,r^2}, \cdots, v_{n-4,r^2})$$ on $B_{10r}(y)$. There is a similar estimate for $u_i:=v_{i,r^2}$.  Since the constant matrix $D$ satisfies $|D|<C(n,B,\tau)$, we can assume $v_{i,r^2}=w_{i,r^2}$ and by definition $w_{i,r^2}= \left(f_{i,r^2}-f_{0,r^2}-d(x_0,x_i)^2\right)/2d(x_0,x_i)$. Notice by choice $d(x_0,x_1)\geq \rho(n,B)r$, then direct calculation shows that
\begin{align*}
\nabla^2 w_{i,r^2}=\nabla^2 (f_{i,r^2}-f_{0,r^2})/2d(x_0,x_i)\le (\nabla^2 f_{i,r^2}-2g)r^{-1}+(2g-\nabla f_{0,r^2})r^{-1}.
\end{align*}
By this we conclude using Lemma \ref{l:cHfunctional} (specifically (\ref{e:hessian_parabolic_integral})) again
\begin{gather}\label{a}
\begin{split}
& \sum_{i=1}^{n-4}\left(\fint_{B_{r/2}(y)}r^2|\nabla^2u_i|^2+r^4\fint_{B_{r/2}(y)}\phi|\nabla^3u_i|^2(z)\,dz\right)\\ 
& \qquad \le C(n,\rv,B)\left(\eta+\delta r^2+r^4\fint_{B_{r/2}(y)}\Big(\brs{\Ric^g}^2+\frac{1}{10}\phi \brs{\Rm}^2\Big)(z)\,dz\right)\,.
\end{split}
\end{gather}

Moreover, by the pointwise nondegeneration of $A(x)$ in Claim 2 and Lemma \ref{l:curvature_level_sets}, for any $z\in B_{r/2}(y)$, we have scale invariant estimates
\begin{align*}
r^4\brs{\Rm}^2(z)\le C(n,\rv,B)\left(\brs{\Ric^g}^2r^4+r^2|\nabla^3h|^2+\sum_{i=1}^{n-4}r^4|\nabla^3 u_i|^2+\mathcal{F}+\mathcal{F}^2\right)(z)\, ,
\end{align*}
where $\mathcal{F}=|\nabla^2h-2g|^2+\sum_{i=1}^{n-4}r^2|\nabla^2u_i|^2$. By the pointwise Hessian estimates for $u_i,$ and the estimate of $\nabla^2h-2g$ from Lemma \ref{l:cHfunctional}, we have $\mathcal F^2<\mathcal F$, then  for any $z\in B_{r/2}(y)$, we have
\begin{align*}
r^4\brs{\Rm}^2(z)\le C(n,\rv,B)\left(\brs{\Ric^g}^2r^4+r^2|\nabla^3h|^2+\sum_{i=1}^{n-4}r^4|\nabla^3 u_i|^2+|\nabla^2h-2g|^2+\sum_{i=1}^{n-4}r^2|\nabla^2u_i|^2\right)(z)\, .
\end{align*}
Integrating $\phi\brs{\Rm}^2$over $B_{r/2}(y)$, by (\ref{hh}) and (\ref{a}) we get 
\begin{align*}
r^4\fint_{B_{r/2}(y)}\phi\brs{\Rm}^2(z)\,dz&\le C(n,\rv,B)(\eta +\delta r^2+r^4\fint_{B_{r/2}(y)}\brs{\Ric^g}^2(z)\,dz)\\&\le C(n,\rv,B)\Big(\delta r^2+\fint_{B_{r/20}(z)}\big|{\cH}_{10r^2}(x)-{\cH}_{r^2/10}(x)\big|d\mu(x)+r^4\fint_{B_{r/2}(y)}\brs{\Ric^g}^2(z)\,dz\Big)\, .
\end{align*}
where the term $\frac{1}{10}r^4\fint_{B_{r/2}(y)}\phi\brs{\Rm}^2(z)\,dz$ in (\ref{hh}) and (\ref{a}) is absorbed by the left hand side. 
This completes the proof.
\end{proof}

Now, using the local $L^2$ curvature estimate in Proposition \ref{p:local_L2_neck}, we prove Theorem \ref{t:L2_neck}.

\begin{proof}[Proof of Theorem \ref{t:L2_neck}]
Since the estimates are scale invariant, without loss of generality we will assume $r=1$. By the local $L^2$ curvature estimate of Proposition \ref{p:local_L2_neck}, for any $y\in\cN_{K(n,B)}$ with $2s=d(y,\cC)=d(y,z)$ and $B_{s/10}(z)\subset B_2(p)$, which holds in particular for all $y\in B_1(p)\cap \cN_{K(n,B)}$,  we have the estimate
\begin{align*}
\int_{B_{s/2}(y)} \brs{\Rm}^2(z)\,dz \leq C(n,\rv,B)\int_{B_{s}(y)} \brs{H}^4(z)\,dz+C(n,\rv,B)\int_{B_{s/20}(z)}|{\cH}_{10s^2}-{\cH}_{10^{-1}s^2}|(x)\,d\mu(x)\, .	
\end{align*}

In order to use such an estimate, we first construct a Vitali covering. For any $x\in \cN_{K(n,B)}$ with $d(x,\cC)=s_x$, consider the covering $\{B_{s_x/5}(x),x\in \cN_{K(n,B)}\}$ of $\cN_{K(n,B)}$. We can choose a subcovering $\{B_{s_a/5}(x_a)\}$ such that $\{B_{s_a/40}(x_a)\}$ are disjoint and
\begin{align*}
\cN_{K(n,B)}\cap B_1(p)\subseteq \bigcup_a  B_{s_{a}/4}(x_{a})\, ,
\end{align*}
where $s_a=s_{x_a}$. We claim that for the choosen covering,  every point $q$   is covered at most $c(n,\rv)$ times. Suppose $s:=d(q,\cC)$ and let $x$ be the center of those balls $B_{s_x/5}(x)$ which cover $q$. Then $s_x$ will satisfy $s_x/8<s<8s_x$. Therefore volume comparison will give the claimed bound $c(n,
\rv).$
Rearranging $x_a$ such that $d(x_{\alpha,i},\cC)=s_{\alpha,i}=d(x_{\alpha,i},z_{\alpha,i})$ with $s_{\alpha,i}\in (2^{-\alpha-1}, 2^{-\alpha}]$ and $z_{\alpha,i}\in\cC$, then we have
\begin{align*}
\cN_{K(n,B)}\cap B_1(p)\subseteq \bigcup_\alpha \bigcup_{i=1}^{N_\alpha} B_{s_{\alpha,i}/4}(x_{\alpha,i})\, ,	
\end{align*}
 and $\{B_{40^{-1}s_{\alpha,i}}(x_{\alpha,i})\}$ are disjoint.  Moreover, by Ahlfors assumption, for any fixed $\alpha$ we have that $\sharp_i\{B_{s_{\alpha,i}/2}(x_{\alpha,i})\}\le C(n,B,\rv)s_\alpha^{4-n}$ with $s_\alpha=2^{-\alpha}$.
Then we have
\begin{align*}
	& \int_{\cN_{K(n,B)}\cap B_1(p)} \brs{\Rm}^2(z)\,dz\\
	& \quad \le \sum_\alpha \sum_i \int_{B_{s_{\alpha,i}/4}(x_{\alpha,i})} \brs{\Rm}^2(z)\,dz \\
& \quad \leq C(n,\rv,B)\sum_\alpha \sum_i\left(\int_{B_{s_{\alpha,i}/2}(x_{\alpha,i})} (\brs{H}^4(z)+\delta^2)\,dz +\int_{B_{s_{\alpha,i}/20}(z_{\alpha,i})}
|{\cH}_{10s_{\alpha,i}^2}-{\cH}_{10^{-1}s_{\alpha,i}^2}|(x)\,d\mu(x)\right)\notag\\
&\quad \le\epsilon+ C(n,\rv,B)\int_{\cN_{K(n,B)/2}\cap B_1(p)} \brs{H}^4(z)\,dz+C(n,\rv,B)\sum_\alpha\int_{ B_{3/2}(p)}\big|{\cH}_{40s^2_\alpha}-{\cH}_{40^{-1}s^2_\alpha}\big|(x)\,d\mu(x),
\end{align*}
where we have used the fact that $\brs{\Ric^g}^2\le 2\brs{H}^2+ 2\delta^2$ and the fact that in the Vitali covering we choosen above, we can assume every point $p$ in the neck region is covered at most $c(n,\rv)$ times. By the monotonicity of $\cH$-volume, we have
\begin{align*}
& \int_{\cN_{K(n,B)}\cap B_1(p)} \brs{\Rm}^2(z)\,dz\\
&\quad \leq C(n,\rv,B) \int_{\cN_{K(n,B)}\cap B_1(p)} \brs{H}^4(z)\,dz+C(n,\rv,B)\int_{ B_{3/2}(p)}\big|\sum_\alpha({\cH}_{40s^2_\alpha}-{\cH}_{40^{-1}s^2_\alpha})\big|(x)\,d\mu(x)\\
&\quad \leq C(n,\rv,B)	\int_{\cN_{K(n,B)}\cap B_1(p)} \brs{H}^4(z)\,dz+C(n,\rv,B,\tau)\int_{ B_{3/2}(p)}\big|{\cH}_{40}-{\cH}_{40^{-5}r^2_x}\big|(x)\,d\mu(x)\, .
\end{align*}
On the other hand, noting that $B_2(p)$ is $\delta$-close to $\dR^{n-4}\times C(S^3/\Gamma)$, by the $\epsilon$-regularity Lemma \ref{l:neck:distance_harm_rad} on the neck region, for any $\delta'>0$ if $\delta\le \delta(n,\delta')$ we have
\begin{align*}
\int_{(\cN_{10^{-5}}\setminus \cN_{K(n,B)})\cap B_1(p)} \brs{\Rm}^2(z)\,dz\le&\ \sum_{x\in B_1(p)\cap \cC}\int_{B_{K(n,B)r_x}(x)\setminus \cup_{y\in \cC}B_{10^{-5}r_y}(y)}\brs{\Rm}^2(z)\,dz\\
\le&\ \sum_{x\in B_1(p)\cap \cC} C(n,B) \delta' r_x^{n-4}\le C(n,B)\delta'.
\end{align*}
Thus we arrive at 
\begin{gather} \label{f:final5.2}
\begin{split}
\int_{\cN_{10^{-5}}\cap B_1(p)} \brs{\Rm}^2(z)\,dz \le&\ C(n,B)\delta'+C(n,\rv,B)\int_{\cN_{K(n,B)}\cap B_1(p)} \brs{H}^4(z)\,dz\\
&\ \qquad +C(n,\rv,B)\int_{ B_{3/2}(p)}\big|{\cH}_{40}-{\cH}_{40^{-5}r^2_x}\big|(x)\,d\mu(x)\, .
\end{split}
\end{gather}

The final step is to deal with the $\cH$ term.  Since $\cN$ is a neck region we have that both $B_{\delta^{-1}}(x)$ and $B_{\delta^{-1}r_x}(x)$ are Gromov-Hausdorff close to $\dR^{n-4}\times C(S^3/\Gamma)$, and therefore by volume convergence we should expect that this term goes to zero as $\delta\to 0$.  We make this precise in the following claim.
 
\vskip 0.1in
 
\noindent \textbf{Claim:} For any $\epsilon>0$ and $t\le 10$, if $\delta\le \delta(n,\rv,\epsilon)$ and 
\begin{align*}
d_{GH}( B_{\delta^{-1}\sqrt{t}}(x), B_{\delta^{-1}\sqrt{t}}(0^{n-4},y_c))\le \sqrt{t} \delta, ~~\text{ where $(0^{n-4},y_c)$ is a cone vertex of $\dR^{n-4}\times C(S^3/\Gamma)$}
\end{align*}
then $|\cH_t(x)-\bar \cH_{\Gamma}|\le \epsilon$, where the constant $\bar \cH_{\Gamma}$ is the $\cH$-volume of the cone $\dR^{n-4}\times C(S^3/\Gamma)$, which only depends on $|\Gamma|$ and $n$.

\vskip 0.1in

\textbf{Proof of Claim:} By the scaling property of $\cH$-volume in Remark \ref{r:scaling_cHvolume}, it suffices to prove the case $t=1$. For simplicity, let us denote $\cH_1(p)=\int_M L(x,\delta)dx$ and $\bar \cH_{\Gamma}=\int_X L(x)dx$ with $X=\mathbb{R}^{n-4}\times C(S^3/\Gamma)$, which satisfies $L(x,\delta)\to L(x)$ uniformly on compact sets if $\delta\to 0$. Due to the exponential decay of $L(x,\delta)$, $L(x)$ and growth control of volume by the volume comparison, the integral of $L(x,\delta)$ over $M\setminus B_{R}(x)$ and the integral of $L(x)$ over $\dR^{n-4}\times C(S^3/\Gamma)\setminus B_{R}(0^{n-4},y_c))$ both must be smaller than $\epsilon/10$ if $R\geq R(n,\rv,\epsilon)$. For the integral over $B_{R}$, by volume convergence we have if $\delta\le \delta(n,\rv,\epsilon)$ that
\begin{align*}
\brs{\int_{B_R(p)}L(x,\delta)dx-\int_{B_R(0^{n-4},y_c))}L(x)dx} \le \epsilon/10.
\end{align*}
Combining the estimates on $B_R$ and outside $B_R$ we have proved the claim. $\square$

\vskip 0.1in

Choosing any $\delta\leq \delta'(n,\rv)$, we use the estimate on $\cH_t(x)$ from the claim in line (\ref{f:final5.2}) to finish the proof.
\end{proof}

We have the following immediate corollary of Theorem \ref{t:L2_neck}. 
For each $\alpha\in \mathbb N$, let us define a sequence of regions $\widetilde{\mathcal N_\alpha}\subset \mathcal N$ as follows:
$$\widetilde{\mathcal \cN_\alpha}=\mathcal \cN\cap \{d(\cdot,\mathcal C)>2^{-\alpha}\}$$
\begin{cor}
	Let $(M^n,g,H,p)$ satisfy $\Vol(B_1(p))>\rv>0$ and $\brs{\RC}\le \delta$. For any $\epsilon>0$, if  $\delta<\delta'(n,\epsilon,\rv)$ and  $\cN = B_2(p)\setminus \bigcup_{x\in \cC} \overline B_{r_x}(x)$ is a $\delta$-neck region, 
     then we have the curvature estimate 
     \begin{equation*}\label{B}\int_{\widetilde{\cN_{\alpha}}\cap B_1(p)} \brs{\Rm}^2(z)\,dz \leq \epsilon+\int_{\widetilde{\cN_{\alpha}}\cap B_1(p)} \brs{H}^4(z)\,dz.
     \end{equation*}\\
\end{cor}
\begin{proof}
This corollary is exactly the same as Theorem \ref{t:L2_neck} except that the domain is changed from $\cN_{10^{-5}}$ to $\widetilde{\cN_{\alpha}}$. Hence to prove this corollary, we only need to change the integral domain and the proof will be clear. We remark that to control the $L^2$ norm of curvature on a ball near the  boundary of  domain $\widetilde{\cN_{\alpha}}$, the  estimate from Proposition \ref{p:local_L2_neck} shows that it can be controlled by the $L^4$ norm of $H$ plus the $\cH$-volume on the doubled ball. Hence in the claimed estimate of Corollary \ref{B}, the right hand side should be integrated on a domain larger than  $\widetilde{\cN_{\alpha}}$.   But using Lemma \ref{t:neck:annulusRm} and Lemma  \ref{l:neck:scaleinvariantsmall}, we don't have to enlarge the domain.   We also emphasize that for given $\epsilon$, $\delta'$ is independent of $\alpha$.
\end{proof}
\section{Summable $L^4$ Estimate of $H$ on the neck region}
The main theorem of this section is the following summable integral estimate of $H$ on the whole neck region. 
\begin{thm}[$L^1$ Summable  Estimate of H]\label{t:splitting_neck_summable_H}
Let $(M^n,g,H,p)$ satisfy $\Vol(B_1(p))>\rv>0$ with $\brs{\RC}\le \delta$.  Then for every $\epsilon>0,$  if $\delta<\delta'(n,\epsilon,\rv)$ and $\cN = B_2(p)\setminus \overline B_{r_x}(\cC)$ is a $\delta$-neck region, then we have the estimate
\begin{align}\label{aaa}
\int_{b\le R} b^{-3} \brs{H}\phi_{\mathcal N}(z)\,dz < \epsilon+C\left(\int_{b\le R}|b|^{-2}\brs{H}^2\phi_\cN(z)\,dz\right)^{1/2}\left(\int_{b\le R}\brs{\Rm}^2\phi_\cN+\epsilon\,dz\right)^{1/2}\,,
\end{align}
\end{thm}
\begin{rmk}\label{H}  Before  proving the above theorem, we discuss why the above estimate is helpful. On a  $\delta$-neck region, we will show that $b^{-2}|H|^2<\tilde\epsilon^2(\delta) b^{-3}|H|$, where $\tilde\epsilon(\delta)\to 0$ as $\delta\to 0$ and  $\brs{\brs{\Rm}}^2_{L^2}\approx \|H\|^4_{L^4}$. Hence, in \eqref{aaa}, if $\|H\|_{L^4}\leq C(n,\rv)$ then the right hand side of the integral  is  indeed a $\tilde\epsilon(\delta)<<1$ mutiple of the square root of the left hand side integral. A little bit more precisely, $$\int_{b\le R} b^{-3} \brs{H}\phi_{\mathcal N}(z)\,dz < \epsilon+\tilde\epsilon(\delta)C(n,\rv)\left(\int_{b\le R}|b|^{-3}\brs{H}\phi_\cN(z)\,dz\right)^{1/2}.\,$$
Now the point is if we choose $\tilde\epsilon(\delta)<<\epsilon$, we will have an estimate $\int_{b\le R} b^{-3} \brs{H}\phi_{\mathcal N}(z)\,dz <2\epsilon$, which in turn will imply that $\|H\|^4_{L^4}\leq \epsilon$.  

This gives the opportunity to close the loop in the final arguments.  In particular, for the proof of Theorem \ref{t:neck_region}, the estimate of $\|H\|^4_{L^4}\leq \epsilon$ will be done inductively.  We will approximate our neck region by cutting it off below some scale.  Assuming we have proved the estimate on this cutoff neck region, we will want to prove it on a neck region which is a little larger.  Adding this region a priori adds some small extra amount to the $|H|^4$ integral, but by the self improvement estimate above on $\|H\|_{L^4}$, we will eventually get back to $\|H\|_{L^4}\leq\epsilon$.
\end{rmk}
To prove the above we will first discuss a  superconvexity estimate in Section \ref{sss:superconvexity} using similar ideas from  \cite{NV},\cite{JiNa}.  Analysis of the derived ODE together with the ${\bar \mu}$-Green's estimates of Section \ref{ss:splitting:Greens_standard}  will then be applied in order to conclude Theorem \ref{t:splitting_neck_summable_H} itself.  We first show here that the tensor $H$ is small in a scale-invariant sense on the $\delta$-neck region, which relies on the strong rigidity of generalized Ricci-flat cones.
\begin{lemma}\label{l:neck:Hsmall}
	Let $(M^n,g,H,p)$ satisfy $\Vol(B_1(p))>\rv>0$ and $\brs{\RC}\le \delta$.  Then for any  $\epsilon$ small, there exists $\delta'(\rv,n,\epsilon)$ such that if $\mathcal N$ is a $\delta$-neck region with $\gd < \gd'(n,\rv)$, then for each $y\in \cN_{10^{-6}}$ it holds that
	\begin{align*}
	\brs{H}\leq \epsilon {d(y,\cC)}^{-1}\, .
	\end{align*}
\end{lemma}
\begin{proof}
We argue by contradiction.  If no such $\delta'$ exists, there is a sequence of $\delta_i$-neck regions $\mathcal N_i$ and a sequence of points $y_i\in\mathcal N_{i,10^{-6}}$ with $\brs{H}(y_i)d(y_i,\mathcal C_i)>\epsilon_0>0$ and $\delta_i\to 0$.
Then by a standard blow up argument, we may assume $\brs{H}(y_i)>\epsilon_0$ and $\mathcal N_i\to \mathbb R^{n-4}\times C(Y)$ in Gromov-Hausdorff topology. By Lemma \ref{H-control} and the cone rigidity of Proposition \ref{cone}, we have $\brs{H}(y_i)\to 0$ which is a contradiction.
\end{proof}

We also need the following differential inequality for tensor $H$.
\begin{lemma} \label{l:Hellipticineq} Let $(M^n, g)$ be a Riemannian manifold and fix $H \in \Lambda^*$.  Then
\begin{align*}
    \gD \brs{H} \geq - c(n) \brs{\Rm} \brs{H}-\brs{dd_g^*H}-\brs{d_g^*dH}.
\end{align*}
\begin{proof} Note  $\gD_d H =dd^*_gH+d^*_gdH$, thus by the Bochner formula we have
\begin{align} \label{f:Hbochner}
    \gD \brs{H}^2 = \brs{\N H}^2 + \Rm \star H^{*2}+\langle H,dd^*_gH\rangle+\langle H,d^*_gdH\rangle.
\end{align}
Now  Kato's inequality  implies $\sqrt{|H|^2+\epsilon}\Delta\sqrt{|H|^2+\epsilon}\leq |H||\nabla H|$ for any $\epsilon>0$. Then let $\epsilon\to 0$, we get $\brs{\N \brs{H}} \leq \brs{\N H}$ holds in the distribution sense.  Combining this inequality with (\ref{f:Hbochner}) yields the result.
\end{proof}
\end{lemma}
By the  quantitative stratification estimate Theorem \ref{t:quant_strat}, when $p<4$, we can easily show that the $L^p$ norm of $|H|$ is uniformly bounded.
\begin{prop}\label{p:prelim:harmonic_estimates}
	Let $(M^n,g,H,p)$ satisfy $\brs{\RC_{M^n}}\leq n-1$ with $\Vol(B_2(p))>\rv>0$.  Then the following hold:
\begin{enumerate}
\item For each $0<q<4$ we have that $\int_{B_1(p)}\brs{H}^q(x)\,dx\leq C(n,\rv,q)$.
\item For each $0<q<2$ we have that $\int_{B_1(p)}|\nabla H|^q(x)\,dx\leq C(n,\rv,q)$.
\end{enumerate}
\end{prop}

\begin{proof} 
Using $\brs{\RC_{M^n}}\leq n-1$, we have $\brs{H}^2\leq \brs{\Ric^g}+n-1$, thus item $(1)$ follows from Theorem  \ref{t:main_estimate}. Now we prove item (2).  We claim that for $B_{2r}(x)\subseteq B_2(p)$ with $r_{h}(x)\geq 2r$ and $q\leq 2$, then one has $r^{2q-n}\int_{B_r(x)}|\nabla H|^q<C(n,q)$.  

It is enough to prove the claim for $q=2$, the general case then follows from H\"older's inequality. By Lemma \ref{l:neck:cutoff}, choose $\phi:B_{2r}(x)\to \dR$ such that $\phi\equiv 1$ on $B_r(x)$, $\phi\equiv 0$ outside of $B_{3r/2}(x)$ with $r|\nabla\phi|, r^2|\Delta\phi|\leq C(n)$.  Multiplying both sides of (\ref{f:Hbochner}) by $\phi$ and noting the pointwise upper bound $\brs{H}\le C(n,\rv)r^{-1}, \brs{\Rm}\le c(n,\rv) r^{-2}$ , we arrive at
\begin{align*}
2r^{4-n}\int_{B_{2r}(x)} \phi(z) |\nabla H|^2(z)\,dz \leq&\ Cr^{4-n}\int_{B_{2r}(x)} r^{-4} +r^{4-n}\int_{B_{2r}(x)}\phi\langle H,dd^*_gH\rangle\\
&\ +r^{4-n}\int_{B_{2r}(x)}\phi\langle H,d^*_gdH\rangle.
\end{align*}
Furthermore we can estimate
\begin{align*}
\int_{B_{2r}(x)}\phi\langle H,dd^*_gH\rangle =&\ \int_{B_{2r}(x)}\langle d_g^*(\phi H),d^*_gH\rangle\le 2\int_{B_{2r}(x)}\brs{d^*_gH}^2+ \int_{B_{2r}(x)}r^{-2}\brs{H}^2,\\
\int_{B_{2r}(x)}\phi\langle H,d^*_gdH\rangle =&\ \int_{B_{2r}(x)}\langle d (\phi H),d H\rangle\le 2\int_{B_{2r}(x)}\brs{dH}^2+ \int_{B_{2r}(x)}r^{-2}\brs{H}^2.
\end{align*}
Using $\brs{d^*_gH}\le n-1+\brs{H}^2,\brs{dH}\le n-1+\brs{H}^2$
will  finish the proof of the claim, and item (2) follows by the same argument in Theorem \ref{t:main_estimate}.  
\end{proof}

\subsection{The Superconvexity Equation}\label{sss:superconvexity}
In this subsection, we prove the key superconvexity estimate for the tensor $H$. Similar estimates are used in \cite{NV,JiNa} in different context  or for different purposes. 
\begin{prop}\label{p:superconvexity}
Let $(M^n,g,H,p)$ satisfy $\Vol(B_1(p))>\rv>0$ and $\brs{\RC}\le \delta$.  Then for every $\epsilon>0$  if $\delta\leq \delta'(n,\epsilon,\rv)>0$ and $\cN = B_2(p)\setminus \overline B_{r_x}(\cC)$ is a $\delta$-neck region, fix a cutoff function $\phi_{\cN}$ via Lemma \ref{l:neck:cutoff}, and set 
\begin{align*}
F(r)=&\ r^{-3}\int_{b=r}\sqrt{\brs{H}^2+{\epsilon}^2}\phi_\cN|\nabla b|(z)\,dz, \qquad Q(r)=\ r F(r).
\end{align*} 
Then we have
\begin{align*}
r^2Q''(r)&+rQ'(r)-Q(r)\\
\geq&\  e(r)=\int_{b=r}\left(-\sqrt{\brs{H}^2+{\epsilon}^2}~|\Delta \phi_\cN|-c(n)|Rm|\cdot\brs{H}\phi_\cN+U_1+U_2\right)|\nabla b|^{-1}(z)\,dz,
\end{align*}
where $U_1=2\langle\nabla\phi_\cN,\nabla \sqrt{\brs{H}^2+{\epsilon}^2}\rangle, U_2= -\phi_\cN\frac{\langle H,dd^*_gH\rangle+\langle H,d^*_gdH\rangle}{\sqrt{\brs{H}^2+{\epsilon}^2}}$.
\end{prop}

\begin{proof}
As in Lemma \ref{l:Hellipticineq} we have
\begin{align*}
    \gD \brs{H}^2 = 2\brs{\N H}^2 + \Rm \star H^{*2}+2\langle dd_g^*H,H\rangle.
\end{align*}
  Let $f=\phi_\cN \sqrt{\brs{H}^2+{\epsilon}^2}$. Then
\begin{align*}
\Delta f&=\Delta \phi_\cN ~\sqrt{\brs{H}^2+{\epsilon}^2}+\Delta \sqrt{\brs{H}^2+{\epsilon}^2}\, \phi_\cN+2\langle\nabla\phi_\cN,\nabla \sqrt{\brs{H}^2+{\epsilon}^2}\rangle\\ \nonumber
&\geq -\sqrt{\brs{H}^2+{\epsilon}^2}~|\Delta \phi_\cN|-c(n)|Rm|\cdot\brs{H}\phi_\cN+U_1+U_2,
\end{align*}
where $U_1,U_2$ are  defined as in the lemma.
Note that $Q(r)=rF(r)=r^{-2}\int_{b=r}f|\nabla b|(z)\,dz$ with $b^{-2}=G_{\bar \mu}$, so by the Green formula and Lemma \ref{l:green_formula}, we have
\begin{align*}
Q''(r)+\frac{1}{r}Q'-\frac{1}{r^2}Q&=r^{-2}\int_{b=r}\Delta f~|\nabla b|^{-1}(z)\,dz\\ \nonumber
&\geq r^{-2}\int_{b=r}\left(-\sqrt{\brs{H}^2+{\epsilon}^2}~|\Delta \phi_\cN|-c(n)\brs{\Rm}\cdot\brs{H}\phi_\cN+U_1+U_2\right)|\nabla b|^{-1}(z)\,dz.
\end{align*}
Thus we finish the proof.
\end{proof}

Of course, in order to exploit the above superconvexity we will want to apply a maximum principle in order to obtain bounds for $H$.  To accomplish this we first need to find special solutions of the above ODE to compare against.

\begin{lemma} \label{l:special_solution3}
Consider on the interval $(0,R)$ the differential equation
$$h''(r)+\frac{1}{r}h'(r)-\frac{1}{r^2}h(r)={g(r)},$$
for some bounded $g(r)$ which may change signs. Then we have a solution $h(r)$ such that
\begin{align*}
\int_0^R\frac{h(r)}{r}dr=-\int_0^Rsg(s)ds+\frac{1}{2R}\int_0^Rs^2g(s)ds,
\end{align*}
with $h(R)=-\frac{1}{2R}\int_0^Rs^2g(s)ds$ and $h(0)=0$.
\end{lemma}
\begin{proof}Let us define $h(r)=\frac{1}{2}\left(-r\int_{r}^R{g}ds-r^{-1}\int_{0}^rs^2{g}ds\right)$, which one can easily check is a solution satisfying the claimed boundary conditions.  We furthermore observe that
\begin{align*}
-2\int_0^R r^{-1}h(r)dr&=\int_0^R\int_{r}^R{g(s)}dsdr+\int_0^Rr^{-2}\int_{0}^rs^2{g}dsdr\\ \nonumber
&=\int_0^Rg(s)\int_0^sdr ds+\int_0^Rs^2g(s)\int_s^Rr^{-2}drds\\ \nonumber
&=2\int_0^Rsg(s)ds-\frac{1}{R}\int_0^Rs^2g(s)ds,
\end{align*}
as claimed.
\end{proof}

We will now use the above comparison solution to provide estimates on the $L^1$ norm of $H$ on a neck region:

\begin{thm}\label{t:superconvexity_estimate}
Let $(M^n,g,H,p)$ satisfy $\Vol(B_1(p))>\rv>0$ and $\brs{\RC}\le \delta$.  Then for every $\epsilon>0$  if $\delta\leq \delta'(n,\epsilon,\rv)>0$ and $\cN = B_2(p)\setminus \overline B_{r_x}(\cC)$ is a $\delta$-neck region,
fix a cutoff function $\phi_{\cN}$ via Lemma \ref{l:neck:cutoff}, and set 
\begin{align*}
F(r)=&\ r^{-3}\int_{b=r}\sqrt{\brs{H}^2+{\epsilon}^2}\phi_\cN|\nabla b|(z)\,dz, \qquad Q(r)= r F(r).
\end{align*} 
Then we have the Dini estimate:
\begin{align*}
\int_0^\infty \frac{1}{r}Q(r)dr\leq C(n,\rv,B)\epsilon+C\left(\int_{b\le R}|b|^{-2}\brs{H}^2\phi_\cN(z)\,dz\right)^{1/2}\left(\int_{b\le R}\brs{\Rm}^2\phi_\cN(z)+\epsilon\,dz\right)^{1/2}.
\end{align*}
\end{thm}
\begin{proof}
By the definition of $\phi_\cN$ and Proposition \ref{p:green_function_distant_function}, there exists $R=R(n,\rv,B)$ such that $\phi_\cN\equiv 0$ on $\{b\geq R\}$. Therefore we have
$$\int_0^\infty \frac{1}{r}Q(r)\,dr=\int_0^R\frac{1}{r}Q(r)\,dr\, .$$
We first use Lemma \ref{l:special_solution3} to choose $h_1(r)$ such that
\begin{align*}
h_1''+\frac{1}{r}h_1'-\frac{1}{r^2}h_1= r^{-2}\int_{b=r}U_1 |\nabla b|^{-1}(z)\,dz=r^{-2}\int_{b=r}2\langle\nabla\phi_\cN,\nabla \sqrt{\brs{H}^2+{\epsilon}^2}\rangle|\nabla b|^{-1}(z)\,dz,
\end{align*}
 with $h_1(0) = 0, h_1(R)=-\frac{1}{R} \int_{b\le R}\langle\nabla\phi_\cN,\nabla \sqrt{\brs{H}^2+{\epsilon}^2}\rangle(z)\,dz $, and
\begin{align}\label{e:H_1_integral}
\int_0^R\frac{h_1(r)}{r}dr=-2\int_{b\le R}b^{-1}\langle\nabla\phi_\cN,\nabla \sqrt{\brs{H}^2+{\epsilon}^2}\rangle(z)\,dz+\frac{1}{R}\int_{b\le R}\langle\nabla\phi_\cN,\nabla \sqrt{\brs{H}^2+{\epsilon}^2}\rangle(z)\,dz.
\end{align}

Again apply Lemma \ref{l:special_solution3} to produce $h_2(r)$ satisfying
\begin{align}\label{e:H_2_integral}
h_2''+\frac{1}{r}h_2'-\frac{1}{r^2}h_2= r^{-2}\int_{b=r}\left(-\sqrt{\brs{H}^2+{\epsilon}^2}~|\Delta \phi_\cN|-c(n)|Rm|\cdot\brs{H}\phi_\cN\right) |\nabla b|^{-1}(z)\,dz,
\end{align}
with $h_2(0) = 0, h_2(R)=-\frac{1}{2R} \int_{b\le R}\left(-\sqrt{\brs{H}^2+{\epsilon}^2}~|\Delta \phi_\cN|-c(n)|Rm|\cdot \brs{H}\phi_\cN\right)(z)\,dz$ and
\begin{align}\label{e:H_3_integral}
\int_0^R\frac{h_2(r)}{r}dr=\int_{b\le R}\left(-b^{-1}+\frac{1}{2R}\right)\left(-\sqrt{\brs{H}^2+{\epsilon}^2}~|\Delta \phi_\cN|-c(n)|Rm|\brs{H}\phi_\cN\right)(z)\,dz\, .
\end{align}

Again apply Lemma \ref{l:special_solution3} to produce $h_3(r)$ satisfying (We skip the term $\langle H,d_g^*dH\rangle$, which will can be dealt the same as $\langle H,dd_g^*H\rangle$)
\begin{align}\label{e:H_2_integral1}
h_3''+\frac{1}{r}h_3'-\frac{1}{r^2}h_3= r^{-2}\int_{b=r}\phi_\cN\frac{\langle H,dd^*_gH\rangle}{\sqrt{\brs{H}^2+{\epsilon}^2}}|\nabla b|^{-1}(z)\,dz,
\end{align}
with $h_3(0) = 0, h_3(R)=-\frac{1}{2R} \int_{b\le R} \phi_\cN\frac{\langle H,dd^*_gH\rangle}{\sqrt{\brs{H}^2+{\epsilon}^2}}\,dz$ and
\begin{align}\label{e:H_3_integral1}
\int_0^R\frac{h_3(r)}{r}dr=\int_{b\le R}\left(-b^{-1}+\frac{1}{2R}\right)\phi_\cN\frac{\langle H,dd^*_gH\rangle}{\sqrt{\brs{H}^2+{\epsilon}^2}}\,dz\, .
\end{align}

Finally let us choose $h_4(r)$ such that
\begin{align}
h_4''+\frac{1}{r}h_4'-\frac{1}{r^2}h_4=0
\end{align}
with $h_4(0)=0$ and $h_4(R)=|h_1(R)|+|h_2(R)|+|h_3(R)|+|Q(R)|=|h_1(R)|+|h_2(R)|+|h_3(R)|\equiv A$.  Note that we can explicitly solve $h_4(r)=Ar/R$. On the other hand, we now have that
\begin{align}\label{e:H1H2H3H4}
\left(Q-\sum_{i=1}^4 h_i\right)''(r)+\frac{1}{r}\left(Q-\sum_{i=1}^4 h_i\right)'(r)-\frac{1}{r^2}\left(Q-\sum_{i=1}^4 h_i\right)(r)\geq 0,
\end{align}
with $\left(Q-\sum_{i=1}^4 h_i\right)(0)=0$ and $\left(Q-\sum_{i=1}^4 h_i\right)(R)\le 0$.  Therefore, 
we have $Q- \sum_{i=1}^4 h_i\le 0$. Actually, assume $(Q- \sum_{i=1}^4 h_i)(r_0)=\max_{0\le r\le R} (Q- \sum_{i=1}^4 h_i)(r)>0$, we can compute
\begin{align*}
\left(Q-\sum_{i=1}^4 h_i\right)''(r_0)+\frac{1}{r_0}\left(Q-\sum_{i=1}^4 h_i\right)'(r_0)-\frac{1}{r^2_0}\left(Q-\sum_{i=1}^4 h_i\right)(r_0)\le -\frac{1}{r^2_0}\left(Q-\sum_{i=1}^4 h_i\right)(r_0)<0,
\end{align*}
which contradicts to \eqref{e:H1H2H3H4}. Thus 
we have $Q- \sum_{i=1}^4 h_i\le 0$.  Observe also that $Q\geq 0$ while $h_i$ may change signs.  We may now estimate
\begin{align}\label{e:H_bound}
\int_0^R \frac{Q(r)}{r}dr\le \sum_{i=1}^4\int_0^R \frac{h_i(r)}{r}dr.
\end{align}
Therefore, to estimate the Dini integral for $Q(r)$, we only need to control the Dini integral of each $h_i$.  Beginning with $h_1$, we have by (\ref{e:H_1_integral}) and integrating by parts that
\begin{align*}
\int_0^R\frac{h_1(r)}{r}dr&=-2\int_{b\le R}b^{-1}\langle\nabla\phi_\cN,\nabla \sqrt{\brs{H}^2+{\epsilon}^2}\rangle(z)\,dz +\frac{1}{R}\int_{b\le R}\langle\nabla\phi_\cN,\nabla \sqrt{\brs{H}^2+{\epsilon}^2}\rangle(z)\,dz\\ \nonumber
&\le C(n,\rv,B)\left(\int_{b\le R}\Big(b^{-2}|\nabla \phi_\cN|+b^{-1}|\Delta\phi_\cN|\Big) \sqrt{\brs{H}^2+{\epsilon}^2}\right)(z)\,dz\, ,
\end{align*}
where we have used that $\phi_\cN\equiv 0$ on $\{b=R\}$ and $|\nabla b|\le C(n,\rv,B)$ on $\text{ supp}\, \phi_\cN$.  Now using Lemma \ref{l:neck:Hsmall} and Proposition \ref{p:green_function_distant_function} with $\delta,\delta'\leq \delta(\epsilon,n,\rv)\le \epsilon^5$ we can obtain the a priori estimate $\brs{H}\le \epsilon b^{-1}$ on $\text{supp }\phi_\cN$.  Plugging this in gives us
\begin{align*}
\int_0^R\frac{h_1(r)}{r}dr
&\le C(n)\epsilon\left(\int_{b\le R}\big(b^{-3}|\nabla \phi_\cN|+b^{-2}|\Delta\phi_\cN|\big)\right)(z)\,dz\, .
\end{align*}

Finally for $h_2$, by (\ref{e:H_3_integral}) and the estimate of $H$ on the neck region by Lemma \ref{l:neck:Hsmall}, we have
\begin{align*}
\int_0^R\frac{h_2(r)}{r}dr&=\int_{b\le R}\left(-b^{-1}+\frac{1}{2R}\right)\left(-\sqrt{\brs{H}^2+{\epsilon}^2}~|\Delta \phi_\cN|-c(n)\brs{\Rm}\cdot \brs{H}\phi_\cN\right)(z)\,dz\\ \nonumber
&\le C(n)\epsilon\int_{b\le R} b^{-2}|\Delta\phi_\cN|(z)\,dz+C\left(\int_{b\le R}|b|^{-2}\brs{H}^2\phi_\cN(z)\,dz\right)^{1/2}\left(\int_{b\le R}\brs{\Rm}^2\phi_\cN(z)\,dz\right)^{1/2}.
\end{align*}
Similarly,
\begin{equation}
\begin{aligned}\label{e:H_3_integral2}
\int_0^R\frac{h_3(r)}{r}dr&\le 2\int_{b\le R}b^{-1}\phi_\cN\frac{\langle d^*_gH,d^*_gH\rangle}{\sqrt{\brs{H}^2+{\epsilon}^2}}\,dz\,+ 2\int_{b\le R}(b^{-2}|\nabla b|\phi_\cN+b^{-1}|\nabla\phi_\cN|)\frac{\brs{H}\brs{d_g^*H}}{\sqrt{\brs{H}^2+{\epsilon}^2}}\,dz\,\\&+2\int_{b\le R}b^{-1}\phi_\cN\frac{\brs{H}^2|\nabla\brs{H}|\brs{d_g^*H}}{(\sqrt{\brs{H}^2+{\epsilon}^2})^3}\,dz\,
\end{aligned}\end{equation}
By choosing $\delta<\epsilon^2$ such that $|d_g^*H|\le\epsilon^2+\brs{H}^2,|dH|\le\epsilon^2+\brs{H}^2$ and using Lemma \ref{l:neck:Hsmall}, we have $ \brs{H}<\epsilon b^{-1}$, and so the right hand side of
the above inequality is smaller than $$2\epsilon\int_{b\le R}b^{-3}\phi_\cN\brs{H}\,dz\,+2\epsilon\int_{b\le R}\Big(b^{-2}\phi_\cN+b^{-3}|\nabla\phi_\cN|\Big)\,dz\,+2\int_{b\le R}b^{-1}\brs{H}\phi_\cN|\nabla\brs{H}|\,dz\,.$$
For the last term of above inequality, by H\"older inequality
\begin{align}
 \int_{b\le R}b^{-1}\brs{H}\phi_\cN|\nabla\brs{H}|\,dz\,<\Big(\int_{b\le R}b^{-2}\brs{H}^2\phi_\cN dz\,\Big)^{1/2}\Big(\int_{b\le R}|\nabla\brs{H}|^2\phi_\cN dz\,\Big)^{1/2}
 \end{align}
 Then we integrate by parts again to obtain
 \begin{align*}
 \int_{b\le R}|\nabla\brs{H}|^2\phi_\cN dz\,\le&\ \int_{b\le R}\Delta\phi_{\mathcal N}\brs{H}^2+\int_{b\le R}\phi_{\mathcal N}(\brs{\Rm}\brs{H}^2+|dH|^2+|d^*_gH|^2)\\
 &\ \qquad +\int_{b\le R}|\nabla\phi_{\mathcal N}|\brs{H}(|dH|+|d^*_gH|)\\
 \le&\ 2\epsilon\int_{b\le R}\Big(b^{-2}|\Delta\phi_\cN|+b^{-3}|\nabla\phi_\cN|\,dz\,\Big)+4\int_{b\le R}\Big(\brs{\Rm}^2\phi_\cN+\brs{H}^4\phi_\cN+\epsilon\Big)\,dz\,.
 \end{align*}
We remark that the term $\int_{b\le R}b^{-2}|\Delta\phi_\cN|\,dz\,+b^{-3}|\nabla\phi_\cN|\,dz\,$will be shown to bounded by $C(n,\rv)$.
 
 To estimate $h_4$, it suffices to estimate $A=|h_1(R)|+|h_2(R)|+|h_3(R)|$. Actually, from the formulas of $h_1(R), h_2(R),h_3(R)$ and the formulas \eqref{e:H_1_integral}, \eqref{e:H_2_integral}, \eqref{e:H_3_integral},
\eqref{e:H_2_integral1} and the above estimates, we have already proved
\begin{align} 
|A|\le  C(n,B,\rv)\epsilon\left(\int_{b\le R}b^{-2}|\Delta \phi_\cN|(z)+b^{-3}|\nabla\phi_\cN|(z)+b^{-1}|\nabla H|(z)\phi_\cN\,dz\right).
\end{align}

Combining all these with \eqref{e:H_bound} we get
\begin{align}\label{e:L1_summable}
\int_0^R \frac{Q(r)}{r}dr &\le C(n,B,\rv)\epsilon\left(\int_{b\le R}b^{-2}|\Delta \phi_\cN|(z)+b^{-3}|\nabla\phi_\cN|(z)+b^{-1}|\nabla H|(z)\phi_\cN\,dz\right)\\
 \nonumber&+C(n,\rv)\left(\int_{b\le R}|b|^{-2}\brs{H}^2\phi_\cN(z)\,dz\right)^{1/2}\left(\int_{b\le R}\brs{\Rm}^2\phi_\cN(z)\,dz+\epsilon\right)^{1/2}.
\end{align}

Hence, to get the Dini estimate it suffices to estimate each term of (\ref{e:L1_summable}).  In fact, from the definition of $\phi_\cN$ in Lemma \ref{l:neck:cutoff}, the definition of neck region, and by the comparison of $b$ and $d_\cC$ in Proposition \ref{p:green_function_distant_function}, one can easily get 
\begin{gather}\label{neckestimate}
\begin{split}
\int_{b\le R} & b^{-2}|\Delta \phi_\cN|(z)\,dz+\int_{b\le R}b^{-3}|\nabla\phi_\cN|(z)\,dz\\
&\le C(n,B,\rv)\left(\int_{d_\cC\le 3}d_\cC^{-2}|\Delta \phi_\cN|(z)\,dz+\int_{d_\cC\le 3}d_\cC^{-3}|\nabla\phi_\cN|(z)\,dz\right)
\\ 
&\le C(n,B,\rv)\left(\sum_{x\in\cC\cap B_{19/10}(p)}\int_{B_{r_x}(x)}\left(d_\cC^{-2}|\Delta \phi_\cN|+d_\cC^{-3}|\nabla \phi_\cN|\right)(z)\,dz+\int_{\text{ supp }\phi_\cN}d_\cC^{-3}(z)\,dz\right)\\ 
&\le  C(n,B,\rv)\left(\sum_{x\in \cap B_{19/10}(p)}r_x^{n-4}+\int_{\{d_\cC\le 3\}\cap \text{ supp }\phi_\cN}d_\cC^{-3}(z)\,dz\right)\\ 
&\le C(n,B,\rv)\mu_\cN(B_{19/10}(p))+C(n,B,\rv)\int_{\{d_\cC\le 3\}\cap \text{ supp }\phi_\cN}d_\cC^{-3}(z)\,dz\\
&\le C(n,B,\rv)+C(n,B,\rv)\int_{\{d_\cC\le 3\}\cap \text{ supp }\phi_\cN}d_\cC^{-3}(z)\,dz\, ,
\end{split}
\end{gather}
where the term $\int_{\text{ supp }\phi_\cN}d_\cC^{-3}(z)\,dz$ comes from the derivative of $\phi_1$ in Lemma \ref{l:neck:cutoff} ($\phi_2=1$ on supp $\phi_\cN$).  We can compute that 
\begin{gather}\label{e:dcC-31}
\begin{split}
\int_{\text{ supp }\phi_\cN}d_\cC^{-3}(z)\,dz&\le \sum_{i=0}^\infty\int_{\text{ supp }\phi_\cN \cap \{2^{-i}\le d(x,\cC)\le 2^{-i+1}\}}d_\cC^{-3}(z)\,dz\\
&\le \sum_{i=0}^\infty 2^{3i}\Vol(\text{supp }\phi_\cN \cap \{2^{-i}\le d(x,\cC)\le 2^{-i+1}\}).
\end{split}
\end{gather}
In Lemma \ref{l:neck:scaleinvariantsmall}, we have proved that,
\begin{gather}\label{e:dcC-3211}
\begin{split}
\Vol(\text{supp }\phi_\cN \cap \{2^{-i}\le d(x,\cC)\le 2^{-i+1}\})\le&\ \sum_{\alpha=1}^N\Vol(B_{2^{-i+3}}(x_\alpha))\\
\le&\ C(n,B)2^{(n-4)i} C(n)2^{n(-i+3)}\le C(n,B)2^{-4i}.
\end{split}
\end{gather}
Combining \eqref{e:dcC-31} and \eqref{e:dcC-3211} we get 
\begin{align}
\int_{\text{supp }\phi_\cN}d_\cC^{-3}(z)\,dz\le C(n,B).
\end{align}
Using this in (\ref{neckestimate}) we arrive at 
\begin{align}
\int_{b\le R}b^{-2}|\Delta \phi_\cN|(z)\,dz+\int_{b\le R}b^{-3}|\nabla\phi_\cN|(z)\,dz\le C(n,\rv,B).
\end{align}

On the other hand, we know from Theorem \ref{p:prelim:harmonic_estimates} that $\int_{B_{1}(p)}|\nabla  H|^{3/2}\le C(n,\rv)$. Therefore, by the H\"older inequality we have
\begin{align*}
\int_{b\le R}b^{-1}|\nabla H|\phi_\cN(z)\,dz\le \left(\int_{b\le R}b^{-3}\phi_\cN(z)\,dz\right)^{1/3}\left(\int_{b\le R}|\nabla H|^{3/2}\phi_\cN(z)\,dz\right)^{2/3}\le C(n,B,\rv)\, .
\end{align*}
Combining the above estimates proves the theorem.
\end{proof}

\begin{proof}[Proof of Theorem \ref{t:splitting_neck_summable_H}]
In order to finish the proof let us apply Theorem \ref{t:superconvexity_estimate} with the given $\epsilon$ in order to conclude
\begin{align*}
	\int_0^\infty r^{-3}\int_{b=r} \brs{H}\phi_\cN(z)\,dz dr  \le&\
\int_0^\infty \frac{1}{r}Q(r)dr\\
\leq&\ C(n,\rv,B)\epsilon+C\left(\int_{b\le R}|b|^{-2}\brs{H}^2\phi_\cN(z)\,dz\right)^{1/2}\left(\int_{b\le R}\brs{\Rm}^2\phi_\cN(z)+\epsilon\,dz\right)^{1/2}.
\end{align*}
By using the coarea formula, it holds that
\begin{align*}
\int_{M} b^{-3}\,|\nabla b|\, \brs{H} \phi_\cN(z)\,dz=\int_0^\infty r^{-3}\int_{b=r} \brs{H}\phi_\cN(z)\,dz dr.	
\end{align*}
Finally, using the Green's function estimates $|\nabla b|> C^{-1}(n,B)$
from Proposition \ref{p:green_function_distant_function} will finish the proof.
\end{proof}
We have the immediate corollary of Theorem \ref{t:splitting_neck_summable_H}.
For each $\alpha\in \mathbb N$, let us define a sequence of regions $\mathcal N_\alpha\subset \mathcal N$ as follows:
$$\mathcal \cN_\alpha=\mathcal \cN\cap \{b>2^{-\alpha}\}$$
We also define a sequence of cutoff functions $$\phi_\alpha=\psi(b/2^{-\alpha})\phi_{\mathcal N}$$ where $\psi$ is a fixed smooth one-variable function with $\psi(x)=0$ when $x\leq\frac{1}{2}$ and $\psi(x)=1$ when $x\geq1$. 

\begin{cor}\label{C}

Let $(M^n,g,H,p)$ satisfy $\Vol(B_1(p))>\rv>0$ and $\brs{\RC}\le \delta$.  Then for every $\epsilon>0$, if $\delta\leq \delta'(n,\epsilon,\rv)>0$ and $\cN = B_2(p)\setminus \overline B_{r_x}(\cC)$ is a $\delta$-neck region, then we have the estimate
\begin{align*}
\int_{b\leq R} b^{-3} \brs{H}\phi_{\alpha}(z)\,dz < \epsilon+C(n,\rv)\left(\int_{b\le R}|b|^{-2}\brs{H}^2\phi_\alpha(z)\,dz\right)^{1/2}\left(\int_{b\le R}\brs{\Rm}^2\phi_\alpha(z)+\epsilon\,dz\right)^{1/2}\,.
\end{align*}\\
\end{cor}
\begin{proof}
The corollary is exactly the same as Theorem \ref{t:splitting_neck_summable_H} except that here we use a different cutoff function. In the proof of Theorem  \ref{t:splitting_neck_summable_H}, we have proved the key estimate  for the cutoff function $\phi_{\cN}$
\begin{align}\label{abcd}
\int_{b\le R}b^{-2}|\Delta \phi_\cN|(z)\,dz+\int_{b\le R}b^{-3}|\nabla\phi_\cN|(z)\,dz\le C(n,\rv,B).
\end{align}
Hence to prove Corollary \ref{C} holds, we only need to show  that similar inequality as above holds when the cutoff function $\phi_{\mathcal N}$ is replaced by $\phi_\alpha$. Then the rest of the proof will be verbatim as the proof of Theorem \ref{t:splitting_neck_summable_H}.

\vskip 0.1in
\textbf{Claim:} The cutoff function $\phi_\alpha$ satisfies the following estimate
\begin{align}
\int_{b\le R}b^{-2}|\Delta \phi_\alpha|(z)\,dz+\int_{b\le R}b^{-3}|\nabla\phi_\alpha|(z)\,dz\le C(n,\rv,B).
\end{align}

\vskip 0.1in

\textbf{Proof of Claim:} First we have the simple fact 
\begin{align*}
|\nabla\phi_\alpha|\leq&\ C(n,\rv,B)|\nabla\phi_\cN|+|\nabla\psi(b/2^{-\alpha})|,\\
|\Delta\phi_\alpha|\leq&\ C(n,\rv,B)|\Delta\phi_\cN|+|\Delta_\psi(b/2^{-\alpha})|+2|\nabla\phi_\cN||\nabla\psi(b/2^{-\alpha})|.
\end{align*}
By (\ref{neckestimate}), to prove the claim, it suffices to prove that 
\begin{align}
\int_{b\le R}b^{-2}\Big(|\nabla \phi_\cN||\nabla\psi(b/2^{-\alpha})|+\phi_\cN|\Delta\psi(b/2^{-\alpha})|\Big)\,dz+\int_{b\le R}b^{-3}\phi_\cN|\nabla\psi(b/2^{-\alpha})|\,dz\le C(n,\rv,B).
\end{align}
Using 
$\supp\{|\nabla\psi|,|\Delta\psi\}\subseteq \{2^{-\alpha}\geq b\geq 2^{-\alpha-1}\}$ and $|\nabla b|\le C(n,\rv)$ on $\supp\phi_\cN,$ we have

$$\nabla\psi(b/2^{-\alpha})\leq C(n,\rv)2^\alpha$$ on $\supp\nabla\phi_\alpha$.
Moreover since $\Delta b=b^{-1}|\nabla b|^2$, we have
$$\Delta\psi(b/2^{-\alpha})\leq C(n,\rv)2^{2\alpha},$$
on $\supp\nabla\phi_\alpha$.
Noticing that $b^{-1}\approx 2^{\alpha}$ on the support of $\nabla\psi(b/2^{-\alpha}),\Delta\psi(b/2^{-\alpha})$, we have 
$$\int_{b\le R}b^{-2}\phi_\cN|\Delta\psi(b/2^{-\alpha})|+b^{-3}\phi_\cN|\nabla\psi(b/2^{-\alpha})|\,dz\,\le C(n,\rv,B)\int_{2^{-\alpha-1}\le b\le 2^{-\alpha} }b^{-4}\phi_\cN\,dz\,\le C(v,\rv, B)$$
where the last inequality is due to Ahlfors regularity. Using $\nabla\psi(b/2^{-\alpha})\leq C(n,\rv)2^\alpha$ and (\ref{neckestimate}) again, we have  
$$\int_{b\le R}b^{-2}|\nabla \phi_\cN||\nabla\psi(b/2^{-\alpha})|(z)\leq C(n,\rv,N).$$ This finishes the proof of the claim, and hence the corollary.
\end{proof}
\begin{rmk}We remark that the reason why we use function $b$ instead of function $d(\cdot,\cC)$ to define a exhaustion of $\cN$ is that $b$ satisties the nice equation $\Delta b=b^{-1}|\nabla b|^2$. And it is important that these 2 functions can bound each other on the whole neck region. 
\end{rmk}

\subsection{Proof of Theorem \ref{t:neck_region}}\label{ss:neck_region_proof:smooth}

With Theorem \ref{t:L2_neck} in hand, in this subsection we compete the proof of item (4) of Theorem \ref{t:neck_region}.  To show item (4), we only need to show that for any $\epsilon$,  the $\delta<\delta'(n,\rv,\epsilon)$ neck region satisfies
$$\int_{\cN\cap B_1(p)} \brs{H}^4(z)\,dz \leq \epsilon.$$

Recall that, for each $\alpha\in \mathbb N$, we have  defined a sequence of regions $\mathcal N_\alpha\subset \mathcal N$ before Corollary \ref{B} as follows:
$$\mathcal \cN_\alpha=\mathcal \cN\cap \{b>2^{-\alpha}\}$$
we will prove inductively that for any $\alpha\in\mathbb N$, the $L^4$ integral of $|H|$ is $\epsilon$ small on  $\cN_\alpha$.\\
\vspace{.1cm}

\textbf{Theorem $A_\alpha$:}
Let $(M^n,g,H,p)$ satisfy $\Vol(B_1(p))>\rv>0$ 
and $|\RC|\le \delta$.  
For any $\epsilon>0$, if  $\delta<\delta'(n,\epsilon,\rv)$ are such that $\cN = B_2(p)\setminus \bigcup_{x\in \cC} \overline B_{r_x}(x)$ is a $\delta$-neck region. Then we have 	$$\int_{\cN_{\alpha}\cap B_1(p)} \brs{H}^4(z)\,dz \leq \epsilon^4.$$\\
\begin{proof}
We prove Theorem $A_\alpha$ for any $\alpha\in\mathbb N$ inductively. 
When $\alpha=1$, the region $\cN_1$ only has a finite number of scales. Using the fact that function $b$ and $d(\cdot,\mathcal C)$ are comparable, Theorem $A_1$ holds by Lemma \ref{l:neck:scaleinvariantsmall}.  Now we assume $A_{\alpha}$ holds and aim to prove that $A_{\alpha+1}$ holds. Using Corollary \ref{B} with $\epsilon=\epsilon_1$ to be determined, we have 
\begin{equation}\int_{\tilde\cN_{\alpha}\cap B_1(p)} \brs{\Rm}^2(z)\,dz \leq \int_{\tilde\cN_{\alpha}\cap B_1(p)} \brs{H}^4(z)\,dz+\epsilon_1\leq \epsilon^4+\epsilon_1\end{equation}
We remark that since $C(n,\rv)^{-1}\,d(x,\cC)\leq b(x)\leq C(n,\rv)\,d(x,\cC)$, $\cN_\alpha$ and $\tilde\cN_\alpha$ can only differ by finitely many scales. Using this and Lemma \ref{l:neck:scaleinvariantsmall} with $\epsilon=\epsilon_2$, we have  \begin{equation}\int_{\cN_{\alpha}\cap B_1(p)} \brs{\Rm}^2(z)\,dz \leq \int_{\tilde\cN_{\alpha}\cap B_1(p)} \brs{\Rm}^2(z)\,dz+C(n,\rv)\epsilon_2\leq \epsilon^4+ \epsilon_1+C(n,\rv)\epsilon_2.
\end{equation}

Using the scale invariant estimate Lemma \ref{l:neck:scaleinvariantsmall} again, we have  \begin{equation}\int_{\cN_{\alpha+2}\cap B_1(p)} \brs{\Rm}^2(z)\,dz \leq \int_{\cN_{\alpha}\cap B_1(p)} \brs{\Rm}^2(z)\,dz+C(n,\rv)\epsilon_2\leq \epsilon_1+C(n,\rv)\epsilon_2+\epsilon^4.
\end{equation}

Hence \begin{equation}\int_{b\leq R} \phi_{\alpha+1}\brs{\Rm}^2(z)\,dz \leq \epsilon_1+C(n,\rv)\epsilon_2+\epsilon^4.
\end{equation}
Using Corollary $\ref{C}$, with  $\epsilon=\epsilon_3$ to be determined, we have
\begin{align}
\int_{b\leq R} b^{-3} \brs{H}\phi_{\alpha+1}(z)\,dz < \epsilon_3+C(n,\rv)\epsilon^{1/2}\left(\int_{b\le R}|b|^{-3}\brs{H}\phi_{\alpha+1}(z)\,dz\right)^{1/2}\left(\epsilon^4+\epsilon_1+C(r,\rv)\epsilon_2+\epsilon_3\right)^{1/2}\, ,
\end{align}
where we have used the inequality $\brs{H}\le\epsilon C(n,\rv)b^{-1}$ by Lemma \ref{l:neck:Hsmall}. Now let $\epsilon_1=\epsilon^4,\epsilon_2=\epsilon^4,\epsilon_3=\epsilon^4$ (we assume $\epsilon<< C(n,\rv)$), elementary inequalities yield

\begin{equation}\int_{\cN_{\alpha+1}\cap B_{18/10}(p)} b^{-3} \brs{H}\phi_{\mathcal N}(z)\,dz < \epsilon C(n,\rv)^{-3}.\end{equation}
This implies $A_{\alpha+1}$ by using $\brs{H}\le\epsilon C(n,\rv)b^{-1}$ again. Now choosing $\delta\leq \delta'(n,\rv, \epsilon_1,\epsilon_2,\epsilon_3,\epsilon)$ will complete
the inductive argument. 
\end{proof}

\subsection{\texorpdfstring{Proof of the $L^2$ curvature estimate Theorem \ref{t:L2curvature}}{}} 

In this subsection we finish the proof of Theorem \ref{t:L2curvature}.  We first state the neck decomposition theorem which is a slight modification of the neck decomposition Theorem 2.12 of \cite{CJN}. 
\begin{thm} \label{t:neck_decomposition}
		Let $(M^n_i,g_i,H_i, p_i)\to (X,d,p)$ satisfy $\Vol(B_1(p))>\rv>0$ and $\brs{\RC}\le \delta$.  Then for each $0<\delta<\delta(n,v)$ and $\tau(n)$ we can write
	\begin{align}
		&B_1(p)\cap \cR(X)\subseteq \bigcup_a \big(\cN_a\cap B_{r_a}\big) \cup \bigcup_b B_{r_b}(x_b)\, , \notag\\
		&B_1(p)\cap \cS(X) \subseteq \bigcup_a \big(\cC_{0,a}\cap B_{r_a}\big)\cup \tilde \SS(X)\, , \notag
	\end{align}
	where
	\begin{enumerate}
		\item $\cN_a\subseteq B_{2r_a}(x_a)$ are $\delta$-neck regions.
		\item $B_{r_b}(x_b)$ satisfy $r_h(x_b)>2r_b$, where $r_h$ is the harmonic radius.
		\item $\cC_{0,a}\subseteq B_{2r_a}(x_a)\setminus \cN_a$ is the singular set associated to $\cN_a$.
		\item $\tilde\cS(X)$ is a singular set of $n-4$ measure zero.
		\item $\sum_a r_a^{n-4} + \sum_b r_b^{n-4} + H^{n-4}\big(\cS(X)\big) \leq C(n,\rv,\delta)$.
	\end{enumerate}
\end{thm}
\begin{proof} We remark that the Neck Decomposition Theorem  in \cite{CJN} is proved under the condition that $\Ric^g\geq-\delta$ rather than that $\brs{\RC}\le \delta$. Also in their definition of neck region, they use the ball is Gromov Hausdorff close to $\mathbb R^{n-4}\times C(Y)$, where $Y$ is a smooth metric space.  We point out that under the conditions $\brs{\RC}\le \delta$ and $\Vol(B_1(p))>\rv>0$, when a unit ball is sufficiently close to $\mathbb R^{n-4}\times C(Y)$, $C(Y)$ is automatically close to $C(S^3/\Gamma)$ by the cone rigidity Theorem \ref{c:rigidity}. With this remark, it is clear that the Neck Decomposition stated above holds. \end{proof}
With the neck decomposition theorem above and the $L^2$ curvature estimate Theorem \ref{t:L2_neck} on the neck region in hand, the proof of Theorem \ref{t:L2curvature} follows the argument in \cite{JiNa}. We include the argument for reader's convenience. 

\begin{proof}[Proof of Theorem \ref{t:L2curvature}]
 Indeed, for $\delta>0$ let us use Theorem \ref{t:neck_decomposition} and consider the covering
\begin{align}
	B_1(p)\subseteq \bigcup_a \big(\cN_a\cap B_{r_a}\big) \cup \bigcup_b B_{r_b}(x_b)\, ,
\end{align}
where $\cN_a\subseteq B_{2r_a}(x_a)$ is a $(\delta,\tau)$-neck, $r_h(x_b)\geq 2r_b$, and
\begin{align}\label{e:L2proof:content}
	\sum_a r_a^{n-4} + \sum_b r_b^{n-4}< C(n,\rv,\delta)\, .
\end{align}

Now since the generalized Ricci curvature is uniformly bounded and we have a harmonic radius lower bound, by Lemmas \ref{H-control-regular} and \ref{H-control}, we have scale invariant $W^{2,2}$-estimates on the metric on $B_{r_b}(x_b)$.  In particular, we have the estimate
\begin{align}
r^{4-n}_b \int_{B_{r_b}(x_b)} \brs{\Rm}^2(z)\,dz < C(n)\, .	
\end{align}

On the other hand, if we let $\epsilon=1$ and fix a $\delta=\delta(n,\rv,\epsilon=1)$, then by Theorem \ref{t:neck_region} we have the scale invariant $L^2$ estimate on each neck region $\cN_a$ given by:
\begin{align}
r^{4-n}_a \int_{\cN_a\cap B_{r_a}} \brs{\Rm}^2(z)\,dz < C(n)\, .	
\end{align}

Combining these estimates with the content estimate of \eqref{e:L2proof:content} we get the estimate
\begin{align}
\int_{B_1(p)} \brs{\Rm}^2(z)\,dz &\leq \sum_a \int_{\cN_a\cap B_{r_a}} \brs{\Rm}^2(z)\,dz + \sum_b \int_{B_{r_b}(x_b)}	\brs{\Rm}^2(z)\,dz\, \notag\\
&\leq C(n)\Big( \sum_a r_a^{n-4} + \sum_b r_b^{n-4}\Big) \leq C(n,\rv)\, ,
\end{align}
which completes the proof of item (1) of Theorem \ref{t:L2curvature}. Item (2) of Theorem \ref{t:L2curvature} is a direct consequence of item (3).  Now item (3) can be derived by using the Bochner formula \ref{f:Hbochner} and integration by parts.
\end{proof}

\section{Applications}
\subsection{Finite diffeomorphism type of manifold with bounded generalized Ricci tensor}

\begin{proof}[Proof of Corollary \ref{t:main_finite_diffeo}:] The finiteness of diffeomorphism types of manifold $(M,g)$ satisfying 
$|\Ric^g|\leq3, \Vol(B_1)>\rv, \diam(M)<D$ is proved by Cheeger and the second author in \cite{Codim4} as an application of codimension 4 regularity. We now state a theorem generalizing Theorem  8.6 of \cite{Codim4}. 

\begin{thm}\label{t:annulus2}
 For every $\epsilon>0$ there exists $\delta(\rv,\epsilon)>0$ 
such that if $M^4$ satisfies $|\RC_{M^4}|\leq 3\delta$, $\Vol(B_1(p))>\rv>0$ and $|\cV^\delta_{4}(p)-\cV^\delta_{1/4}(p)|<\delta$, then there exists a 
discrete subgroup $\Gamma\subseteq {\rm O}(4)$ with $|\Gamma|\leq N(\rv)$ such that the following hold:
\begin{enumerate}
\item For each $x\in A_{\epsilon,2}(p)$ we have the harmonic radius lower bound $r_h(x)>r_0(\rv)\epsilon$.
\vskip1mm

\item There exists a subset $A_{\epsilon,2}(p)\subseteq U\subseteq A_{\epsilon/2,2+\epsilon}(p)$ and a 
diffeomorphism $\Phi:A_{\epsilon,2}(0)\to U$, with $0\in \dR^4/\Gamma$, such that
if $g_{ij}=\Phi^*g$ is the pullback metric then
\begin{align}
||g_{ij}-\delta_{ij}||_{C^{0}} + ||\partial_k g_{ij}||_{C^{0}}<\epsilon\, .
\end{align}
\end{enumerate}
\end{thm}
\begin{proof}
 The proof is by contradiction.  So let us assume for some $\epsilon>0$ there is no such $\delta(\rv,\epsilon)>0$.  
Thus, we have a sequence of spaces $(M^4_j,g_j,p_j)$ with $\Vol(B_1(p_j))>\rv>0$, $|\RC_{M^4_j}|\leq \delta_j\to 0$ and 
$|\cV_{4}(p_j)-\cV_{1/4}(p_j)|<\delta_j\to 0$, but the conclusions of the theorem fail.  After passing to a subsequence 
we can take a limit
\begin{align}
(M^4_j,d_j,p_j)\stackrel{d_{GH}}{\longrightarrow}(X,d,p)\, .
\end{align}
Using the almost volume cone implies almost metric cone Theorem of \cite{ChC1}, we then have 
\begin{align}
B_4(p) = B_4\big(y_0)\, ,
\end{align}
where $y_0\in C(Y)$ is the cone vertex and $Y$ some metric space of diameter $\leq\pi$.  

Now using Theorem \ref{t:codim4}, we know that away from a set of codimension $4$ in $C(Y)$,  the 
harmonic radius $r_h>0$ is bounded uniformly from below.  Assume there is some point $y\in Y$ such that $r_h(y)=0$ 
and consider the ray $\gamma_y$ in $C(Y)$ through the point $y$. In that case, it would follow that for
 every point of $\gamma_y$, the harmonic radius $r_h=0$ vanishes.  The ray $\gamma$ 
has Hausdorff dimension $1$, and therefore its
existence would contradict Theorem \ref{t:codim4}.  
Thus, we conclude that $r_h>0$ and that $Y=(Y,g_Y)$ 
is a $C^{1,\alpha}\cap W^{2,q}$ manifold for every $\alpha<1$ and $q<\infty$.  

Now by writing the formula for the Ricci tensor in harmonic coordinates and using
$|\RC_{M^4_j}|\to 0$, it follows that $C(Y)$ is smooth and Ricci flat away from the vertex.  In particular, since $C(Y)$
 is a metric cone over $Y$, we must have $\Ric^g_{Y^3} = 3g^Y$.  Since in dimension $3$, constant Ricci curvature
 implies constant sectional curvature, it follows  $Y=S^3/\Gamma$ has constant sectional curvature $\equiv 1$.
  Additionally, we know from the volume bound, $\Vol(B_1(p))>\rv>0$, 
that the order $|\Gamma|<N(\rv)$ is uniformly bounded.  
In particular, we have that $C(Y)=\dR^4/\Gamma$ is an orbifold with an isolated singularity. 

 It now follows that there exists 
$r_0(\rv)>0$ such that for $y\in \dR^4/\Gamma$ with $|y|=1$, we have
\begin{align}
B_{2r_0}(y) = B_{2r_0}(0^4)\, ,
\end{align}
where $0^4\in \dR^4$.  In particular, for all $j$ sufficiently large, we have from the $\epsilon$-regularity theorem, 
Theorem \ref{s:eps_reg}, that for all $x\in A_{\epsilon,2}(p_j)$, the harmonic radius,
$r_h(x)>r_0(\rv,\epsilon)=r_0(\rv)\epsilon$ is bounded uniformly from below independent of $j$.  
Thus, if there exists $\epsilon$ as above, for which there is no $\delta(\rv,\epsilon)$, it must be (2) that fails to hold.

However, by  using the first part of  Theorem 8.1 of \cite{Codim4}  ,
we have that for $j$ sufficiently large, there exist diffeomorphisms
\begin{align}
\Phi_j:A_{\epsilon,2}(0)\to M^4_j\, ,
\end{align}
such that
\begin{align}
\Phi_j^*g_j\stackrel{C^{0}} {\longrightarrow} dr^2+r^2g_Y\, 
\end{align}
on $A_{\epsilon,2}(0)$. Moreover, we have $\brs{\RC}<n-1$ and $r_h(x)>r_0(\rv)\epsilon$, thus using Lemma \ref{H-control-regular} and Lemma \ref{H-control}, we have
\begin{align}
\Phi_j^*g_j\stackrel{C^{1,\alpha}\cap W^{2,q}}{\longrightarrow} dr^2+r^2g_Y\, .
\end{align}
For $j$ sufficiently large, this implies that (2) holds; a contradiction.
\end{proof}
\noindent Using Theorem \ref{t:annulus2}, the proof of Corollary \ref{t:main_finite_diffeo} will be verbatim as in \cite{Codim4}.
\end{proof}

\subsection{Proof of Corollary \ref{c:rigidity}}

\begin{proof}[Proof of Corollary \ref{c:rigidity}] 
As we have assumed $H \in \Lambda^3$, consider the one-form $\theta = \star H$.  By the generalized Einstein equation, $\theta$ is harmonic for the Hodge Laplacian, so by the Bochner formula for one-forms, the generalized Einstein equation, and the fact that $\theta^{\sharp} \hook H = 0$ we obtain
\begin{align*}
    \gD \brs{\theta}^2 = \brs{\N \theta}^2 + \Rc^g(\theta^{\sharp}, \theta^{\sharp}) = \brs{\N \theta}^2 + \tfrac{1}{4} H^2(\theta^{\sharp}, \theta^{\sharp}) = \brs{\N \theta}^2 + \tfrac{1}{4} \brs{\theta^{\sharp} \hook H}^2 = \brs{\N \theta}^2\ .
\end{align*}
Fix $p \in M$, $R > 0$, and let $\eta$ denote a cutoff function for $B_R(p)$.  Using the formula above we obtain after integration by parts
\begin{align*}
    0 =&\ \int_M \left[ - \gD \brs{\theta}^2 + \brs{\N \theta}^2 \right] \eta^2\\
    =&\ \int_M \brs{\N \theta}^2 \eta^2 + \theta \star \N \theta \star \eta \N \eta\ .
\end{align*}
We can furthermore estimate
\begin{align*}
    \brs{\int_M \theta \star \N \theta \star \eta \N \eta} \leq&\ \brs{\brs{\eta \brs{\N \theta}}}_{L^2(M)} \brs{\brs{\theta}}_{L^4(\supp(\eta))} \brs{\brs{ \brs{\N \eta}}}_{L^4(\supp(\eta))}\\
    \leq&\ \tfrac{1}{2} \int_M \brs{\N \theta}^2 \eta^2 + C \brs{\brs{ \theta}}_{L^4(\supp(\eta))}^2 \brs{\brs{\N \eta}}_{L^4(\supp(\eta))}^2\ .
\end{align*}
Since $g$ has nonnegative Ricci curvature, we know that $\Vol(B_R(p)) \leq C R^4$, thus since $\brs{\N \eta} \leq C R^{-1},$ it follows that the $L^4$ norm of $\brs{\N \eta}$ is uniformly bounded for an arbitrary choice of $R$.  Furthermore, by the Euclidean volume growth hypothesis, we can rescale arbitrarily large balls to unit size and apply Theorem \ref{t:L2curvature} to conclude that $|H|$ is bounded in $L^4$, thus $\theta$ is bounded in $L^4$, it follows then that 
\begin{align*}
    \lim_{R \to \infty} C \brs{\brs{ \theta}}_{L^4(\supp(\N \eta))}^2 \brs{\brs{\N \eta}}_{L^4(\supp(\N \eta))}^2 = 0\ .
\end{align*}
It follows that $\N \theta \equiv 0$.  Again since $\theta$ is in $L^4$, it follows that $\theta = 0$, hence $H = 0$ as claimed.
\end{proof}

\bibliographystyle{plain}

\end{document}